\documentclass[12pt,reqno]{article}

\usepackage[usenames]{color}
\usepackage{amssymb}
\usepackage{amsmath}
\usepackage{amsthm}
\usepackage{amsfonts}
\usepackage{amscd}
\usepackage{enumerate}

\usepackage[colorlinks=true,
linkcolor=webgreen,
filecolor=webbrown,
citecolor=webgreen]{hyperref}

\definecolor{webgreen}{rgb}{0,.5,0}
\definecolor{webbrown}{rgb}{.6,0,0}

\usepackage{color}
\usepackage{fullpage}
\usepackage{float}
\usepackage{maple}

\usepackage{graphics}
\usepackage{latexsym}
\usepackage{epsf}
\usepackage{breakurl}

\numberwithin{equation}{section}

\newcommand\numberthis{\addtocounter{equation}{1}\tag{\theequation}}

\begin{document}
\parskip 8pt
\parindent 20pt


\theoremstyle{plain}
\newtheorem{theorem}{Theorem}
\newtheorem{cor}[theorem]{Corollary}
\newtheorem{lemma}[theorem]{Lemma}
\newtheorem{prop}[theorem]{Proposition}

\theoremstyle{definition}
\newtheorem{definition}[theorem]{Definition}
\newtheorem{example}[theorem]{Example}
\newtheorem{conjecture}[theorem]{Conjecture}

\theoremstyle{remark}
\newtheorem{rem}[theorem]{Remark}

\begin{center}
\vskip 1cm{Five Parameter Hypergeometric 3F2(1) when One or more Parameters are Integers or Separated by Integers:  Derivations, Review, Exotics and More
}
\vskip 1cm
\large
 
Michael Milgram, Consulting Physicist, Geometrics Unlimited, Ltd.\\Box 1484, Deep River, Ont., Canada.\\

\href{mailto:email}{\tt  mike@geometrics-unlimited.com} \\
September 30, 2025\\
\vskip .2 in
\end{center}
MSC Categories: 33C20, 33c90, 33-01\\\\
Keywords: Hypergeometric, 3F2(1), Saalsch\"{u}tz, exotic 
\vskip .2 in

{\bf Revisons: November 18,2025}.
\newline\\
1. Minor corrections throughout;\\
2. Major revision to Appendix \ref{sec:Thomae} and discussion about terminating series added throughout;\\
3. Addition of \eqref{T1Cx};\\
4. Addition of Sections (\ref{sec:ApplContig}), (\ref{sec:CWExotic})  and (\ref{sec:XWang});\\
5. Addition of several new references;

\begin{abstract}
 This work was intended to be all about, and only about, hypergeometric 3F2(1). The initial goal was to revisit many identities from the literature that have been derived over the years and show that they can be obtained in a simpler way armed, with only a minimum of elementary identities. That goal has been achieved as a (patient) reader will discover. In another sense, this work is a partial review of the last half-century's worth of progress in the evaluation of a particular set of 3F2(1), in particular those cases where at least one parameter is an integer or two or more parameters are separated by an integer. The result is a collection of very general identities (or techniques) that an analyst seeking to evaluate a particular 3F2(1), might want to consider as a starting point. That is the secondary goal.

Along the way however, the temptation arose to investigate at least one of the unanswered questions that others have raised. This led to a few digressions, and some possibly new results. The reader is invited to follow where curiosity led me to depart from a straightforward review of the state-of-the-art.
\end{abstract}

\section{Introduction and Background} \label{I&B}

The infinite series defined by
\begin{equation}
_{3}F_{2}(a,b,c;e,f;1)\equiv\moverset{\infty}{\munderset{k =0}{\sum}}\frac{\Gamma \left(k +a \right) \Gamma \left(k +b \right) \Gamma \left(k +c \right) \Gamma \left(e \right) \Gamma \left(f \right)}{\Gamma \left(a \right) \Gamma \left(b \right) \Gamma \left(c \right) \Gamma \left(k +1\right) \Gamma \left(k +e \right) \Gamma \left(k +f \right)}
\label{3F2Def}
\end{equation}
arises frequently in mathematical physics and analysts commonly seek a closed form evaluation for special combinations and values of the parameters $\{a,b,c,e,f\}$. In particular, the combination of parameters  $s\equiv\Re(e+f-a-b-c)$, at times referred to as the ``parametric excess", defines the {\bf s-balanced} sum that converges in general only if $\Re(s)>0$ or one or more of the parameters $a,b,c$ is a non-positive integer. This sum is of special interest when $s=n$, a positive integer.

\begin{rem} {\bf Regarding notation:} The above is the notation that will be used here with the caveat that if the referenced inequality fails, any corresponding summable or closed form evaluation of the right-hand side of \eqref{3F2Def} represents the left-hand side by the principle of analytic continuation, (with caveats -- see the discussion preceding \eqref{T2a} below). 

Since an alternate typographical notation exists, that is $_{3}F_{2}{(\overset {a,b,c}{\underset {e,f}{\null}}}|1)$, which will be used infrequently here,  any of the corresponding parameters may be labelled {\it top} or {\it bottom}. {\bf Throughout, the symbols $j,k,m,n,p,q$ are non-zero positive integers} unless specified otherwise, and, for brevity, the left-hand side of \eqref{3F2Def} may be written $3F2(1)$ (or $pFq(1)$ in general). All other symbols are implicitly complex, but are explicitly treated as real unless required by the context. {\bf Any summation vanishes if the lower limit exceeds the upper limit}; the symbol $:=$ indicates symbolic replacement and the symbol $a \leftrightarrow b$ means symbolic interchange of symbols $a$ and $b$. The word ``summable" means that a sum is terminating or can be transformed into a terminating sum; the word ``closed" means that an expression has been reduced to a form that no longer includes a summation operation (in spite of the fact that the digamma function $\psi(x)$ itself, by definition, usually includes a buried summation operation).
The $n^{th}$ derivative of digamma functions is indicated by $\psi(n,x)$, $G$ is Catalan's constant and the Pochhammer symbol $(x)_{n}\equiv\Gamma(x+n)/\Gamma(x)$. Many of the identities derived here are based on the Gamma function reflection formula $\Gamma(1-s)=\pi/\sin(\pi s)/\Gamma(s)$; in particular for the case of integer arguments, $\Gamma(k-n)/\Gamma(-n)=(-1)^k\Gamma(1+n)/\Gamma(1-k+n)$. 
\end{rem} 
The earliest 3F2(1) evaluations, due to Dixon, Whipple and Watson (DWW), reproduced in Appendix \ref{sec:Contig} when $m=n=0$, were summarized in 1935 by Bailey \cite{Bailey1935}, who also presented the 1879 Thomae two-part transformations plus Whipple's analysis of three-part transformations. For a listing of the Thomae transformations as well as Hardy's ingenious derivation, see Appendix \ref{sec:Thomae}. Using the Thomae transformations, any of the DWW identities (and only one of the DWW identities) can be obtained starting from any one of the others. Over the years, numerous efforts were made  to extend the DWW identities to contiguous values of their parameters mostly by recursion (see \cite{Awadetal2025} for references), until finally in 2012, Chu \cite{ChuW} resolved that issue by giving a complete summary of all possible transformations contiguous to the basic DWW identities, thus rendering these identities summable (for notation, see Appendix \ref{sec:Contig}; for examples, see \eqref{AOT} and \eqref{AOTP} below). Notwithstanding Chu's encompassing solution, a recent paper by Awad, Rakha and Mohammed \cite{Awadetal2025} has revisited the older recursion methods.   

With respect to the Thomae transformations, it was commonly believed that these 120 transformations, of which only 9 are independent, represent the only general five-parameter, two-part transformations between  3F2(1), although a few curious exceptions were known. However, in 2006, Krattenthaler and Rivoal \cite{KrattRiv} proved the existence of many exceptions to this belief, labelling them ``exotic". Shortly  thereafter, Chu and Wang \cite{ChuWang2007} elaborated on these identities, but, in both cases, the proofs were limited to those 3F2(1) that involved only non-negative integer parameters. In Section \ref{sec:Exotic1}, examples of exotic terminating cases involving non-integer parameters will be studied.

In addition to the three well-known DWW identities, numerous collections of special cases with varying numbers of independent parameters have appeared in the literature over the years; perhaps the general summaries of Prudnikov et. al. \cite[Section 7.4.4]{prudnikov} and Slater \cite[Appendix III]{Slater} are the best known, recognizing that collections of very specialized evaluations also exist (e.g. \cite{KrupKol} where Krupnikov and K{\"o}lbig list specialized cases that emphasize fractional parameters). Adding to the collections of listed identities with specialized parameters, Rainville \cite[pages 81 ff]{Rainville} provides four master contiguity relations among instances of 3F2(1) with arbitrary parameters; others have added significantly to this list (e.g. \eqref{HidbKR} below), particularly Chu and Wang \cite{ChuWang2007}. A few of these are collected in Appendix \ref{sec:select}. 

To summarize, in this work, I attempt to deal mostly with cases involving five independent parameters, where the parameters are either integers or separated by integers. So, in Section \ref{sec:History} are presented the classical identities that transform an infinite into a terminating sum, where it is shown that it is possible to derive (and extend) each of them in a simpler way, starting only from the basic definition \eqref{3F2Def} and the Thomae identities. In Section \ref{sec:Hunt}, I attempt to utilize the methods developed in the previous  section to confront a special case challenge posed elsewhere, with limited success. However, by reducing the challenge case to a special subset, a number of new general identities are found and simpler derivations of known identities are developed. These identities take us far afield, involving some 4F3(1), a special Euler sum and a generalization of the well-known Saalsch\"{u}tz theorem \cite{Bailey1935}. In Section \ref{sec:Exotic1}, the transformations that have been deduced throughout are examined from the point of view of their ``exoticness". These identities are then employed in Section \ref{sec:LimSpec} to obtain general transformations, as well as special cases that are either new, or more easily derived than previously thought. Finally, in Section \ref{sec:Applics}, I focus on a selection of identities extracted from the recent literature to show how it is possible to use the tools developed here to obtain or generalize given identities in a simpler way.

This summary now ends with a collection identifying some of the more interesting identities developed throughout, labelled by their equation reference:\\
\\
\eqref{MintonCor}:  a $p-$balanced extension of Minton's terminating sum, where a top parameter exceeds a bottom parameter by an integer:
\begin{equation} \nonumber
{}_{3}^{}{\moversetsp{}{\mundersetsp{}{F_{2}^{}}}}\left(-m ,a ,c +n ;c ,p +n +a -m ;1\right);
\end{equation}
\eqref{S1b}: an m-balanced variation of the KRGT identities studied in Section \ref{sec:History};
\begin{equation} \nonumber
{}_{3}^{}{\moversetsp{}{\mundersetsp{}{F_{2}^{}}}}\left(a ,b ,c +n ;c ,a +b+n +m ;1\right)
\end{equation}
\eqref{Hnew2}: an $m-$balanced, non-terminating analogue of the Saalsch\"{u}tz theorem (c.f. \eqref{Saal}):
\begin{equation} \nonumber
{}_{3}^{}{\moversetsp{}{\mundersetsp{}{F_{2}^{}}}}\left(a ,b ,n ;c ,a+b -c  +n +m ;1\right);
\end{equation}
\eqref{Cm12}: a contiguity relation for the challenge problem:
\begin{align} \nonumber
{}_{3}^{}{\moversetsp{}{\mundersetsp{}{F_{2}^{}}}}&\left(a ,b ,c ;a +1,e;1\right)
+A\, {}_{3}^{}{\moversetsp{}{\mundersetsp{}{F_{2}^{}}}}\left(a ,b ,c ; a +2,e;1\right)+B=0\,.
\end{align}
\eqref{Cx1(4)}: a non-terminating sum with only positive integer parameters:
\begin{equation} \nonumber
{}_{3}^{}{\moversetsp{}{\mundersetsp{}{F_{2}^{}}}}\left(j ,m ,n ;p ,1;1\right);
\end{equation}
\eqref{HidZ}: a zero-balanced, terminating variant of the Saalsch\"{u}tz theorem (c.f. \eqref{Saal}).
\begin{equation} \nonumber
{}_{3}^{}{\moversetsp{}{\mundersetsp{}{F_{2}^{}}}}\left(1-n,a ,m  ;1,a -n +m ;1\right);
\end{equation}
\eqref{Xidb}: an $n-$balanced, (terminating) generalization of the Saalsch{\"u}tz theorem (c.f. \eqref{Saal}):
\begin{equation} \nonumber
{}_{3}^{}{\moversetsp{}{\mundersetsp{}{F_{2}^{}}}}\left(-m,a ,b  ;c ,a+b -c+n-m   ;1\right);
\end{equation}
\eqref{Excb4}: a summable form frequently encountered: 
\begin{equation} \nonumber
{}_{3}^{}{\moversetsp{}{\mundersetsp{}{F_{2}^{}}}}\left(a,b ,m  ;1,a +n;1\right);
\end{equation}
\eqref{B3e0}: an Euler-type sum:
\begin{equation}
\moverset{m}{\munderset{k =0}{\sum}}\left(-1\right)^{k}\,{\binom{m}{k}}\,\psi \left(k +1\right)
\end{equation}

\section{A historical perspective } \label{sec:History}
\subsection{The familiar identities} \label{sec:familiar}
An important identity is characterized solely by the presence of a negative integer top parameter, leading to the elementary identification of a terminating sum where convergence is not an issue:
\begin{equation}
{}_{3}^{}{\moversetsp{}{\mundersetsp{}{F_{2}^{}}}}\left(-n,a ,b  ;c ,e ;1\right)
 = 
\frac{\Gamma \left(c \right)\Gamma \left(e \right)\,\Gamma \left(n +1\right)}{\Gamma \left(a \right) \Gamma \left(b \right)}\moverset{n}{\munderset{k =0}{\sum}}\frac{\Gamma \left(a +k \right) \Gamma \left(b +k \right)   \left(-1\right)^{k}}{ \Gamma \left(c +k \right) \Gamma \left(e +k \right) \Gamma \left(k +1\right) \Gamma \left(n -k +1\right)}.
\label{NegTop}
\end{equation}
A similar special case generally attributed  to Sheppard \cite[Corollary 3.3.4]{AndAskRoy}, involves the transformation
\begin{align} \nonumber
{}_{3}^{}{\moversetsp{}{\mundersetsp{}{F_{2}^{}}}}\left(-n ,a ,b ;c ,e ;1\right)
 &= 
\frac{\Gamma \left(n +c -a \right) \Gamma \left(n +e -a \right) \Gamma \left(c \right) \Gamma \left(e \right) }{\Gamma \left(c -a \right) \Gamma \left(e -a \right) \Gamma \left(c +n \right) \Gamma \left(n +e \right)}
\\ &
\times{}_{3}^{}{\moversetsp{}{\mundersetsp{}{F_{2}^{}}}}\left(-n ,a ,1-e -c +a +b -n ;1+a -c -n ,1+a -e -n ;1\right),
\label{Sheppard}
\end{align}
which, if we set $e=a+b-n+1-c$, reduces to the closed, (one-balanced) Saalsch\"{u}tz theorem \cite[page 9]{Bailey1935}
\begin{align} \nonumber
&\hspace{-40pt}{}_{3}^{}{\moversetsp{}{\mundersetsp{}{F_{2}^{}}}}\left(-n,a ,b ;c ,1+a +b -c -n ;1\right)
 \\ 
&=\frac{\Gamma \left(1+b -c \right) \Gamma \left(1+a -c \right) \Gamma \left(1+a +b -c -n \right) \Gamma \left(1-c -n \right)}{\Gamma \left(1+b -c -n \right) \Gamma \left(1+a -c -n \right) \Gamma \left(1-c +b +a \right) \Gamma \left(1-c\right)}\,.
\label{Saal}
\end{align}
For a derivation of \eqref{Sheppard} see \eqref{Eq16}ff, and for a generalization of \eqref{Saal} see \eqref{R3}, both below. Among the other few closed identities known, we find the important terminating, zero-balanced,  special case that arises frequently, merits individual attention and is characterized by a top parameter exceeding a bottom parameter by a positive integer. It is due to Minton \cite[Eq. (10)]{Minton}, and, in a generalized form, reads
\begin{align} \nonumber
{}_{3}^{}{\moversetsp{}{\mundersetsp{}{F_{2}^{}}}}&\left(-m ,a ,c +n ;c ,n +a -m ;1\right)
 = 
\left(-1\right)^{m}\,\frac{\Gamma \left(c \right) \Gamma \left(n +a -m \right) \Gamma \left(1+m \right)}{\Gamma \left(c +n \right) \Gamma \left(a \right)}\\
&+\frac{\Gamma \left(c \right) \Gamma \left(n +a -m \right) \Gamma \left(1+n \right) }{\Gamma \left(c +n \right) \Gamma \left(a -m -c \right)}\moverset{n}{\munderset{k =1+m}{\sum}}\frac{\Gamma \left(-c +a -m +k \right) \left(-1\right)^{k}}{\Gamma \left(1+n -k \right) k\,\Gamma \left(a +k \right) \Gamma \left(k-m  \right)}\,.
\label{Minton}
\end{align}
\begin{rem}
For a proof of \eqref{Minton} see Appendix \ref{sec:Proofs}, Theorem (\ref{sec:MintsExTheorem}). Minton's original derivation (for the zero balanced sum), only considered the case $m\geq n$ where the sum in \eqref{Minton} vanishes. For an extension of \eqref{Minton} that is $p-$balanced, see \eqref{MintonCor} below.
\end{rem}

In the case of a non-terminating sum where a top parameter exceeds a bottom parameter by a positive integer, Karlsson \cite[Eq. (10)]{Karlsson}, using contour integration, 
obtained the transformation
\begin{equation}
{}_{3}^{}{\moversetsp{}{\mundersetsp{}{F_{2}^{}}}}\left(a ,b ,c +n ;c ,e ;1\right)
 = 
\frac{\Gamma \left(e \right) \Gamma \left(e -a -b \right)}{\Gamma \left(e -b \right) \Gamma \left(e -a \right)}\; {}_{3}^{}{\moversetsp{}{\mundersetsp{}{F_{2}^{}}}}\left(-n,a ,b  ;c ,1-e +a +b ;1\right).
\label{Karl8}
\end{equation}
Here, the right-hand side terminates by \eqref{NegTop} and therefore \eqref{Karl8} sums a specialized infinite series where $e-a-b-n>0$ and analytically continues the left-hand side if the inequality fails. By setting $e=b+1$, the right-hand side of \eqref{Karl8} reduces to a (summable) 2F1(1) with $1-a-n>0$, yielding a closed-form generalization of \eqref{Minton}:
\begin{equation}
{}_{3}^{}{\moversetsp{}{\mundersetsp{}{F_{2}^{}}}}\left(a ,b ,c +n ;c ,b +1;1\right)
 = 
\frac{\Gamma \left(b +1\right) \Gamma \left(1-a \right) \Gamma \left(c +n -b \right) \Gamma \left(c \right)}{\Gamma \left(b +1-a \right) \Gamma \left(c +n \right) \Gamma \left(c -b \right)}.
\label{K10}
\end{equation}

As Karlsson also points out \cite[Eq. (11)]{Karlsson}, a special case mixture of \eqref{Minton} and \eqref{Karl8} is 
\begin{equation}
{}_{3}^{}{\moversetsp{}{\mundersetsp{}{F_{2}^{}}}}\left(-(m +n),a +m ,b +n ;a ,b ;1\right)
 = 
\frac{\left(-1\right)^{m +n}\,\Gamma \left(m +n +1\right) \Gamma \left(a \right) \Gamma \left(b \right)}{\Gamma \left(a +m \right) \Gamma \left(b +n \right)}\,.
\label{Karl11}
\end{equation}

In 1981, by reducing $q-$analogues, Gasper \cite{Gasper_1981}, as rederived\footnote{For two alternative derivations of \eqref{Gasper} see the commentary following \eqref{GaRosen} and \eqref{Thm1Rev}.} by Karp and Prilepkina  \cite[Eq. (2.15)]{Karp&Prep2018}, found an alternative representation of \eqref{Karl8}, with the transformation 
\begin{align} \nonumber
{}_{3}^{}{\moversetsp{}{\mundersetsp{}{F_{2}^{}}}}\left(a ,b ,c +n ;c,e ;1\right)
 &= 
\frac{\Gamma \left(e \right) \Gamma \left(1-a \right) \Gamma \left(n +c -b \right) \Gamma \left(c \right) }{\Gamma \left(c -b \right) \Gamma \left(b +1-a \right) \Gamma \left(e -b \right) \Gamma \left(c +n \right)}\\&
\times{}_{3}^{}{\moversetsp{}{\mundersetsp{}{F_{2}^{}}}}\left(b ,1 +b-c,1+b -e ;b +1-a ,1+b -c -n ;1\right)\,.
\label{Gasper}
\end{align}
If we now choose $e=b+1+m$, it is clear that the right-hand side of \eqref{Gasper} terminates (see \eqref{NegTop}), yielding the specialization
\begin{align} \nonumber
{}_{3}^{}{\moversetsp{}{\mundersetsp{}{F_{2}^{}}}}&\left(a ,b ,c +n ;c,b +1+m  ;1\right)
 = 
\frac{\Gamma \left(b +1+m \right) \Gamma \left(1-a \right) \left(-1\right)^{n}\,\Gamma \left(c \right) }{\Gamma \left(c +n \right) \Gamma \left(b \right)}\\
&\times\moverset{m}{\munderset{k =0}{\sum}}\frac{\left(-1\right)^{k}\,\Gamma \left(b +k \right) \Gamma \left(1+b -c +k \right)}{\Gamma \left(1-k +m \right) \Gamma \left(k +1\right) \Gamma \left(b +1-a +k \right) \Gamma \left(1+b -c -n +k \right)},
\label{G81M}
\end{align}
reducing to \eqref{K10} if $m=0$. Also see \eqref{KarpKalmy} below.

Somewhat later, defining ${\bf c}\equiv \{c_{1}\dots c_{r}\}$, Rosengren \cite{Rosengren2004}, working within a multvariable framework, obtained an identity for ${}_{r+2}F_{r+1}(a,b,{\bf c+n};e,{\bf c};1)$, which, specialized to the case $r=1$ yields \cite[Eq. (2)]{Rosengren2004}
\begin{align} \nonumber
{}_{3}{\moversetsp{}{\mundersetsp{}{F_{2}}}}\left(a ,b ,c +n ;c,e ;1\right) \nonumber
&=\frac{\Gamma \left(e \right) \Gamma \left(c +1-e +n \right) \Gamma \left(c \right) \Gamma \left(e -n -a -b \right) \Gamma \left(1-a \right)}{\Gamma \left(e -a \right) \Gamma \left(c +n \right) \Gamma \left(e -b \right) \Gamma \left(1-n -a \right) \Gamma \left(c +1-e \right)}\\ 
&\times {}_{3}^{}{\moversetsp{}{\mundersetsp{}{F_{2}^{}}}}\left(-n ,c-a  ,1 +b-e ;1+c -e ,1-n -a ;1\right),
\label{RosenG}
\end{align}
a terminating sum to be compared with \eqref{Karl8} and \eqref{Gasper} -- see Section \ref{sec:versus}. In summation form, Rosengren's terminating transformation reads
\begin{align} \nonumber
{}_{3}{\moversetsp{}{\mundersetsp{}{F_{2}}}}&\left(a ,b ,c +n ;c,e ;1\right)\\ \nonumber
&= 
\frac{\Gamma \left(e \right) \Gamma \left(c +1-e +n \right) \Gamma \left(c \right) \Gamma \left(1+n \right) \Gamma \left(e -n -a -b \right) \Gamma \left(1-a \right) \sin \left(\pi  \left(e -b \right)\right) }{\Gamma \left(e -a \right) \Gamma \left(c +n \right) \pi \,\Gamma \left(c -a \right)}\\
&\times \moverset{n}{\munderset{k =0}{\sum}}\frac{\Gamma \left(b +1-e +k \right) \Gamma \left(c -a +k \right) \left(-1\right)^{k}}{\Gamma \left(1-n -a +k \right) \Gamma \left(c +1-e +k \right) \Gamma \left(1+n -k \right) \Gamma \left(k +1\right)}\,,
\label{RosenGSum}
\end{align}
and, if $e=b+m+1$, it is straightforward to find, using the Gamma function reflection property, that \eqref{RosenGSum} reduces to
\begin{align} \nonumber
{}_{3}^{}{\moversetsp{}{\mundersetsp{}{F_{2}^{}}}}&\left(a ,b ,c +n ;c ,b +m +1;1\right)\\ \nonumber
& = 
\frac{\Gamma \left(b +m +1\right) \Gamma \left(n-b -m +c  \right) \Gamma \left(c \right) \Gamma \left(m +1-a -n \right) \Gamma \left(1-a \right) \Gamma \left(1+n \right)}{\Gamma \left(c-a \right) \Gamma \left(b +m -a +1\right) \Gamma \left(c +n \right)}\\
& \times\moverset{\min \left(m , n\right)}{\munderset{k =0}{\sum}}\frac{\Gamma \left(c -a +k \right)}{\Gamma \left(1-k +n \right) \Gamma \left(1-k+m \right) \Gamma \left(k +1\right) \Gamma \left(c -b -m +k \right) \Gamma \left(1-a -n +k \right)}\,.
\label{RosengEb}
\end{align}

  
In the case of terminating series, two of the Thomae transformations generate a transformation between terminating series. From Thomae transform \eqref{Thom8} with $c=-n$, we have
\begin{equation}
{}_{3}^{}{\moversetsp{}{\mundersetsp{}{F_{2}^{}}}}\left(-n ,a ,b ;e ,f ;1\right)
 = 
\frac{ \Gamma \left(e +f -a -b +n \right) \Gamma \left(e \right)}{\Gamma \left(e +n \right) \Gamma \left(-a -b +e +f \right)}{}_{3}^{}{\moversetsp{}{\mundersetsp{}{F_{2}^{}}}}\left(-n ,f-b ,f -a ;e +f-a -b  ,f ;1\right)
\label{Thom8n}
\end{equation}
and, following the reversal of the right-hand side and redefinition of variables, we find
\begin{align} \nonumber
{}_{3}^{}{\moversetsp{}{\mundersetsp{}{F_{2}^{}}}}&\left(-n ,a ,b ;e,f ;1\right)
 \\&= 
\frac{\left(-1\right)^{n}\,(a)_{n}\,(b)_{n} }{ (e)_{n} (f)_{n}} 
\;{}_{3}^{}{\moversetsp{}{\mundersetsp{}{F_{2}^{}}}}\left(-n ,1-f -n,1-e -n;1-b -n,1-a -n ;1\right)\,,
\label{R8Rev}
\end{align}
a seemingly weighty, but fundamentally elementary, transformation. Similarly, the Thomae transformation \eqref{Thom9} with $c=-n$ generates the equivalent of \eqref{Thom8n} and \eqref{R8Rev} with $e\leftrightarrow f$. It is noteworthy that although \eqref{Thom8n} loses the anticipated symmetry between interchange of $e$ and $f$, it is recovered in \eqref{R8Rev}. See \eqref{We1} and \eqref{H12h} where similar situations arise.

In addition to the above, in 1952, Weber and Erdelyi \cite[Eq.(16)]{WebErdelyi} have used the simple integral representation of 3F2(1) 
\begin{equation}
{}_{3}^{}{\moversetsp{}{\mundersetsp{}{F_{2}^{}}}}\left(-n,a ,b  ;e ,f ;1\right)
 = 
\frac{\Gamma \left(f \right)}{\Gamma \left(b \right) \Gamma \left(f -b \right)}
 \int_{0}^{1}{}_{2}^{}{\moversetsp{}{\mundersetsp{}{F_{1}^{}}}}\left(a ,-n ;e ;t \right) t^{b -1} \left(1-t \right)^{f -b -1}d t 
\label{Eq16}
\end{equation}
followed by a simple linear transform and term-by-term integration, to obtain the transformation 
\begin{equation}
{}_{3}^{}{\moversetsp{}{\mundersetsp{}{F_{2}^{}}}}\left(-n ,a ,b ;e ,f ;1\right)
 = 
\frac{\Gamma \left(f \right) \Gamma \left(n+f -a \right) }{\Gamma \left(n +f \right) \Gamma \left(f-a \right)}\,{}_{3}^{}{\moversetsp{}{\mundersetsp{}{F_{2}^{}}}}\left(-n ,a ,e-b ;e ,a -f -n +1;1\right).
\label{We1}
\end{equation}
Noting that the symmetry between $e$ and $f$ was lost on the right-hand side, as in \eqref{Thom8n}, they rewrote \eqref{We1} with that interchange. After equating the right-hand sides of both and with a suitable redefinition of variables they found the transformation 
\begin{align} \nonumber
{}_{3}^{}{\moversetsp{}{\mundersetsp{}{F_{2}^{}}}}\left(-n,a ,b  ;e ,f ;1\right)
 &= 
\frac{\Gamma \left(e \right) \Gamma \left(e -a +n \right)  \Gamma \left(a -f +1\right) \Gamma \left(1-f -n \right)}{\Gamma \left(n +e \right) \Gamma \left(e -a \right) \Gamma \left(a -f -n +1\right) \Gamma \left(1-f \right)}\\
&\times{}_{3}^{}{\moversetsp{}{\mundersetsp{}{F_{2}^{}}}}\left(-n ,a ,1-\sigma;a -f -n +1,1-e +a -n ;1\right)\,
\end{align}
where  $\sigma\equiv -a - b + e + f + n$, coinciding with that of Sheppard - see \eqref{Sheppard}. In 1992, Rao, Van der Jeugt, Raynal, Jagannathan and Rajeswari \cite{RaoJeu1992} carried this program further by recursively reproducing the Weber/Erdelyi transformation employing symmetries between exchanges among other variables in the same manner as that outlined in Remark \ref{sec:RemXfm} (Appendix \ref{sec:nonTerm}), to find that the transformations of terminating 3F2(1) yield a group of 72 transformations of which 18 are non-trivial. These are listed in \cite[Appendix]{RaoJeu1992}; the basic ones are reproduced here -- see Appendix \ref{sec:termCase} -- and are labelled RJRJR throught this work.

In addition to the two-part Thomae transformations, Bailey also presented three-part transformations due to Whipple. These were summarized again by both Luke \cite{Luke} and Slater \cite{Slater}, but remain rarely employed by analysts, perhaps because the notation is largely incomprehensible, although Milgram \cite[Appendix C]{MilgVar} attempted to provide a clearer interpretation, unaware that both Rao, Van der Jeugt, Raynal, Jagannathan and Rajeswari\cite{RaoJeu1992} (see above) and Beyer, Louck and Stein \cite{BeyLouStein1987} had done the same years earlier, by developing a group theoretical basis for both generalized and terminating cases of 3F2(1). 

Along with collecting general identities, in the book {\it A=B} \cite{AeqB}, new algorithms for evaluating general hypergeometric functions, usually labelled by the acronym WZ, were introduced. These largely rely on recursion between contiguous elements and generally yield identities that cannot be foretold -- for recent examples, see \cite{CampLevrie2024} and \cite{Campbell2022}. 
In practice, the WZ algorithms fail in those exceptional cases where integer parameters generate the presence of digamma functions, due to limiting operations that render terms to become non-hypergeometric.

In 2010, faced with the need to evaluate a particular 3F2(1), following Whipple and Watson, Milgram \cite{Milg3F2Rev} collected all available results from the literature and tested each one against all nine of the independent Thomae transformations in a search for new identities. 62 were found and a computerized database of 3F2(1) evaluations was established. Following that effort, in 2011, Milgram \cite{Milgram447} attempted to isolate more new evaluations by recursively scanning all two-part, and a selection of three-part transformations, starting from the previous database, searching for identities where two of the three terms in the corresponding transformation were known in summable or closed form. As a consequence of that study, 259 new evaluations were found, the older database was extended to 469 entries and made available in machine searchable form (\cite[Appendix E]{Milgram447}). 

Since that time, several purportedly new identities have appeared in the literature, many of which are based on either the WZ algorithm or integration identities. Although many of these are now labelled ``exotic", most are already known to the database (e.g. \cite[Eq. (6)]{Campbell2022} - see \eqref{LukeGen} below), usually recognizable, but not recognized, as a variation of the Minton/Karlsson/Rosengren/Gasper/Thomae (MKRGT) identities or Chu's extension of DWW identities \cite{ChuW}. In fact, in \cite{Milgram447}, I opined that no new closed form evaluations of 3F2(1) would be found by studying transformations until a fundamentally new and independent evaluation was discovered. 

In 2012 such an evaluation was found by Miller and Paris \cite{Miller&Paris}, closely followed in 2015 by a second evaluation developed by Shpot and Srivastava\footnote{They independently rediscovered a 1988 identity of Gottschalk and Maslen \cite{Gottschalk}.}  \cite{Shpot&Sriv}. One aim of this work is to revisit these two identities to show that they are related and easily derived (see Section \ref{sec:relax}), to investigate the consequences and extend the derivation to similar identities of other related cases, all of which are characterized by the fact that they include one or more parameters that are integers, or are separated by integers. Throughout, extensive use is made of the Thomae transformations (Appendix \ref{sec:nonTerm}), those of Rao et. al. (RJRJR - see Appendix \ref{sec:termCase}) plus the simplest (see \eqref{HAll1}) of the known three-part transformations among 3F2(1). Concurrently, I also study a collection of ``new" identities extracted from the recent literature, in particular special case identities developed by Karp and Kalmykov \cite{KarpAndKalmy} and Karp and Prilepkina \cite{Karp&Prep2018}, where $c=f+n$ and $b:=b+p$ in \eqref{3F2Def} and  ``exotic" identities introduced by Chen and Chu \cite{ChuChen1}. The goal is to demonstrate that many new and interesting evaluations of 3F2(1) are lurking among the existing known transformation sets.

\subsection{Recent extensions of the familiar terminating identities} \label{extensions}
 A general $p-$balanced extension of Minton's zero-balanced identity \eqref{Minton} is
\begin{align} \nonumber
{}_{3}^{}{\moversetsp{}{\mundersetsp{}{F_{2}^{}}}}&\left(-m ,a ,c +n ;c ,p +n +a -m ;1\right)
 = 
\frac{\Gamma \left(c \right) \Gamma \left(p +n +a -m \right) \Gamma \left(1+n \right) }{\Gamma \left(c +n \right) \Gamma \left(a -m +p -c \right)}\\
&\times\moverset{n}{\munderset{0\leq k =1-p+m}{\sum}}\frac{\Gamma \left(p -c +a -m +k \right) \Gamma \left(p +k \right) \left(-1\right)^{k}}{\Gamma \left(1+n -k \right) \Gamma \left(k +1\right) \Gamma \left(p +a +k \right) \Gamma \left(p -m +k \right)}\,,
\label{MintonCor}
\end{align}
a proof of which can be found in Appendix \ref{sec:Proofs}, Corollary (\ref{sec:MintonExCor}). Notice that if $m\geq p+n$, the right-hand side of \eqref{MintonCor} vanishes, reproducing Karlsson's result \cite[Eq. (12)]{Karlsson}.

A further extension of the previous work is due to Karp and Prilepkina \cite{Karp&Prep2018} who, in 2018, studied extensions of the MKRG-type identities, motivated by the omission of the second right-hand side term (added here) in Minton's original work (cf. \eqref{Minton} i.e. $n>m$), and the further observation that, although the $3F2(1)$ in the MKRG identities (e.g. \eqref{K10}) diverges if $1-a-n\leq 0$ and the laws of analytic continuation state that the right-hand side is the continuation of the left-hand side, if $a=m$, the left-hand side terminates (and thus doesn't diverge) if the inequality is true, but the purported identity fails numerically. Based on this observation, Karp and Prilepkina derive a generalized form of \eqref{K10} valid for zero or negatively balanced terminating sums; reduced to the case of a $3F2(1)$, they find \cite[Eq. (2.3)]{Karp&Prep2018}, for $0\leq n\leq m-1$,
\begin{align} \nonumber
{}_{3}^{}{\moversetsp{}{\mundersetsp{}{F_{2}^{}}}}&\left(-m ,a ,c +n ;a +1,c ;1\right)
  \\ &=
\frac{\Gamma \left(1+m \right) \Gamma \left(a +1\right) \Gamma \left(n +c -a \right) \Gamma \left(c \right)}{\Gamma \left(a +1+m \right) \Gamma \left(c -a \right) \Gamma \left(c +n \right)}
-\frac{\left(-1\right)^{n}\,\Gamma \left(1+m \right) \Gamma \left(c \right) q_{m}}{\Gamma \left(c +n \right)}
\label{T2a}
\end{align}
where they define the constants $q_{m}$ analytically based on properties of Meijer's G-function (see \eqref{QM} below). In the case that $m:=m-1,~ n=m$, they present an example that exactly coincides with \eqref{Minton} for equivalent parameters. Of more interest is the negatively balanced, specialized example $m:=m-2,n=m$, where they find, for $m\geq 2$,
\begin{align} \nonumber
{}_{3}^{}{\moversetsp{}{\mundersetsp{}{F_{2}^{}}}}&\left(2-m,a ,c +m  ;c ,a +1;1\right)
=\frac{\Gamma \left(m -1\right) \Gamma \left(c \right)}{\Gamma \left(c +m \right)} \\ 
&\times\left(\frac{\Gamma \left(c -a +m \right) \Gamma \left(a +1\right)}{\Gamma \left(c -a \right) \Gamma \left(a +m -1\right)}+\left(-1\right)^{m}\,a \left(1-a +m \left(c +m -2\right)\right)\right) 
\label{KPX2}
\end{align} 
thereby extending \eqref{MintonCor} to the range $p<0$.
\begin{rem}
It is possible, and left as an exercise for the reader, to independently obtain \eqref{KPX2} from the $p=1,~n=m$ version of \eqref{MintonCor} by utilizing recursion formulas (e.g. \eqref{RainvEq15a}). Hint: utilize the negatively balanced, but terminating sum:
\begin{equation}
{}_{2}^{}{\moversetsp{}{\mundersetsp{}{F_{1}^{}}}}\left(2-m ,c +m -1;c ;1\right)
 = 
\left(-1\right)^{m}\,\frac{\Gamma \left(m \right) \Gamma \left(c \right)}{\Gamma \left(c +m -2\right)}\,.
\label{Sc2F1}
\end{equation}
\end{rem}

\subsection{Simple alternative derivations} \label{sec:Altern}

In 2010, Miller and Srivastava \cite{MillerSriv2010} employed ``elementary identities for binomial coefficients ${\binom{j}{k}}$ and Stirling numbers of the second kind $\genfrac{\lbrace}{\rbrace}{0pt}{}{j}{k}$" to rederive the Karlsson-Minton summation formulas in a manner that would be ``accessible to non-specialists". Here, an alternative simpler way to achieve that same goal is presented, by applying the Thomae identities to the common left-hand side and observing that two variations produce terminating sums. Explicitly, in addition to the KRG identities, we find the following two variants by applying Thomae identities \eqref{Thom1}, and \eqref{Thom4} respectively:
\begin{align} \nonumber
{}_{3}^{}{\moversetsp{}{\mundersetsp{}{F_{2}^{}}}}\left(a ,b ,c +n ;c ,e ;1\right)
 &= 
\frac{ \Gamma \left(e -n -a -b \right) \Gamma \left(e \right) \Gamma \left(c \right)}{\Gamma \left(c +n \right) \Gamma \left(-b +e -n \right) \Gamma \left(e -a -n \right)}\;\\&\times{}_{3}^{}{\moversetsp{}{\mundersetsp{}{F_{2}^{}}}}\left(-n,e -n -a -b ,e -c -n  ;e-a  -n ,e -b -n ;1\right)
\label{X1}\\
&=
\frac{ \Gamma \left(e -n -a -b \right) \Gamma \left(e \right)}{\Gamma \left(e -a \right) \Gamma \left(-b +e -n \right)}\;{}_{3}^{}{\moversetsp{}{\mundersetsp{}{F_{2}^{}}}}\left(-n ,a,c-b ;c,e-b -n ;1\right)\,.
\label{X4}
\end{align}
Setting $e=b+1$ in either of \eqref{X1} or \eqref{X4}, yields Karlsson's identity \eqref{K10}. Similarly, setting $e=b+1+m$ in \eqref{X1} yields an equivalent variation of Gasper's identity \eqref{G81M}
\begin{align} \nonumber
{}_{3}^{}{\moversetsp{}{\mundersetsp{}{F_{2}^{}}}}&\left(a ,b ,c +n ;c ,b +m +1;1\right) = 
\frac{\Gamma \left(1+n \right) \Gamma \left(b +m +1\right) \Gamma \left(c \right)}{\Gamma \left(c +n \right) \Gamma \left(b +m +1-c -n \right)}\\
&\times \moverset{n}{\munderset{0\leq k=n-m}{\sum}}\frac{\left(-1\right)^{k}\,\Gamma \left(m +1-a -n +k \right) \Gamma \left(b +m +1-c -n +k \right) }{\Gamma \left(1-k +n \right) \Gamma \left(k +1\right) \Gamma \left(m +1-n +k \right) \Gamma \left(b +m +1-a -n +k \right)}\,.
\label{X1b}
\end{align}
Rather than setting $m=0$ to reproduce \eqref{K10} as in \eqref{G81M}, set $n=1$ in \eqref{X1b} to obtain a closed  form special case equivalent to the $n=1$ instance of Rosengren's identity \eqref{RosengEb}:
\begin{equation}
{}_{3}^{}{\moversetsp{}{\mundersetsp{}{F_{2}^{}}}}\left(a ,b ,c +1;c ,b +m +1;1\right)
 = 
\frac{\Gamma \left(m -a \right) \left(\left(m -a \right) c +a\,b \right) \Gamma \left(b +m +1\right)}{c\,\Gamma \left(m +1\right) \Gamma \left(b +m -a +1\right)}\,.
\label{X1c}
\end{equation}
%
It is also instructive to consider the Thomae transformation \eqref{X1} with $e=a+1$ giving
\begin{align} \nonumber
{}_{3}^{}{\moversetsp{}{\mundersetsp{}{F_{2}^{}}}}&\left(a ,b ,c +n ;c ,a +1;1\right)
 = 
\frac{\Gamma \left(1+n \right)  \Gamma \left(a +1\right) \Gamma \left(c \right)}{\Gamma \left(c +n \right) \Gamma \left(1-n -c +a \right)}\\
&\times\moverset{\infty}{\munderset{k =0}{\sum}}\frac{ \left(-1\right)^{k}\,\Gamma \left(1-b -n +k \right) \Gamma \left(1-n -c +a +k \right)}{\Gamma \left(1-k +n \right) \Gamma \left(k +1\right) \Gamma \left(-b +a +1-n +k \right) \Gamma \left(1-n +k \right)}\,,
\label{X1s}
\end{align}
where only the term associated with the index $k=n$ contributes, unless $b$ is an integer. In that case, set $b=-m$ and thereby connect Minton's transformation \eqref{Minton} with KRG type transformations. If $m>n-1$, then \eqref{X1s} becomes
\begin{equation}
{}_{3}^{}{\moversetsp{}{\mundersetsp{}{F_{2}^{}}}}\left(-m ,a ,c +n ;c ,a +1;1\right)
 = 
\frac{\left(-1\right)^{n}\,\Gamma \left(m +1\right) \Gamma \left(1-c +a \right) \Gamma \left(a +1\right) \Gamma \left(c \right)}{\Gamma \left(c +n \right) \Gamma \left(1-n -c +a \right) \Gamma \left(m +a +1\right)},
\label{Case1}
\end{equation}
corresponding to the special case \eqref{MintonCor} when $p=m-n+1>0$. If $m\leq n-1$ then the sum in \eqref{X1s} must be split into parts $k=0\dots n-m-1$ and $k=n$, since all terms between $n-m\leq k<n$ vanish, eventually producing
\begin{align} \nonumber
{}_{3}^{}{\moversetsp{}{\mundersetsp{}{F_{2}^{}}}}&\left(-m ,a ,c +n ;c ,a +1;1\right)
 = \frac{\left(-1\right)^{n}\,\Gamma \left(m +1\right) \Gamma \left(1-c +a \right) \Gamma \left(a +1\right) \Gamma \left(c \right)}{\Gamma \left(c +n \right) \Gamma \left(1-n -c +a \right) \Gamma \left(m +a +1\right)}\\
&+\frac{\Gamma \left(1+n \right) \Gamma \left(a +1\right) \Gamma \left(c \right) }{\Gamma \left(c +n \right) \Gamma \left(1-n -c +a \right)}
\moverset{n-m -1}{\munderset{k =0}{\sum}}\frac{\left(-1\right)^{m +k}\,\Gamma \left(1-n -c +a +k \right)}{\Gamma \left(n-m -k \right) \left(n-k \right) \Gamma \left(k +1\right) \Gamma \left(m +a +1-n +k \right)}\,
\label{Case2a}
\end{align}
where the identity
\begin{equation}
\munderset{b \rightarrow -m}{\mathrm{lim}}\frac{\Gamma \left(1-b -n +k \right)}{\Gamma \left(1-n +k \right)}
 = 
\frac{\left(-1\right)^{-m}\,\Gamma \left(n-k  \right)}{\Gamma \left(-m +n -k \right)}\,,\hspace{20pt} m+1\leq n,~0\leq k\leq n-m-1,
\label{Leg}
\end{equation}
has been employed. Comparison of \eqref{Case2a} and \eqref{T2a} produces a simple identification of Karp and Prilepkina's constants $q_{m}$ (see \eqref{T2a}) as they relate to 3F2(1):
\begin{equation}
q_{m} = 
\frac{\left(-1\right)^{n}\,\Gamma \left(a +1\right) \Gamma \left(1+n \right)}{\Gamma \left(1-n -c +a \right) \Gamma \left(m +1\right)}\moverset{n-m -1}{\munderset{k =0}{\sum}}\frac{\left(-1\right)^{m +k}\,\Gamma \left(1-n -c +a +k \right)}{\left(k-n \right) \Gamma \left(m +a +1-n +k \right) \Gamma \left(n-m  -k \right) \Gamma \left(k +1\right)}\,.
\label{QM}
\end{equation}
Setting $m:=m-2$ followed by $n:=m$ in \eqref{Case2a} reproduces Karp and Prilepkina's example \eqref{KPX2}.

In addition to providing simple derivations, if we apply the Thomae identities to Gasper's identity \eqref{Gasper}, it turns out that three variants (\eqref{Thom1}, \eqref{Thom3} and \eqref{Thom6} respectively) resolve as terminating sums. Those produce the following new transformations that can be added to our rapidly growing collection:
\begin{align} \nonumber
{}_{3}^{}{\moversetsp{}{\mundersetsp{}{F_{2}^{}}}}&\left(a ,b ,c +n ;c ,e ;1\right)
 = {}_{3}^{}{\moversetsp{}{\mundersetsp{}{F_{2}^{}}}}\left(-n ,c -a ,e -n -a -b ;e-a  -n ,1-n -a ;1\right)\\ 
&\times\frac{\Gamma \left(e \right) \Gamma \left(1-a \right) \Gamma \left(n +c -b \right) \Gamma \left(c \right)  \Gamma \left(e -n -a -b \right) \Gamma \left(1+b -c -n \right)}{\Gamma \left(1-c +b \right) \Gamma \left(e-a  -n \right) \Gamma \left(1-n -a \right) \Gamma \left(c -b \right) \Gamma \left(e -b \right) \Gamma \left(c +n \right)}\label{GS1}\\ \nonumber
=&{}_{3}^{}{\moversetsp{}{\mundersetsp{}{F_{2}^{}}}}\left(-n ,b ,e-c  -n ;e-a  -n ,1+b -c -n ;1\right)\\
&\times\frac{\Gamma \left(e \right) \Gamma \left(n +c -b \right) \Gamma \left(c \right)  \Gamma \left(e -n -a -b \right)}{\Gamma \left(e-a  -n \right) \Gamma \left(c -b \right) \Gamma \left(e -b \right) \Gamma \left(c +n \right)}
\label{GS3}\\ \nonumber
=& {}_{3}^{}{\moversetsp{}{\mundersetsp{}{F_{2}^{}}}}\left(-n ,1-e +b,1-c -n ;1-n -a ,1+b -c -n ;1\right)\\
&\times\frac{\Gamma \left(e \right) \Gamma \left(1-a \right) \Gamma \left(n +c -b \right) \Gamma \left(c \right) \Gamma \left(e -n -a -b \right)}{\Gamma \left(e -a \right) \Gamma \left(1-n -a \right) \Gamma \left(c -b \right) \Gamma \left(e -b \right) \Gamma \left(c +n \right)}\,.
\label{GS6}
\end{align}
These can be used  to study interesting specializations; for example if we set $e=n+a+b+m$ in \eqref{GS1}, we obtain the $m-$balanced transformation
\begin{align} \nonumber
{}_{3}^{}{\moversetsp{}{\mundersetsp{}{F_{2}^{}}}}&\left(a ,b ,c +n ;c ,a +b+n +m ;1\right)\\ \nonumber
& = 
\frac{\Gamma \left(n +a +b +m \right) \Gamma \left(1-a \right) \Gamma \left(n +c -b \right) \Gamma \left(c \right) \Gamma \left(1+n \right) \Gamma \left(1+b -c -n \right)}{\Gamma \left(c -a \right) \Gamma \left(1-c +b \right) \Gamma \left(c -b \right) \Gamma \left(n +a +m \right) \Gamma \left(c +n \right)}\\
& \times\moverset{n}{\munderset{k =0}{\sum}}\frac{\Gamma \left(m +k \right) \left(-1\right)^{k}\,\Gamma \left(c -a +k \right)}{\Gamma \left(1-k +n \right) \Gamma \left(k +1\right) \Gamma \left(b +m +k \right) \Gamma \left(1-n -a +k \right)}\,,
\label{S1b}
\end{align}
a variation of \eqref{G81M} and \eqref{RosengEb}. Other specializations obviously follow; for example, if $a=j$:
\begin{align}
{}_{3}^{}{\moversetsp{}{\mundersetsp{}{F_{2}^{}}}}&\left(b ,j ,n +c ;c ,n +j +b +m ;1\right)
 = 
\frac{\Gamma \left(n +j +b +m \right) \Gamma \left(c \right) \Gamma \left(1+n \right) }{\Gamma \left(j \right) \Gamma \left(c -j \right) \Gamma \left(n +j +m \right) \Gamma \left(n +c \right)}\\
& \times
 \moverset{n}{\munderset{k =0}{\sum}}\frac{\Gamma \left(m +k \right) \Gamma \left(c -j +k \right) \Gamma \left(n +j -k \right)}{\Gamma \left(1-k +n \right) \Gamma \left(k +1\right) \Gamma \left(b +m +k \right)}\,.
\label{S1B}
\end{align}

\section{A Challenge} \label{sec:Hunt}

In the closing statement of their paper, Miller and Paris write: {\it ``...there remains the open problem of deducing a summation formula for the series $_{3}F_{2}(a,b,c;e,a+n;1)$... We hope that the developments presented here will stimulate interest in this problem."} The intent in this section is to employ an elementary, three-part transformation, to confront the Miller/Paris challenge.  We begin with the simple, well-known  three-part transformation among hypergeometric $3F2(1)$, to be found (with $a\leftrightarrow c$) in Luke's book \cite[Eq. 3.13.3(11)]{Luke}:
\begin{align} \nonumber
{}_{3}^{}{\moversetsp{}{\mundersetsp{}{F_{2}^{}}}}&\left(a ,b ,c ;e ,f ;1\right) = \Gamma \left(f \right) \Gamma \left(e \right)  \Gamma \left(1-c \right)\\ \nonumber
 &
\times\left(\frac{\Gamma \left(b -a \right)}{\Gamma \left(1-c +a \right) \Gamma \left(f -a \right) \Gamma \left(e -a \right) \Gamma \left(b \right)}\,{}_{3}^{}{\moversetsp{}{\mundersetsp{}{F_{2}^{}}}}\left(a ,1-e +a ,1-f +a ;1-b +a ,1-c +a ;1\right)\right. \\
&\left.+\,\frac{\Gamma \left(a -b \right)}{\Gamma \left(1-c +b \right) \Gamma \left(f -b \right) \Gamma \left(e -b \right) \Gamma \left(a \right)}\,{}_{3}^{}{\moversetsp{}{\mundersetsp{}{F_{2}^{}}}}\left(b ,1-e +b ,1-f +b ;1-c +b ,1-a +b ;1\right)\right)\,. 
\label{HAll1}
\end{align}
In order to obtain a ``summable" evaluation (as opposed to a transformation), the right-hand side of \eqref{HAll1} must reduce to a representation that does not include an infinite summation, since the general form of the left-hand side certainly does -- see \eqref{3F2Def}. A simple attempt to achieve this goal involves reducing the right hand side to terminating sums. With the Miller/Paris challenge in mind, in \eqref{HAll1}, let $f=a+n$ to find
\begin{align} \nonumber
{}_{3}^{}{\moversetsp{}{\mundersetsp{}{F_{2}^{}}}}&\left(a ,b ,c ;e ,a +n ;1\right)
 =  \Gamma \left(1-c \right) \Gamma \left(a +n \right) \Gamma \left(e \right)\\ \nonumber
 &\times
\left(\Gamma \left(a -b \right)\frac{ {}_{3}^{}{\moversetsp{}{\mundersetsp{}{F_{2}^{}}}}\left(b ,1 +b-e ,1+b-a -n  ;1+b-c  ,1+b-a  ;1\right)}{\Gamma \left(1+b-c  \right) \Gamma \left(a +n -b \right) \Gamma \left(e -b \right) \Gamma \left(a \right)}\right.\\ 
&\left. \hspace{20pt}
+\,\Gamma\left(b-a  \right)\frac{ {}_{3}^{}{\moversetsp{}{\mundersetsp{}{F_{2}^{}}}}\left(a ,1-n ,1+a-e  ;1+a-b  ,1+a-c  ;1\right) }{\Gamma \left(1+a-c  \right) \Gamma \left(n \right) \Gamma \left(e -a \right) \Gamma \left(b \right)}\right)
\label{HAll2}
\end{align}
and notice that the second term within the parentheses is a terminating sum when $n\geq 1$ due to \eqref{NegTop}. To satisfy the Miller/Paris challenge, it remains to determine if the first hypergeometric term within those parentheses either terminates, or can be transformed into an expression that terminates, for general values of $a,b,c,e$ and $n>0$. There are several properties of \eqref{HAll2} worth noting:
\begin{itemize}
\item{}
There exists a certain degree of symmetry between the terms under the exchange $a\Leftrightarrow b$;
\item{}
On the left-hand side, there is complete symmetry under the exchange $b\leftrightarrow c$ and this symmetry must somehow extend to the right-hand side;;
\item{}
With respect to the first term, there exists a parametric difference equal to $-n$ between the third top parameter $1+b-a-n$ and the second bottom parameter $1+b-a$, indicating that that term is itself of the same form as that which appears on the left-hand side;
\item{}
Relaxing the ground-rules slightly, setting $e=b+1$ produces the known identity \cite{MilgRaoComment}
\begin{align} \nonumber
{}_{3}^{}{\moversetsp{}{\mundersetsp{}{F_{2}^{}}}}&\left(a ,b ,c ;b +1,a +n ;1\right)
 = 
\frac{\Gamma \left(b +1\right) \Gamma \left(a -b \right) \Gamma \left(1-c \right) \Gamma \left(a +n \right)}{\Gamma \left(1-c +b \right) \Gamma \left(a +n -b \right) \Gamma \left(a \right)}\\
&-\frac{ b\,\Gamma \left(a +n \right) \Gamma \left(1-c \right)}{\left(a -b \right) \Gamma \left(1-c +a \right) \Gamma \left(n \right)}\;{}_{3}^{}{\moversetsp{}{\mundersetsp{}{F_{2}^{}}}}\left(a ,a -b ,1-n ;1-b +a ,1-c +a ;1\right)\,,
\label{Cor2p2}
\end{align}
a terminating sum. A slightly greater relaxation of the ground-rules produces the identity \eqref{Hnew3} below, if we set $e=a+m$.
\item{}
Finally, in the remainder of this work, we study the case where the ground-rules are relaxed by setting $c=m,m>1$, since the parametric difference between the first top parameter $b$ and the first bottom parameter $1+b-c$ becomes a positive integer of value $m-1$ and so the first term transmutes into a terminating form through the use of either \eqref{Karl8} or \eqref{RosenG}, recognizing that, the presence of the overall factor $\Gamma(1-c)$ negates the simplicity of such an approach. 
\end{itemize}

\subsection{Four attempts}
\subsubsection{A functional equation}
Since the first term in \eqref{HAll2} is of the same form as the left-hand side, consider redefining the parameters in that identity, effectively treating the identity as a functional equation. To verify that this leads either to a tautology or another equivalent three-part identity is left as an exercise for the reader.

\subsubsection{Relaxing the constraint again yields simple derivations}\label{sec:relax}

By again relaxing the constraint that \eqref{HAll2} must be expressed in the form of a terminating sum for general values of $a,b,c,e$ and $n>0$, it becomes possible to make use of the observed symmetry by letting $e=b+m$ in \eqref{HAll2} and obtain
\begin{align} \nonumber
{}_{3}^{}{\moversetsp{}{\mundersetsp{}{F_{2}^{}}}}&\left(a ,b ,c ;a+n, b +m ;1\right)
 = \Gamma \left(a +n \right) \Gamma \left(b +m \right) \Gamma \left(1-c \right)\\ \nonumber
&\times\left(\frac{\Gamma \left(b-a \right){}_{3}^{}{\moversetsp{}{\mundersetsp{}{ \,F_{2}^{}}}}\left(1-n ,a ,1+a-b -m  ;1+a-b  ,1+a-c  ;1\right)}{\Gamma \left(1+a-c  \right) \Gamma \left(n \right) \Gamma \left(b +m -a \right) \Gamma \left(b \right)}\right.\\ &\left.
+\hspace{25pt}\frac{\Gamma \left(a -b \right){}_{3}^{}{\moversetsp{}{\mundersetsp{}{F_{2}^{}}}}\left(1-m ,b ,1+b-a -n  ;1-c +b ,1 +b-a ;1\right)}{\Gamma \left(1+b-c  \right) \Gamma \left(a +n -b \right) \Gamma \left(m \right) \Gamma \left(a \right)}\right) \,,
\label{hnew3}
\end{align}
in which case both terms on the right-hand side of \eqref{hnew3} are terminating sums and therefore \eqref{HAll1} has been reduced to finite summations. Explicitly, in summation form
\begin{align} \nonumber
{}_{3}^{}{\moversetsp{}{\mundersetsp{}{F_{2}^{}}}}&\left(a ,b ,c ;a+n,b +m  ;1\right)
 = \frac{\Gamma \left(1-c \right) \Gamma \left(b +m \right) \Gamma \left(a +n \right)}{\Gamma \left(b \right) \Gamma \left(a \right)}\\ \nonumber
 &\times
\left( \left(-1\right)^{n}\moverset{m-1}{\munderset{k =0}{\sum}}\frac{\Gamma \left(b +k \right) \Gamma \left(1 +b +k-a -n \right) \left(-1\right)^{k}}{\Gamma \left(m -k \right) \Gamma \left(k +1\right) \Gamma \left(1 +b +k-a \right) \Gamma \left(1+b +k -c \right)} \right. \\
& \left.\hspace{25pt} +\left(-1\right)^{m} \moverset{n -1}{\munderset{k =0}{\sum}}\frac{\Gamma \left(a +k \right) \Gamma \left(1-b -m +a +k \right) \left(-1\right)^{k}}{\Gamma \left(n -k \right) \Gamma \left(k +1\right) \Gamma \left(1-b +a +k \right) \Gamma \left(1-c +a +k \right)} \right)\,,
\label{Hnew3}
\end{align}
exactly reproducing the identity first obtained by Gottschalk and Maslen \cite[Eq. 23]{Gottschalk} in 1988, later independently rediscovered in 2015 by Shpot and Srivastava \cite{Shpot&Sriv} and misinterpreted by myself \cite{MilgVar} somewhat later. If $n=1$, \eqref{Hnew3} reduces to the known  identity \eqref{Cor2p2}
\begin{align} \nonumber
{}_{3}^{}{\moversetsp{}{\mundersetsp{}{F_{2}^{}}}}&\left(a ,b ,c ;1+a  ,b+m ;1\right)=\frac{\left(-1\right)^{m}\,\Gamma \left(1+a -b -m \right) \Gamma \left(1+a \right) \Gamma \left(b +m \right) \Gamma \left(1-c\right)}{\Gamma \left(1+a -b \right) \Gamma \left(b \right) \Gamma \left(1+a -c \right)}
\\ 
& + 
\frac{\Gamma \left(b -a \right)  \Gamma \left(1+m -c \right) \Gamma \left(1+a \right) \Gamma \left(b +m \right)}{\Gamma \left(c \right) \Gamma \left(m -c +b \right)}
\moverset{m -1}{\munderset{k =0}{\sum}}\frac{\left(-1\right)^{k}\Gamma \left(c -1-k \right) }{\Gamma \left(a -k \right) \Gamma \left(m -k \right) \Gamma \left(b -a +1+k \right)},
\label{L2p2}
\end{align}
written in sum form. Seizing the moment, apply the two-part Thomae transformation \eqref{Thom9} to the left hand side of \eqref{Hnew3} and find
%
\begin{align} \nonumber
{}_{3}^{}{\moversetsp{}{\mundersetsp{}{F_{2}^{}}}}&\left(a ,b ,m ;m+n  ,c ;1\right)
 = 
\frac{\pi  \left(-1\right)^{m}\,\Gamma \left(n -a -b +c \right) \Gamma \left(n +m \right) \Gamma \left(c \right)}{\Gamma \left(c -m \right) \Gamma \left(-a +n +m \right) \Gamma \left(-b +n +m \right)}
\\ \nonumber
&\times\left(\frac{1}{\sin \left(\pi \,c \right) \Gamma \left(c -b \right) \Gamma \left(c -a \right) \Gamma \left(n \right)}
\moverset{m-1}{\munderset{k =0}{\sum}}\frac{\Gamma \left(1-b+k \right) \Gamma \left(1-a +k \right) \Gamma \left(n +m -1-k \right)}{\Gamma \left(k +1\right) \Gamma \left(m -k \right) \Gamma \left(2-c +k \right)}
\right. \\
&\left.
-\frac{ \left(-1\right)^{n} \sin \left(\pi  \left(b -c +a \right)\right)}{\sin \left(\pi \,a \right) \sin \left(\pi \,b \right) \Gamma \left(m \right)}\moverset{n -1}{\munderset{k =0}{\sum}}\frac{\Gamma \left(m +k \right) \Gamma \left(1+a +b -c -n +k \right) \left(-1\right)^{k}}{\Gamma \left(b +1+k -n \right) \Gamma \left(1+a -n +k \right) \Gamma \left(k +1\right) \Gamma \left(n -k \right)}\right)\,, 
\label{Hnew}
\end{align}
equivalent to the identity that first appeared in the 2012 Miller/Paris paper \cite{Miller&Paris}.
%
Again, following the footsteps of Whipple and Watson, apply the Thomae transformation \eqref{Thom3} to the identity \eqref{Hnew3}, revealing another, possibly new, identity:
\begin{align} \nonumber
{}_{3}^{}{\moversetsp{}{\mundersetsp{}{F_{2}^{}}}}&\left(a ,b ,n ;c ,a+b -c  +n +m ;1\right)=\\ \nonumber
&  
\frac{ \Gamma \left(1-b \right) \Gamma \left(c \right) \Gamma \left(b -c +a +n +m \right) \Gamma \left(m \right)}{\Gamma \left(c -b \right) \Gamma \left(a -c +n +m \right) \Gamma \left(c -n \right) \Gamma \left(b +m +a -c \right) \Gamma \left(n \right)}\\ \nonumber
&\times \moverset{m -1}{\munderset{k =0}{\sum}}\frac{\Gamma \left(c -b -n -m +1+k \right) \Gamma \left(a -c +1+k \right) \Gamma \left(n +m -1-k \right)}{\Gamma \left(1+k \right) \Gamma \left(m -k \right) \Gamma \left(2-b -m +k \right)}
\\ \nonumber
&-\frac{ \left(-1\right)^{n}\,\Gamma \left(1+b -c \right) \sin \left(\pi \,a \right) \Gamma \left(c \right) \Gamma \left(b -c +a +n +m \right)}{\sin \left(\pi  \left(c -a \right)\right) \Gamma \left(b \right) \Gamma \left(a -c +n +m \right)}\\
&\times \moverset{n -1}{\munderset{k =0}{\sum}}\frac{\left(-1\right)^{k}\Gamma \left(m +k \right) \Gamma \left(1-a +k \right) }{\Gamma \left(b-c  +m +1+k \right) \Gamma \left(c +1-a -n +k \right) \Gamma \left(k +1\right) \Gamma \left(n -k \right)}\,.
\label{Hnew2}
\end{align}
%
\begin{rem} {\bf Remarks}
\begin{itemize}
\item{}
The $m$-balanced identity \eqref{Hnew2}, a non-terminating analogue of the Saalsch\"{u}tz theorem \eqref{Saal}, may be new\footnote{The terminating case \eqref{Saal} is 1-balanced; a closed result for the corresponding  2-balanced case can be found in \cite[Eq.5]{MathWSaal} or see \eqref{XidbTwo} below. Setting $e=a+b+m-n-c$ in \eqref{Sheppard} sums the $m-$balanced terminating generalization.}. Setting $m=n=1$ reduces \eqref{Hnew2} to the known $1-$balanced identity \cite[Eq. 7.4.4(28)]{prudnikov}; if $n=1$ and $m>0$ we find the generalization for the $m$-balanced case with four independent parameters
\begin{align} \nonumber
{}_{3}^{}{\moversetsp{}{\mundersetsp{}{F_{2}^{}}}}&\left(1,a ,b ;c ,a -c +b +1+m ;1\right)=\frac{\Gamma \left(m \right) \Gamma \left(1+b -c \right) \Gamma \left(1-c +a \right) \Gamma \left(a -c +b +1+m \right) \Gamma \left(c \right)}{\Gamma \left(b \right) \Gamma \left(a -c +1+m \right) \Gamma \left(1-c +b +m\right) \Gamma \left(a \right)}\\
&+ 
\frac{\left(b +m +a -c \right) \Gamma \left(m \right) \left(c -1\right) \Gamma \left(1-b \right) }{\Gamma \left(a -c +1+m \right) \Gamma \left(c-b \right)}\moverset{m -1}{\munderset{k =0}{\sum}}\frac{\Gamma \left(c -b -m +k \right) \Gamma \left(1+a -c +k \right)}{\Gamma \left(k +1\right) \Gamma \left(-b +2-m +k \right)}\,;
\label{Hnew2n1}
\end{align}
\item{}
In the case that $b=-n,~m=1$, \eqref{Hnew2} reduces to a $b=n$ instance of the Saalsch\"{u}tz theorem\eqref{Saal}. 
\item{}
As in the case of DWW, the identities \eqref{Hnew3}, \eqref{Hnew} and \eqref{Hnew2} are closed under the action of the Thomae transformations; i.e. given any one of the triplet, the other two, and only the other two, can be found from the set of those transformations. See Appendix \ref{sec:Thomae}, Corollary \ref{sec:CorTom}.
\item{}
In any of these identities, the variables $a,b,c,e$ are assumed to take on arbitrary complex values such that the infinite series converges. In exceptional cases, a limiting calculation must be performed - see example \eqref{AOT} below.
\item{}
The identity \eqref{Hnew} generalizes well-known identities that appear in textbook form (e.g. Luke \cite[page 111] {Luke}) or \eqref{LukeGen} below, when $n$ and $m$ are small integers.
\end{itemize}
\end{rem}

\subsubsection{Check for useful Thomae equivalents} \label{sec:Equivs}

Here, we investigate the possibility that the non-terminating sum on the right-hand side of \eqref{HAll2} is reducible to a terminating sum by applying one of the Thomae transformations. The nine non-trivial Thomae progeny of that term are given below:
%
\begin{align} \nonumber
&{}_{3}^{}{\moversetsp{}{\mundersetsp{}{F_{2}^{}}}}\left(b ,1+b-e  ,1+b-a -n  ;1+b-a  ,1+b-c  ;1\right)
 \\ 
&=\frac{ \Gamma \left( e +n -b -c\right) \Gamma \left(1+b-c  \right) \Gamma \left(1+b-a  \right)}{\Gamma \left(1 +b-a -n \right) \Gamma \left(e +n-c  \right) \Gamma \left(1+n-c  \right)}\; _{3}F_{2}  \left(\right.\overset{ \displaystyle {n,e+n-b-c,a+n-c~}}{_ {\displaystyle {e+n-c,1+n-c}}}|1\left.\right)
 \label{Cx1}\\
&=\frac{ \Gamma \left(e +n-b -c  \right) \Gamma \left(1+b-c  \right) \Gamma \left(1+b-a  \right)}{\Gamma \left(1-e +b \right) \Gamma \left(e -c +n \right) \Gamma \left(1+e-a  -c \right)}
\; _{3}F_{2}  \left(\right.\overset{ \displaystyle {e+n-b-c,e-c,e-a~}}{_ {\displaystyle {e+n-c,1+e-a-c}}}|1\left.\right)
\label{Cx2}\\
&=\frac{ \Gamma \left(e +n-b -c  \right) \Gamma \left(1+b-a  \right)}{\Gamma \left(1-a \right) \Gamma \left(e +n -c \right)}
\; _{3}F_{2}  \left(\right.\overset{ \displaystyle {b,a+n-c,e-c~}}{_ {\displaystyle {e+n-c,1+b-c~}}}|1\left.\right)
\label{Cx3}\\
&=\frac{ \Gamma \left(e +n-b -c  \right) \Gamma \left(1+b-c  \right)}{\Gamma \left(1-c \right) \Gamma \left(e +n-c  \right)}\;
_{3}F_{2}  \left(\right.\overset{ \displaystyle {n,b,e-a~}}{_ {\displaystyle {e+n-c,1+b-a~}}}|1\left.\right)
\label{Cx4}\\
&=\frac{ \Gamma \left(e +n-b -c  \right) \Gamma \left(1+b-c  \right) \Gamma \left(1+b-a  \right)}{\Gamma \left(b \right) \Gamma \left(1+n-c  \right) \Gamma \left(1+e-a  -c \right)}\;
_{3}F_{2}  \left(\right.\overset{ \displaystyle {e+n-b-c,1-c,1-a~}}{_ {\displaystyle {1+n-c,1+e-a-c~}}}|1\left.\right)
\label{Cx5}\\
&=\frac{ \Gamma \left(e +n -b -c \right) \Gamma \left(1 +b-a \right)}{\Gamma \left(e -a \right) \Gamma \left(1+n-c  \right)}\;
_{3}F_{2}  \left(\right.\overset{ \displaystyle {a +n -c,1-c ,1+b-e~}}{_ {\displaystyle {1+n -c ,1 +b-c~}}}|1\left.\right)
\label{Cx6}\\
&=\frac{ \Gamma \left(e +n-b -c  \right) \Gamma \left(1+b-c  \right)}{\Gamma \left(e -c \right) \Gamma \left(1 +n-c \right)}\;
_{3}F_{2}  \left(\right.\overset{ \displaystyle {n ,1-a ,1+b-e ~}}{_ {\displaystyle {1+n-c  ,1+b-a~}}}|1\left.\right)
\label{Cx7}\\
&=\frac{ \Gamma \left( e +n-b -c \right) \Gamma \left(1+b-a  \right)}{\Gamma \left(n \right) \Gamma \left(1+e-a  -c \right)}\;
_{3}F_{2}  \left(\right.\overset{ \displaystyle {e -c ,1-c ,1+b-a -n~}}{_ {\displaystyle {1+e-a  -c ,1+b-c ~}}}|1\left.\right)
\label{Cx8}\\
&=\frac{ \Gamma \left(e +n-b -c  \right) \Gamma \left(1 +b-c \right)}{\Gamma \left(a +n-c  \right) \Gamma \left(1+e-a  -c \right)}\;
_{3}F_{2}  \left(\right.\overset{ \displaystyle {e -a ,1-a ,1 +b-a -n~}}{_ {\displaystyle {1+e-a  -c ,1 +b-a ~}}}|1\left.\right)\,,
\label{Cx9}
\end{align}
and it is clear that none of these reduce to a finite sum for arbitrary values of  the parameters. However, if one is willing to pursue a more relaxed goal, consider that the case \eqref{Cx8} appears propitious, since it does reduce to a terminating sum if we let $c= m$, a positive integer. However, as noted, the factor $\Gamma(1-c)$ appearing in \eqref{HAll2} requires that a complicated limit $c\rightarrow m$ must be evaluated, and this yields some surprising identities. First of all, it is convenient to write, for $m>0$,
\begin{align} 
{}_{3}^{}{\moversetsp{}{\mundersetsp{}{F_{2}^{}}}}\left(e -c ,1-c ,1-a -n +b ;1-a +e -c ,1-c +b ;1\right) = \moverset{m -1}{\munderset{k =0}{\sum}}\{\cdots\}
+\moverset{\infty}{\munderset{k =m}{\sum}}\{\cdots\}
\label{H1}
\end{align}
where the contents of $\{\cdots\}$ are simply the corresponding symbols appearing in \eqref{3F2Def}. The evaluation of the limit $c\rightarrow m$ is lengthy, but eventually yields two significant results:
\begin{itemize}
\item{}
First, a term corresponding to $(c-m)^{-1}$ appears in the series expansion, and since it is clear that the series $_{3}F_{2}(a,b,m;e,a+n;1)$ is convergent for at least some range of its parameters, it must be that the coefficient of that term vanishes identically. When isolated, this yields the identity
\begin{align} \nonumber
\moverset{m -1}{\munderset{k =0}{\sum}}&\frac{\Gamma \left(e -m +k \right) \left(-1\right)^{k}\,\Gamma \left(1-a -n +b +k \right)}{\Gamma \left(m -k \right) \Gamma \left(k +1\right) \Gamma \left(1-a +e -m +k \right) \Gamma \left(1-m +b +k \right)}\\ \nonumber
 &= \,
-\frac{ \left(-1\right)^{n}\,\Gamma \left(e-m  \right) \Gamma \left(e -b \right) \Gamma \left(n \right) \sin \left(\pi  \left(a -e \right)\right) \csc \left(\pi  \left(a -b \right)\right)}{\Gamma \left(e +n-b -m  \right) \Gamma \left(m \right) \Gamma \left(b \right)}\\
&\times
\moverset{n -1}{\munderset{k =0}{\sum}}\frac{\left(-1\right)^{k}\,\Gamma \left(a +k \right) \Gamma \left(1-e +a +k \right)}{\Gamma \left(n -k \right) \Gamma \left(k +1\right) \Gamma \left(1-b +a +k \right) \Gamma \left(1-m +a +k \right)},
\label{Hid}
\end{align}
that is, two equal terminating sums possess differing numbers of terms. For more on \eqref{Hid} see Section \ref{sec:Hidb}.
\item{}
The evaluation of the leading term $\mathrm{O}((c-m)^0)$ in the series expansion is lengthy and requires significant computer manipulation to reduce it to a simplified form, among which, the identity
\begin{equation}
\munderset{x \rightarrow k}{\mathrm{lim}}\frac{\psi \left(-x \right)}{\Gamma \left(-x \right)}
 = \left(-1\right)^{k +1}\,\Gamma \left(k +1\right)
\label{Limx}
\end{equation}
\end{itemize}
is required during the evaluation of the second term on the right-hand side of \eqref{H1}. The final result produces the identity
\begin{align} \nonumber
{}_{3}^{}{\moversetsp{}{\mundersetsp{}{F_{2}^{}}}}&\left(a ,b ,m ;e ,a +n ;1\right)=\frac{\Gamma \left(a +n \right) \Gamma \left(e \right) \sin \left(\pi  \left(a -e \right)\right) }{\Gamma \left(a \right) \Gamma \left(m \right) \Gamma \left(b \right) \sin \left(\pi  \left(a -b \right)\right)}\\ \nonumber
& \times 
\moverset{n -1}{\munderset{k =0}{\sum}}\frac{\left(-1\right)^{m +k} \left(\psi \left(1-m +a +k \right)-\psi \left(m \right)\right) \Gamma \left(1 +a +k-e \right) \Gamma \left(a +k \right)}{\Gamma \left(1-m +a +k \right) \Gamma \left(1+a +k-b  \right) \Gamma \left(n -k \right) \Gamma \left(k +1\right)}\\ \nonumber
&-\frac{\Gamma \left(a +n \right) \Gamma \left(e \right) \Gamma \left(-b -m +e +n \right)}{\Gamma \left(a \right) \Gamma \left( e-m \right) \Gamma \left(e -b \right) \Gamma \left(n \right)}\\ \nonumber
&\times
 \moverset{m -1}{\munderset{k =0}{\sum}}\frac{\left(-1\right)^{n +m +k}\,\Gamma \left(e -m +k \right) \Gamma \left(1-a -n +b +k \right)}{\Gamma \left(k +1\right) \Gamma \left(m -k \right) \Gamma \left(1-m +b +k \right) \Gamma \left(1-a +e -m +k \right)}\\ \nonumber
&\hspace{20pt}\times \left(\psi \left(m -k \right)-\psi \left(1+e+k-a  -m  \right)-\psi \left(1 +b +k-m \right)\right.\\ \nonumber & \left. \hspace{50pt}+\psi \left(e+k -m  \right)-\psi \left(e-m  \right)+\psi \left( e +n-b -m \right)\right)\\ \nonumber
&+ \left(-1\right)^{n}\,\frac{\Gamma \left(a +n \right) \Gamma \left(1+b +m-a -n  \right) \Gamma \left(e +n-b -m  \right) \Gamma \left(e \right)^{2}}{\Gamma \left(1+e-a  \right) \Gamma \left(1+m \right) \Gamma \left(e -m \right) \Gamma \left(e -b \right) \Gamma \left(n \right) \Gamma \left(1+b \right) \Gamma \left(a \right)}\\
&\hspace{20pt}\times {}_{4}^{}{\moversetsp{}{\mundersetsp{}{F_{3}^{}}}}\left(1,1,e ,1-a -n +b +m ;1+m ,1+b ,1-a +e ;1\right)
\label{HAll2Ha}
\end{align}
where the last term $_{4}F_{3}(1)$, arises from the application of \eqref{Limx} as discussed. 

Although 3F2(1) hypergeometric series of the form $_{3}F_{2}(1,1,a;1+m,e,f;1)$ are known to be reducible to summable series (see \eqref{Hnew}, \eqref{LukeGen} below and \cite{Luke}), I am not aware that a similar property is known for the 4F3(1) generalization appearing in \eqref{HAll2Ha}. However, the particular 4F3(1) appearing in \eqref{HAll2Ha} is $n-$balanced and for $n=1$, a general, (but unsourced and somewhat mysterious) closed expression can be found in the literature (see \cite{Wolf4F3}). In Section  \ref{sec:LimSpec}, a few special case identities will be explored.

\subsubsection{Applied contiguity} \label{sec:ApplContig}
Somewhat after the previous investigations were performed, I became aware that, in 2021, Chen had attacked the Miller/Paris challenge problem using partial fraction decomposition of \eqref{3F2Def} to obtain  the contiguity relation \cite[Lemma 1]{Chensym}
\begin{align} \nonumber
{}_{3}^{}{\moversetsp{}{\mundersetsp{}{F_{2}^{}}}}\left(a ,b ,c ;a +1,e;1\right)
& = 
\frac{\left(a -c +1\right) \left(a -b +1\right) }{\left(1+a -e \right) \left(a +1\right)}{}_{3}^{}{\moversetsp{}{\mundersetsp{}{F_{2}^{}}}}\left(a+1,b ,c ;a +2,e;1\right)\\
&-\frac{\Gamma \left(e \right) \Gamma \left(e -b +1-c \right)}{\left(1+a -e \right) \Gamma \left(e -c \right) \Gamma \left(e -b \right)}\,,
\label{CL1}
\end{align} 
reproducing a known identity previously derived by Krattenthaler and Rivoal \cite[Eq. (12.5)]{KrattRao} and cited by Chu and Wang  \cite[Proposition 20]{ChuWang2007}. Chen then extended \eqref{CL1} using induction \cite[Theorem 1]{Chensym} to read
\begin{align}\nonumber
&{}_{3}^{}{\moversetsp{}{\mundersetsp{}{F_{2}^{}}}}\left(a +m ,b ,c ;a +m+1,e +1;1\right)
 =
\frac{\Gamma \left(a +m +1\right) \Gamma \left(a +1-e +m \right) \Gamma \left(a -c +1\right) \Gamma \left(a -b +1\right) 
}{\Gamma \left(a +1\right) \Gamma \left(1+a -e \right) \Gamma \left(a +1-c +m \right) \Gamma \left(m +a +1-b \right)}\\ \nonumber
&\times\left({}_{3}^{}{\moversetsp{}{\mundersetsp{}{F_{2}^{}}}}\left(a ,b ,c ;e ,a +1;1\right)
+\frac{\Gamma \left(e \right) \Gamma \left(e -b +1-c \right) \Gamma \left(a +1\right) \Gamma \left(1+a -e \right) }{\Gamma \left(a -c +1\right) \Gamma \left(a -b +1\right) \Gamma \left(e -b \right) \Gamma \left(e -c \right)}
\right. \\ &\left . \hspace{5cm}\times\moverset{m -1}{\munderset{k =0}{\sum}}\frac{\Gamma \left(a +1-c +k \right) \Gamma \left(a +1-b +k \right)}{\Gamma \left(a +k +1\right) \Gamma \left(a +2-e +k \right)}\right)\,,
\label{Thm1}
\end{align}
which was then further generalized \cite[Theorem 3]{Chensym} as
\begin{equation}
{}_{3}^{}{\moversetsp{}{\mundersetsp{}{F_{2}^{}}}}\left(a ,b ,c ;a +m,e ;1\right)
 = 
\frac{\Gamma \left(a +m \right) }{\Gamma \left(a \right) }\moverset{m -1}{\munderset{k =0}{\sum}}\frac{ \left(-1\right)^{k}\,{}_{3}^{}{\moversetsp{}{\mundersetsp{}{F_{2}^{}}}}\left(a +k,b ,c  ;a +k +1,e;1\right)}{\Gamma \left(k +1\right) \Gamma \left(m -k \right) \left(a +k \right)}\,.
\label{Thm3}
\end{equation}
Chen also provided an analogue to \eqref{Thm1} for the case ${}_{3}^{}{\moversetsp{}{\mundersetsp{}{F_{2}^{}}}}\left(a-m,b ,c  ;a -m +1,e;1\right)$ as well as an algorithm to evaluate many special cases related to the above, with examples. If one now sets $m=2$ in \eqref{Thm3} to obtain
\begin{align}
{}_{3}^{}{\moversetsp{}{\mundersetsp{}{F_{2}^{}}}}\left(a ,b ,c ;a +2,e;1\right)
 = 
\left(a +1\right) {}_{3}^{}{\moversetsp{}{\mundersetsp{}{F_{2}^{}}}}\left(a ,b ,c ;e ,a +1;1\right)-a\,{}_{3}^{}{\moversetsp{}{\mundersetsp{}{F_{2}^{}}}}\left(a+1,b ,c ;a +2,e
;1\right),
\label{Cm2}
\end{align}
and inverts \eqref{CL1} (equivalent to setting $m=1$ in \eqref{Thm1}) to similarly obtain
\begin{align} \nonumber
{}_{3}^{}{\moversetsp{}{\mundersetsp{}{F_{2}^{}}}}&\left(a+1,b ,c ;a +2,e;1\right)
 = 
\frac{\Gamma \left(e -b +1-c \right) \left(a +1\right) \Gamma \left(e \right)}{\left(a -b +1\right) \left(a -c +1\right) \Gamma \left(e -b \right) \Gamma \left(e -c \right)}\\
&+\frac{ \left(a +1-e \right) \left(a +1\right)}{\left(a -b +1\right) \left(a -c +1\right)}\,{}_{3}^{}{\moversetsp{}{\mundersetsp{}{F_{2}^{}}}}\left(a ,b ,c ;e ,a +1;1\right),
\label{Cm1}
\end{align}
then substitution of \eqref{Cm1} into \eqref{Cm2} yields the identity
\begin{align} \nonumber
{}_{3}^{}{\moversetsp{}{\mundersetsp{}{F_{2}^{}}}}&\left(a ,b ,c ;a +1,e;1\right)
 =\, -\frac{a\,\Gamma \left(e \right) \Gamma \left(e -b +1-c \right)}{\Gamma \left(e -b \right) \left(\left(-e +b -1+c \right) a -\left(b -1\right) \left(c -1\right)\right) \Gamma \left(e -c \right)}\\
&-\frac{\left(a -c +1\right) \left(a -b +1\right)}{\left(\left(-e +b -1+c \right) a -\left(b -1\right) \left(c -1\right)\right) \left(a +1\right)}\, {}_{3}^{}{\moversetsp{}{\mundersetsp{}{F_{2}^{}}}}\left(a ,b ,c ; a +2,e;1\right)\,.
\label{Cm12}
\end{align}
In summary, as Chen notes, all of the above still requires a closed expression for the challenge problem ${}_{3}^{}{\moversetsp{}{\mundersetsp{}{F_{2}^{}}}}\left(a ,b ,c  ;a  +1,e;1\right)$. It must be recognized that \eqref{Cm12} comes close to a solution.

\section{Exotics and Special Cases}\label{sec:Exotic1}
Motivated by the work of Krattenthaler and Rivoal \cite{KrattRiv}, who proved that there exist two-part transformations between 3F2(1) that are not included among those of Thomae, in this Section, we present a selection of two-part identities, ask if they satisfy the Thomae relations and study the consequences. By way of example, Krattenthaler and Rivoal began by providing two examples of two-part identities that do not satisfy the Thomae transformations (see \eqref{HidbKR} and \eqref{Eq1p2}) by employing contiguity properties. They then went on to prove the existence of  infinite families of such exceptions to the Thomae relations in the case of non-negative integer parameters and labelled them ``exotic". 

In the first subsection following, we carry Krattenthaler and Rivoal's work further by testing whether comparison among the KRGT identities satisfy the Thomae transformations.  In the remaining subsections we question whether new identities obtained herein satisfy the Thomae transformations, and study the consequences. We also study the role of both the Thomae transformation group (Appendix \ref{sec:nonTerm}) and that of Rao et. al. (RJRJR - cf. Appendix \ref{sec:termCase}) in the case of terminating series, where a top parameter equals $-n$. 

\subsection{Karlsson versus Rosengren versus Gasper versus Thomae}\label{sec:versus}

To begin, we note that the transformation of Karlsson \eqref{Karl8} is based on contour integration, whereas that of Rosengren \eqref{RosenG} is based on a reduction of a multivariate $q$-series, neither of which coincide with the underlying basis of the Thomae derivation, which is based on symmetric permutations among hypergeometric parameters (see Theorem \ref{sec:T21}).

Written in hypergeometric form (with redefined parameters $e:=1-e+a+b$), by comparing the right-hand side of Karlsson's identity \eqref{Karl8} with that of Rosengren \eqref{RosenG}, we find
\begin{align} \nonumber
{}_{3}^{}{\moversetsp{}{\mundersetsp{}{F_{2}^{}}}}&\left(-n,a ,b ;c ,e ;1\right)
 = 
\frac{\left(-1\right)^{n}\, \Gamma \left(c +e -b -a +n \right) \Gamma \left(1-a\right) \Gamma \left(c \right) \Gamma \left(e \right)}{\Gamma \left(n +e \right) \Gamma \left(c -b +e -a \right) \Gamma \left(1-a -n \right) \Gamma \left(c +n \right)}\\
&\times{}_{3}^{}{\moversetsp{}{\mundersetsp{}{F_{2}^{}}}}\left(-n ,c -a ,e -a ;c -b +e -a ,1-a -n ;1\right)\,.
\label{RK1a}
\end{align}
It is fairly straightforward to establish by a computerized search that neither side of \eqref{RK1a} can be reached from the other by employing any of the Thomae relations from Appendix \ref{sec:nonTerm}, so the transformation \eqref{RK1a} provides an example of an ``exotic" transformation between terminating 3F2(1) with non-integral parameters. This is consistent with Rosengren's assertion that \eqref{RosenG} would ``provide a bridge between the Karlsson-Minton hierarchy and the Bailey hierarchy". 

However, when tested against the RJRJR transformation set, a computerized search finds that the two sides of \eqref{RK1a} are connected in several different ways. Simply, operating by \eqref{Ra1p} with $e\leftrightarrow c$ and $a \leftrightarrow b $ effectively solves \eqref{RK1a}, demonstrating that the right-hand side can be reached from the left-hand side by the appropriate transformation identity. Considered differently, when operated upon by \eqref{Ra8} the left-hand side of \eqref{RK1a} becomes
\begin{align} \nonumber
{}_{3}^{}{\moversetsp{}{\mundersetsp{}{F_{2}^{}}}}&\left(-n ,a ,b ;c ,e ;1\right)
 = 
\frac{\left(-1\right)^{n}\,\Gamma \left(n +e -a \right) \Gamma \left(n-b +e \right) \Gamma \left(c \right) \Gamma \left(e \right) }{\Gamma \left(e -a \right) \Gamma \left(e-b \right) \Gamma \left(c +n \right) \Gamma \left(n +e \right)}\\
&\times{}_{3}^{}{\moversetsp{}{\mundersetsp{}{F_{2}^{}}}}\left(-n ,1+a +b -e -c -n ,1-e -n ;1+a -e -n ,1+b -e -n ;1\right)
\label{A2}
\end{align}and the right-hand side, under the operation of \eqref{Ra3} satisfies
\begin{align} \nonumber
{}_{3}^{}{\moversetsp{}{\mundersetsp{}{F_{2}^{}}}}&\left(-n ,c -a ,e -a ;c +e -a -b ,1-a -n ;1\right)
 =\\  \nonumber
&\frac{\Gamma \left(e-b  +n \right) \Gamma \left(n +e -a \right) \Gamma \left(c +e -a -b \right) \Gamma \left(1-a -n \right) }{\Gamma \left(e-b  \right) \Gamma \left(e -a \right) \Gamma \left(c +e -a -b +n \right) \Gamma \left(1-a \right)}\\
&\times{}_{3}^{}{\moversetsp{}{\mundersetsp{}{F_{2}^{}}}}\left(-n ,1+a +b -e -c -n ,1-e -n ;1+a -e -n ,1+b -e -n ;1\right)\,.
\label{A1}
\end{align}
Eliminating the common element on  the right-hand sides of \eqref{A2} and \eqref{A1} reproduces \eqref{RK1a}.

In a similar vein, equating the corresponding side of  Gasper's transformation \eqref{Gasper} with that of Rosengren \eqref{RosenG}, after redefinition of variables, yields the identity
\begin{align} \nonumber
{}_{3}^{}{\moversetsp{}{\mundersetsp{}{F_{2}^{}}}}\left(a ,b ,c +n ;e ,c ;1\right)
& = 
\frac{\Gamma \left(1-c +a \right) \Gamma \left(e -n -a -b \right)  \Gamma \left(-c -n +1\right) \Gamma \left(e \right)}{\Gamma \left(e -a \right) \Gamma \left( e-b -n \right) \Gamma \left(1-n -c +a \right) \Gamma \left(1-c \right)}\\
&\times {}_{3}^{}{\moversetsp{}{\mundersetsp{}{F_{2}^{}}}}\left(-n,a  ,e -c -n ;1-n -c +a ,e -b -n ;1\right)
\label{GaRosen}
\end{align}
and, by inspection, we find that this transformation also is not Thomae related. However, by reversing the embedded right-hand side sum, the identity \eqref{GaRosen} simply reduces to Rosengren's original identity \eqref{RosenG} itself, thereby demonstrating the equivalence of the identities of Gasper and Rosengren.
\begin{rem}
Setting $e=b+1$ in \eqref{GaRosen}, as was done originally in Section \ref{sec:History}, leads to the need to evaluate the problematic expression
\begin{align} \nonumber
{}_{3}^{}{\moversetsp{}{\mundersetsp{}{F_{2}^{}}}}&\left(a ,b ,c +n ;b +1,c ;1\right)
 = 
\frac{\Gamma \left(1-c +a \right) \Gamma \left(1-a -n \right) \Gamma \left(1-c -n \right) \Gamma \left(b +1\right)}{\Gamma \left(b +1-a \right) \Gamma \left(1-n \right) \Gamma \left(a -c -n +1\right) \Gamma \left(1-c \right)}\\
 &\times {}_{3}^{}{\moversetsp{}{\mundersetsp{}{F_{2}^{}}}}\left(-n  ,a,1+b-c -n ;a -c -n +1,1-n ;1\right)
\label{GaRoSubs}
\end{align}
where the right-hand side 3F2(1) apparently diverges because the second bottom parameter is a negative integer that exceeds the first top parameter, also a negative integer, compensated for by the presence of the denominator factor $\Gamma(1-n)$. However, first reversing the sum order gives \eqref{RosenG},
in which case setting $e=b+1$ immediately yields Karlsson's identity \eqref{K10} and therefore a simple way to interpret
\begin{equation}
\frac{1}{\Gamma \left(1-n \right)}\,{}_{3}^{}{\moversetsp{}{\mundersetsp{}{F_{2}^{}}}}\left(-n ,a ,b ;c ,1-n ;1\right)
 = 
-\frac{\Gamma \left(b +n \right) \Gamma \left(1-a \right) \Gamma \left(c \right)}{\Gamma \left(b \right) \Gamma \left(c +n \right) \Gamma \left(1-a -n \right)}
\label{GaByG}
\end{equation}
after redefining variables $b:=b-1+c+n~,c:=-c+a-n+1$.
\end{rem}
A similar comparison between the transformations of Gasper \eqref{Gasper} and Karlsson \eqref{Karl8}, yields
\begin{align} \nonumber
{}_{3}^{}{\moversetsp{}{\mundersetsp{}{F_{2}^{}}}}\left(-n,a ,b  ;c ,e ;1\right)
 &= 
\frac{\Gamma \left(1-a \right) \Gamma \left(n +c -b \right) \Gamma \left(c \right) \Gamma \left(1-e +b \right)}{\Gamma \left(c -b \right) \Gamma \left(b +1-a \right) \Gamma \left(c +n \right) \Gamma \left(1-e \right)} \\
&\times {}_{3}^{}{\moversetsp{}{\mundersetsp{}{F_{2}^{}}}}\left(b ,e -a ,1+b -c ;b +1-a ,1+b -c -n ;1\right),
\label{GaspKarl}
\end{align}
and again we find that neither side can be reached from the other by any of the Thomae transformations, although they can be reached by several variants of the RJRJR transformation set.

Continuing on  the same course, we compare each of the KRG transformations in turn with the Thomae relation \eqref{X1} to find respectively
\begin{align}
{}_{3}^{}{\moversetsp{}{\mundersetsp{}{F_{2}^{}}}}\left(-n,a ,b ;c ,e ;1\right)
 &= 
\frac{ \Gamma \left(1-e -n \right) \Gamma \left(1-e +a \right)}{\Gamma \left(1+a -e -n \right) \Gamma \left(1-e \right)}{}_{3}^{}{\moversetsp{}{\mundersetsp{}{F_{2}^{}}}}\left(-n ,c -b ,a ;c,1+a -e -n  ;1\right)
\label{X5}
\\\nonumber
&=\frac{\Gamma \left(1-e +b \right) \Gamma \left(c \right) \Gamma \left(1-a \right)  \Gamma \left(e \right)}{\Gamma \left(e+n \right) \Gamma \left(c +n \right) \Gamma \left(1-a -n \right) \Gamma \left(b +1-e -n \right)}\\
&\times{}_{3}^{}{\moversetsp{}{\mundersetsp{}{F_{2}^{}}}}\left(-n ,c -a ,1-e -n ;b +1-e -n ,1-a -n ;1\right)
\label{X7}
\\ \nonumber
&=\frac{\Gamma \left(1-a \right) \Gamma \left(b +n \right) \Gamma \left(c \right) \Gamma \left(c +e -b -a +n \right) \Gamma \left(e \right)}{\Gamma \left(b \right) \Gamma \left(c -b +1-a \right) \Gamma \left(n +e \right) \Gamma \left(c +n \right) \Gamma \left(e -a \right)}\\
& \times{}_{3}^{}{\moversetsp{}{\mundersetsp{}{F_{2}^{}}}}\left(1-e -n ,c -b ,1-b ;c -b +1-a ,1-b -n ;1\right)
\label{X6}
\end{align}
and discover that none of the above satisfy a Thomae relation, although both \eqref{X5} and \eqref{X7} are special cases connected by the RJRJR transformation set as exemplified by \eqref{A2} and \eqref{A1}. 

\begin{rem}
Solving \eqref{X6} and redefining the variables, ($e:=1-e-n,c:=c+b,b:=1-b,a:=1-a+c,b:=b+n,e\leftrightarrow a,b \leftrightarrow c$) yields yet another transformation identity
\begin{align} \nonumber
{}_{3}^{}{\moversetsp{}{\mundersetsp{}{F_{2}^{}}}}&\left(a ,b ,c +n ;c,e ;1\right)
 = 
\frac{ \Gamma \left(1-c -n \right) \Gamma \left(e \right) \Gamma \left(1-a \right) \Gamma \left(1-c +b \right) \Gamma \left(e -n -a -b \right)}{\Gamma \left(e -b \right) \Gamma \left(1-c \right) \Gamma \left(1-c -n +b\right) \Gamma \left(e -a \right) \Gamma \left(1-a -n\right)}\\
&\times{}_{3}^{}{\moversetsp{}{\mundersetsp{}{F_{2}^{}}}}\left(-n ,1-c -n ,1-e +b ;1-a -n ,1-c -n +b ;1\right),
\label{X6a}
\end{align}
reducing to Karlsson's identity \eqref{K10} if $e=b+1$. Further, setting $b=c - a - 1 + e + n$  and redefining the variables $c:=c+a+1-e-n$ reduces \eqref{X7} to the Saalsch\"{u}tz theorem \eqref{Saal}.
\end{rem}


\subsection{Beyond Thomae - order reversal} \label{sec:Beyond}

In addition to the alternative derivations presented in the previous section, here we recognize that the various entities being considered are terminating sums and therefore order-reversible. This introduces a new degree of freedom.  For example, by reversing the right-hand side of Karlssons identity \eqref{Karl8}, we discover that
\begin{align} \nonumber
{}_{3}^{}{\moversetsp{}{\mundersetsp{}{F_{2}^{}}}}&\left(-n, a ,b  ;c ,1-e +a +b ;1\right)
 = 
\frac{\left(-1\right)^{n}\,\Gamma \left(c \right) \Gamma \left(1-e +a +b \right) \Gamma \left(a +n \right) \Gamma \left(b +n \right) }{\Gamma \left(a \right) \Gamma \left(b \right) \Gamma \left(c +n \right) \Gamma \left(1-e +a +b +n \right)}\\
&\times{}_{3}^{}{\moversetsp{}{\mundersetsp{}{F_{2}^{}}}}\left(-n ,1-c -n ,e -n -a -b ;1-n -a ,1-b -n ;1\right)\,.
\label{RkRev}
\end{align}
This yields an analogue of \eqref{Karl8} after substitution of \eqref{RkRev} into \eqref{Karl8}. Equating the right-hand side of that analogue (i.e. \eqref{RkRev}) with the right-hand side of Rosengren's identity \eqref{RosenG} then yields a transformation
\begin{align} \nonumber
{}_{3}^{}{\moversetsp{}{\mundersetsp{}{F_{2}^{}}}}&\left(-n ,a ,b ;c ,e ;1\right)
 = 
\frac{ \left(-1\right)^{n}\Gamma \left( c +e +n -a-b \right)  \Gamma \left(1-e -n\right)}{\Gamma \left(c +e -a -b \right) \Gamma \left(1-e \right)}
\\
&\times{}_{3}^{}{\moversetsp{}{\mundersetsp{}{F_{2}^{}}}}\left(-n ,c-a ,c-b ;c +e-a -b  ,c ;1\right)
\label{RkRosen}
\end{align}
where the variables (from the left) have been redefined ($c:=1-c-n,e:=e+a+b+n,a:=1-a-n,b:=1-b-n,c\leftrightarrow a,e\leftrightarrow b$). A simple inspection then finds that the two sides of \eqref{RkRosen} are related by the set of Thomae transforms. In this sense, \eqref{RkRosen} demonstrates that the identities of Karlsson and Rosengren are related by the Thomae transformations plus reversal of one of the sums, although they are not related by the Thomae transformations themselves. Similarly, reversing the right-hand side of the Thomae transform \eqref{X1} produces
\begin{align} \nonumber
{}_{3}^{}{\moversetsp{}{\mundersetsp{}{F_{2}^{}}}}&\left(-n ,a,c-b  ;c,e-b -n  ;1\right)
 = 
\frac{\left(-1\right)^{n}\,\Gamma \left(c \right) \Gamma \left(c +n -b \right) \Gamma \left(a +n \right) \Gamma \left(e-b  -n \right)}
{\Gamma \left(c-b  \right) \Gamma \left(a \right) \Gamma \left(e -b \right) \Gamma \left(c +n \right)}\\
&\times {}_{3}^{}{\moversetsp{}{\mundersetsp{}{F_{2}^{}}}}\left(-n ,1-e +b ,1-c -n ;1-n -a ,1+b -c -n ;1\right)\,.
\label{Thm1Rev}
\end{align}
Substitution of \eqref{Thm1Rev} into the right-hand side of \eqref{X1} and comparison with the right-hand sides of Gasper's identity \eqref{Gasper} simply reproduces \eqref{Gasper} itself, after redefining the variables. This demonstrates that Gasper's identity is simply the reversal of the right-hand side of the Thomae transformation \eqref{X1}, which itself is simply the basic Thomae transformation \eqref{Thom4} with redefined variables.

The use of reversal transforms also introduces a variety of different ways to obtain closed form identities, much as was used by Karlsson et.al. (Section \ref{sec:History}). For example, reversing the Thomae transformation \eqref{X4} gives
\begin{align} \nonumber
{}_{3}^{}{\moversetsp{}{\mundersetsp{}{F_{2}^{}}}}&\left(a ,b ,c +n ;c ,e ;1\right)
 = 
\frac{\left(-1\right)^{n}\Gamma \left(e -a -b \right) \Gamma \left(e-c  \right) \Gamma \left(e \right) \Gamma \left(c \right)}{\Gamma \left(e-c -n  \right) \Gamma \left(e -b \right) \Gamma \left(e -a \right) \Gamma \left(c +n \right)}\\
&\times{}_{3}^{}{\moversetsp{}{\mundersetsp{}{F_{2}^{}}}}\left(-n ,1-e +b ,1-e+a ;c +1-e ,1-e +a +b ;1\right) 
\label{Thm2Rev}
\end{align} 
and, setting $e=b+1+m$, then $a=b+m$ with $b+n<1$, yields the closed form
\begin{equation}
{}_{3}^{}{\moversetsp{}{\mundersetsp{}{F_{2}^{}}}}\left(b  ,b +m ,c +n;c ,b +1+m ;1\right)
 = 
\frac{\left(-1\right)^{n}\Gamma \left(1-b \right) \Gamma \left(1 +m-c +b \right) \,\Gamma \left(b +1+m \right) \Gamma \left(c \right)}{\Gamma \left(b +1+m -c -n \right) \Gamma \left(1+m \right) \Gamma \left(c +n \right)}\,,
\label{T8X}
\end{equation}
a variation of previous identities (e.g.\eqref{X1b}).
\begin{rem}
Note that the right-hand side of \eqref{T8X} is the analytic continuation of the left-hand side if $b>1-n$ where the sum diverges. This identity also happens  to coincide with a Thomae progeny of the  Saalsch\"{u}tz theorem \eqref{Saal}.
\end{rem}
\subsection{Other 2-part identities as Thomae transformations} \label{sec:THM}
\subsubsection{Based on \eqref{Hnew2}} \label{sec:Hnew2}
In the derivation of \eqref{Hn4b}, where \eqref{Hnew2} is evaluated in the limit $c\rightarrow b+m$ (see Section \ref{sec:LimCase} below), a divergent term of order $(c-(b+m))^{-1}$ arises in the expansion. This coefficient must vanish. The reader is invited to verify that this leads to an apparently new transform identity
\begin{align} \nonumber
{}_{3}^{}{\moversetsp{}{\mundersetsp{}{F_{2}^{}}}}\left(-n ,m ,a ;1,c ;1\right)
 = {}_{3}^{}{\moversetsp{}{\mundersetsp{}{F_{2}^{}}}}\left(-n ,1-m ,1-c -n ;1-n -m,a+1 -n -c  ;1\right)\\
\times\left(-1\right)^{m+1}\,\frac{\Gamma \left(m -c +a -n \right) \Gamma \left(n +m \right)  \Gamma \left(c \right) \Gamma \left( 1-a +c +n-m \right)}{\Gamma \left(c-a  \right) \Gamma \left(n +1\right) \Gamma \left(c +n \right) \Gamma \left(m \right) \Gamma \left(a -n +1-c \right)}
\label{THM}
\end{align} 
where the variables $b:=b+a+n-m-1,b:=c,a:=1-a,n:=n+1$ have been redefined to put \eqref{THM} into canonical form. This raises the question whether the two sides of \eqref{THM} are related by one of the Thomae transformations; a simple survey shows that they are not so-related, although they do satisfy several of the RJRJR transformation identities. In addition, {\bf by numerical experimentation}, it was discovered that \eqref{THM} is numerically valid in the range in which the modified left-hand side converges, when the variable $n$ is extended to become continuous. With $n:=-b$, and $c>1-a-b-m$, this eventually yields the identity
\begin{align} \nonumber
{}_{3}^{}{\moversetsp{}{\mundersetsp{}{F_{2}^{}}}}&\left(a ,b ,m ;1,c ;1\right)
= {}_{3}^{}{\moversetsp{}{\mundersetsp{}{F_{2}^{}}}}\left(1-m ,b ,1-c +b ;1+b -m ,a +1+b -c ;1\right)\\&
\times\left(-1\right)^{m+1}\,\frac{\Gamma \left(m -c +a +b \right) \Gamma \left(m-b  \right) \Gamma \left(c \right) \Gamma \left(1-a +c -b -m \right)}{\Gamma \left(m \right) \Gamma \left(1-b \right) \Gamma \left(a +1+b -c \right) \Gamma \left(c-a  \right) \Gamma \left(c -b \right)}\,,
\label{THM3}
\end{align}
which does not appear to coincide with any of the corresponding identities discussed in Section \ref{sec:History}.

\subsubsection{Based on \eqref{CxA}}
In a similar vein, we consider \eqref{CxA} below and discover that the two sides are not connected by any of the Thomae transformations. In general, both sides of \eqref{CxA} are terminating sums of the form \eqref{NegTop} and, with one exception, all the transformations generated by the application of the Thomae transformations to either of the two terminating sums composing \eqref{CxA} are uninteresting. The exception arises when the transformation \eqref{Thom4} is applied to \eqref{CxA} resulting in the identity
\begin{equation}
{}_{3}^{}{\moversetsp{}{\mundersetsp{}{F_{2}^{}}}}\left(1-n ,1-j, m ;1,p +1-j -n ;1\right)
 = 
\frac{ \Gamma \left(p -m \right) \Gamma \left(p +1-j -n \right)}{\Gamma \left(p -m +1-j -n \right) \Gamma \left(p \right)}\;{}_{3}^{}{\moversetsp{}{\mundersetsp{}{F_{2}^{}}}}\left(j ,n ,m ;p ,1;1\right)\,,
\label{Cx4A}
\end{equation} 
where we have redefined $p:=p-m$, in which case convergence of \eqref{Cx4A} is valid for $p>m$. Inverting \eqref{Cx4A} yields 
\begin{equation}
{}_{3}^{}{\moversetsp{}{\mundersetsp{}{F_{2}^{}}}}\left(j ,m ,n ;p ,1;1\right)
 = 
\frac{\Gamma \left(p \right) \Gamma \left(p -m +1-j -n \right) }{\Gamma \left(p -m \right) \Gamma \left(p +1-j -n \right)}{}_{3}^{}{\moversetsp{}{\mundersetsp{}{F_{2}^{}}}}\left(1-n ,1-j,m ;1,p +1-j -n ;1\right),
\label{Cx1(4)}
\end{equation}
thereby transforming an infinite series where all  parameters are positive integers, into a finite summation with either $n$ or $j$ terms. 
\begin{rem}
The left-hand side of \eqref{Cx1(4)} is a special case of the general summation forms studied by Krattenthaler and Rivoal \cite{KrattRiv} and Chu and Wang \cite[Section 4]{ChuWang2007}.
\end{rem}

\begin{rem}

The left-hand series \eqref{Cx1(4)} converges only if $p+1-j-m-n>0$ in which case it is impossible for $m>p$ or $j+n>p+1$, given the defined ranges of the integers $j,m,n,p$ and so the right-hand side always exists if the left-hand side converges. Written in sum form, \eqref{Cx1(4)} becomes
\begin{align} \nonumber
{}_{3}^{}{\moversetsp{}{\mundersetsp{}{F_{2}^{}}}}\left(j ,m ,n ;p ,1;1\right)
 =& 
\frac{\Gamma \left(n \right) \Gamma \left(j \right)\Gamma \left(p \right) \Gamma \left(p -m +1-j -n \right) }{\Gamma \left(m \right) \Gamma \left(p -m \right)}\\
&\times\moverset{N}{\munderset{k =0}{\sum}}\frac{\Gamma \left(m +k \right)}{\Gamma \left(n -k \right) \Gamma \left(j -k \right) \Gamma \left(k +1\right)^{2}\,\Gamma \left(p +1-j -n +k \right)}\,,
\label{CX4}
\end{align}
where $N=\min \left(j-1 , n -1\right)$. Since all variables are independent, except for the constraint $p>m$, it is intriguing to impose $p\leq m$, leading to an interesting quandry, since  we then have
\begin{align} \nonumber
{}_{3}^{}{\moversetsp{}{\mundersetsp{}{F_{2}^{}}}}\left(j ,m ,n ;1,p ;1\right)
 =& \left(-1\right)^{j+n+1} 
\frac{\Gamma \left(n \right) \Gamma \left(j \right)  \Gamma \left(p \right) \Gamma \left(1-p +m \right)}{\Gamma \left(m \right) \Gamma \left(m -p +j +n \right)}\\
&\times\moverset{N}{\munderset{k =0}{\sum}}\frac{\Gamma \left(m +k \right)}{\Gamma \left(n -k \right) \Gamma \left(j -k \right) \Gamma \left(k +1\right)^{2}\,\Gamma \left(p +1-j -n +k \right)},
\label{HX4}
\end{align}
where the left-hand side is divergent since $p-j-m-n-1<0$, and the right-hand side is not. {\bf If analytic continuation were valid for discrete variables}, it is tempting to speculate that the right-hand side of \eqref{HX4} would represent the analytic continuation of the left-hand side in the variables $p$ and $m$, and, by
replacing integers by reals (i.e. $m:=a,p:=b,j:=m$ ordered from the left), the resulting identity 
\begin{align} \nonumber
{}_{3}^{}{\moversetsp{}{\mundersetsp{}{F_{2}^{}}}}&\left(m ,n ,a ;b ,1;1\right)
 = 
 \left(-1\right)^{1+m+n}\frac{\Gamma \left(n \right) \Gamma \left(m \right)  \Gamma \left(b \right) \Gamma \left(1+a -b \right)}{\Gamma \left(a \right) \Gamma \left(a -b +m +n \right)}\\
 &\times \moverset{N}{\munderset{k =0}{\sum}}\frac{\Gamma \left(a +k \right)}{\Gamma \left(n -k \right) \Gamma \left(m -k \right) \Gamma \left(k +1\right)^{2}\,\Gamma \left(b -n -m +1+k \right)}
\label{Cxab}
\end{align}
where $N=\min(n-1,m-1)$ would be true, where we reiterate that $a,b\in\mathfrak{C}$. In fact, we experimentally find that \eqref{Cxab} {\bf is} numerically valid if $b+1-a-m-n>0$. Since $a,b$ are now continuous variables, the right-hand side of \eqref{Cxab} indeed represents the analytic continuation of the left-hand side in the case that $\Re(a)\geq\Re(b)$ (originally $m\geq p$). In fact, since the left-hand side of \eqref{Cxab} only converges if $b+1-a-m-n>0$, the right-hand side of \eqref{Cxab} is the analytic continuation of the left-hand side if the inequality fails and conjecture morphs into (numerical) truth.

In addition, other representations exist that lend credibility to the above. Letting $b=m$ and $c:=b$ in database entry \eqref{data482}, -- see Appendix \ref{sec:DataNew} -- will yield a representation equivalent to \eqref{Cxab},
noting however, that database entry \eqref{data482} itself is a Thomae progeny of database entry 251, which was experimentally determined. However, all this speculation and experimentation then leads to the apparent paradox that, if $b=a$, the left-hand side of \eqref{Cxab} is divergent and the right-hand side is not. The resolution arises from the fact that the divergent sums symbolized by ${}_{2}^{}{\moversetsp{}{\mundersetsp{}{F_{1}^{}}}}\left(m ,n ;1;1\right)$ and $\munderset{b \rightarrow a}{\mathrm{lim}}{}_{3}^{}{\moversetsp{}{\mundersetsp{}{F_{2}^{}}}}\left(m ,n ,a ;b ,1;1\right)$ mean two different things. See the previous discussion of the work of Karp and Prilepkina \eqref{T2a} for another such caveat.

\end{rem}


Continuing, the difference between the top parameters $m>1$ and $n>1$ and the second bottom parameter (unity) of \eqref{Cxab} is always a positive integer, so the left-hand side can be transformed into a terminating sum through the use of any of the KBGT transformations. For example, after applying the Karlsson transformation \eqref{Karl8}, we obtain 
\begin{align} 
{}_{3}^{}{\moversetsp{}{\mundersetsp{}{F_{2}^{}}}}\left(m ,n ,a ;b,1 ;1\right)
 &= 
\frac{\Gamma \left(n \right) \Gamma \left(b \right)}{\Gamma \left(m \right) \Gamma \left(a \right) \Gamma \left(b -m \right) \Gamma \left(b -a \right)}
 \moverset{n -1}{\munderset{k =0}{\sum}}\frac{\Gamma \left(m +k \right) \Gamma \left(a +k \right) \Gamma \left(b -k -m -a \right)}{\Gamma \left(k +1\right)^{2}\,\Gamma \left(n -k \right)}\,.
\label{Ka2}
\end{align}
Due to symmetry, \eqref{Ka2} is also valid under the interchange $n\leftrightarrow m$. By equating the right-hand sides of \eqref{Ka2} under interchange, this symmetry immediately leads to the subsequent identity
\begin{equation}
\frac{\left(-1\right)^{m}\,{}_{3}^{}{\moversetsp{}{\mundersetsp{}{F_{2}^{}}}}\left(a ,m ,1-n ;1,a -b +m +1;1\right)}{\Gamma \left(b -m \right) \Gamma \left(a -b +m +1\right)}
 = 
\frac{\left(-1\right)^{n}\,{}_{3}^{}{\moversetsp{}{\mundersetsp{}{F_{2}^{}}}}\left(a ,n ,1-m ;1,a -b +n +1;1\right)}{\Gamma \left(b -n \right) \Gamma \left(a -b +n +1\right)}\,.
\label{H12h}
\end{equation}
Setting $b=n+1$ with $n\geq m$ produces the sum
\begin{equation}
\moverset{n -1}{\munderset{k =0}{\sum}}\frac{\Gamma \left(m +k \right) \Gamma \left(a +k \right) \left(-1\right)^{k}}{\Gamma \left(k +m +a -n \right) \Gamma \left(k +1\right)^{2}\,\Gamma \left(n -k \right)}
 = \left(-1\right)^{n+1}\,,
\label{SumId}
\end{equation}
which, in hypergeometric notation identifies the terminating sum
\begin{equation}
{}_{3}^{}{\moversetsp{}{\mundersetsp{}{F_{2}^{}}}}\left(1-n,a ,m  ;1,a -n +m ;1\right)
 = 
\frac{\left(-1\right)^{n+1}\,\Gamma \left(a -n +m \right) \Gamma \left(n \right)}{\Gamma \left(m \right) \Gamma \left(a \right)},\hspace{20pt} n\geq m\,,
\label{HidZ}
\end{equation}
a zero-balanced, terminating variant of the Saalsch\"{u}tz theorem \eqref{Saal}.

\subsection{``Of independent interest"}

Extracted from the literature, and ``of independent interest", for $m+1-a-n>0$, Kalmykov and Karp \cite[Theorem 1]{KarpAndKalmy} derive 
\begin{align} \nonumber
{}_{3}^{}{\moversetsp{}{\mundersetsp{}{F_{2}^{}}}}&\left(a,b ,c +n ;c ,b +1+m ;1\right)\\ \nonumber
& = 
\frac{\left(-1\right)^{n}\,\Gamma \left(1-a \right) \Gamma \left(b +1+m \right) \Gamma \left(m +c -a \right) \Gamma \left(1-b \right) \Gamma \left(1-c -n \right) \Gamma \left(c +n -b -m \right) }{\Gamma \left(c -a \right) \Gamma \left(1-b -m \right) \Gamma \left(1+b +m -a \right) \Gamma \left(1-c \right) \Gamma \left(c -b \right) \Gamma \left(1+m \right)}\\
&\times {}_{3}^{}{\moversetsp{}{\mundersetsp{}{F_{2}^{}}}}\left(-m ,1-c -n ,a-b -m  ;1-b -m ,1-m +a -c ;1\right)\,
\label{KarpKalmy}
\end{align}
corresponding to a special  case of the KRGT transformations and again, independently and effectively, sum a non-terminating series. A simple examination shows that the two sides of \eqref{KarpKalmy} are unrelated by any of the Thomae transformations. Setting $a=b+m$ in \eqref{KarpKalmy} immediately yields the closed form
\begin{equation}
{}_{3}^{}{\moversetsp{}{\mundersetsp{}{F_{2}^{}}}}\left(b ,c +n ,b +m ;c ,b +1+m ;1\right)
 = 
\frac{\left(-1\right)^{n}\,\Gamma \left(b +1+m \right) \Gamma \left(1-b \right) \Gamma \left(1-c -n \right) \Gamma \left(c +n -b -m \right)}{\Gamma \left(c -b -m \right) \Gamma \left(1-c \right) \Gamma \left(1+m \right)}\,,
\label{T1Ax}
\end{equation} 
convergent if $b<1-n$, $n\geq0$, and then if $c=b+1$ we find
\begin{equation}
{}_{3}^{}{\moversetsp{}{\mundersetsp{}{F_{2}^{}}}}\left(b ,b +m ,b +1+n ;b +1,b +1+m ;1\right)
 =\, 
-\frac{\Gamma \left(b +1+m \right) b\,\Gamma \left(-b -n \right) \Gamma \left(m \right)}{\Gamma \left(m -n \right) \Gamma \left(1+m \right)}\,.
\label{T1Cx}
\end{equation} 
A form similar to \eqref{T1Cx} can be obtained by setting $c=b+1$ in Gasper's specialized identity \eqref{G81M}, taking care to consider both cases $n<m$ and otherwise, yielding the (quasi closed-form) identity
\begin{align} \nonumber
{}_{3}^{}{\moversetsp{}{\mundersetsp{}{F_{2}^{}}}}&\left(a ,b ,b +1+n ;b +1+m ,b +1;1\right)
 = \frac{\Gamma \left(b +1+m \right) \Gamma \left(1-a \right) \Gamma \left(1+n \right) \Gamma \left(b +1\right)}{\Gamma \left(b +1-a \right) \Gamma \left(1+m \right) \Gamma \left(b +1+n \right)}\\ \nonumber
&+\frac{\Gamma \left(b +1+m \right) \Gamma \left(-a +1\right) \Gamma \left(b +1\right) \left(-1\right)^{n}}{\Gamma \left(b +1+n \right) \Gamma \left(b \right)}\\
&\times\moverset{m}{\munderset{k =1+n}{\sum}}\frac{\left(-1\right)^{k}\,\Gamma \left(b +k \right) \Gamma \left(k \right)}{\Gamma \left(1-k +m \right)\,\Gamma \left(1+k \right) \Gamma \left(1-a +b +k \right) \Gamma \left(k-n \right)}\,.
\label{T2aX}
\end{align}
%
Of particular interest is the observation that the right-hand side of \eqref{KarpKalmy} does not correspond with the corresponding side of either of the Karlsson/Rosengren transformations \eqref{Karl8} and \eqref{RosenG}, so by equating the respective sides in turn, we find, for \eqref{Karl8},
\begin{align} \nonumber
{}_{3}^{}{\moversetsp{}{\mundersetsp{}{F_{2}^{}}}}&\left(-m,a ,b +n  ;b ,c ;1\right)
 = 
{}_{3}^{}{\moversetsp{}{\mundersetsp{}{F_{2}^{}}}}\left(-n ,1-c -m ,b -a ;1-a -n,1+b -c ;1\right)\\
&\times\frac{\left(-1\right)^{m}\,\Gamma \left(1-b -n \right) \Gamma \left(c -b \right) \Gamma \left(1-a \right) \Gamma \left(1-c -m \right) \Gamma \left(c -a +m -n \right)}{\Gamma \left(c -b -n \right) \Gamma \left(1-a -n \right) \Gamma \left(1-c \right) \Gamma \left(c -a \right) \Gamma \left(1-b \right)}\,,
\label{Extra}
\end{align}
and for \eqref{RosenG}
\begin{align} \nonumber
&{}_{3}^{}{\moversetsp{}{\mundersetsp{}{F_{2}^{}}}}\left(-n ,a ,b ;c ,1-c +b +a +m -n ;1\right)
 \\\nonumber
 &= 
\frac{ \Gamma \left(c -b +n -m \right) \Gamma \left(c -b -a -m \right) \Gamma \left(c -a +n -m \right) \Gamma \left(c \right)}{\Gamma \left(c -b -a -m +n \right) \Gamma \left(c-b  -m \right) \Gamma \left(c -a -m \right) \Gamma \left(c +n \right)} \\ 
&\hspace{50pt}\times{}_{3}^{}{\moversetsp{}{\mundersetsp{}{F_{2}^{}}}}\left(-m ,c -b -a -m ,-n ;c -a -m ,-b +c -m ;1\right)
\label{R3}
\end{align}
after redefining the variables.
\begin{rem}

Setting $m=0$ in \eqref{R3}, recovers the Saalsch\"{u}tz theorem \eqref{Saal} itself. Setting $m=1$, reproduces \eqref{XidbTwo} below, a proof of the (unproven) identity cited in \cite[Eq.5]{MathWSaal}, and resolved in \eqref{XidbTwo} below. Otherwise, \eqref{R3} represents a generalization of the Saalsch\"{u}tz theorem with parametric excess equal to $m+1$, equivalent to \eqref{Xidb} below. It also reproduces the case corresponding to $p=1$ derived by Karp and Prilepkina \cite[Eq. 2.12]{Karp&Prep2018}.\footnote{In that identity, the symbol $(-q)_{q}$ should read $(-1)^q\,\Gamma(q+1)$.}
\end{rem}
There are obviously many more permutations available among the KRGT-like identities, including the nine Thomae transformations plus reversal and those of Rao et. al. In addition, Karp and Prilepkina have shown \cite{KandP2020} that in addition to the traditional three-part transformations, of which only \eqref{HAll1} has been used here, there exist other three-part identities that lie outside the ``Bailey universe".

\section{Limiting, Special Cases and Extensions} \label{sec:LimSpec}

In the following subsections, we consider special, extended and limiting cases of the more general identities that have been developed up to this point.

\subsection{Based on \eqref{Hid}} \label{sec:Hidb}
In Section (\ref{sec:Equivs}), the identity \eqref{Hid} was presented. That result can be written in hypergeometric form:
\begin{align} \nonumber
{}_{3}^{}{\moversetsp{}{\mundersetsp{}{F_{2}^{}}}}&\left(1-m , e-m ,1+b-a -n  ;1+b-m  ,1+e -a -m ;1\right)\\ \nonumber
 =\, 
&\frac{\left(-1\right)^{n+m}\Gamma \left(a \right) \Gamma \left(1+a-e  \right) \Gamma \left( b-a \right) \Gamma \left(e -b \right)  \Gamma \left(1+b-m  \right)}{\Gamma \left(1+a -m \right) \Gamma \left(e +n-b -m  \right) \Gamma \left(b \right) \Gamma \left(1+b-a -n  \right) \Gamma \left(a +m-e  \right)}\\
&\times  {}_{3}^{}{\moversetsp{}{\mundersetsp{}{F_{2}^{}}}}\left(1-n, a ,1+a-e  ;1+a-b  ,1 +a-m ;1\right)\,.
\label{HidX}
\end{align}
%
After redefinition of variables (specifically $e:=e+m,b:=b+a+n-1,a:=-a+e+1,m:=m+1,a:=c,e:=a$ ordered from the left), in canonical form, \eqref{HidX} becomes
\begin{align} \nonumber
&{}_{3}^{}{\moversetsp{}{\mundersetsp{}{F_{2}^{}}}}\left(-m,a ,b  ;c ,a+b+n-m -c   ;1\right)\\ \nonumber
& = 
\left(-1\right)^{n+1}\,\frac{\Gamma \left(c \right)\Gamma \left(1-c +a \right) \Gamma \left(b +n -1\right) \Gamma \left(1+m -n -b +c  \right) \Gamma \left(n-m +b -c +a \right) }{\Gamma \left(m +c \right) \Gamma \left(1-c +a -m \right) \Gamma \left(-b +c \right) \Gamma \left(b -c +a +n \right) \Gamma \left(b \right)}\\
&\times  {}_{3}^{}{\moversetsp{}{\mundersetsp{}{F_{2}^{}}}}\left(1-n ,1-c +a ,1-m -c ;2-b -n ,1-c +a -m ;1\right)\,
\label{Xidb}
\end{align}
and we immediately find an $n-$balanced generalization of the Saalsch\"{u}tz theorem \eqref{Saal}, a transformation analogue of \eqref{R3} that reduces to \eqref{Saal} if $n=1$. If $n=2$, we obtain
\begin{align} \nonumber
{}_{3}^{}{\moversetsp{}{\mundersetsp{}{F_{2}^{}}}}&\left(-m,a ,b  ;c , a+b-m -c  +2;1\right)
 = 
\frac{ \Gamma \left(b-m  -c +a +2\right) \Gamma \left(m -1-b +c \right) \Gamma \left(1-c +a \right) \Gamma \left(c \right)}{\Gamma \left(2-m -c +a \right) \Gamma \left(b -c +a +2\right) \Gamma \left(m +c \right) \Gamma \left(c -b \right)}\\
&\times\left(-c^{2}+\left(-m +b +a +2\right) c -\left(1-m +b \right) a +\left(m -1\right) \left(b +1\right)\right)\,,
\label{XidbTwo}
\end{align}
in agreement with \cite[Eq.5]{MathWSaal}. Specific examples of \eqref{Xidb} with $n=2,3$ and $4$ have been presented by Qureshi and Malik \cite{QuMal}, who acknowledge their use of ``complicated and lengthy algebraic calculations" of a ``cumbersome" nature.

The logic that has led us this far now suggests that we examine the Thomae progeny of each side of \eqref{HidX} in a search for other useful identities. Of the possibilities considered, the third Thomae transform \eqref{Thom3} produces
\begin{align}\nonumber
{}_{3}^{}{\moversetsp{}{\mundersetsp{}{F_{2}^{}}}}&\left(a +m ,b ,c ;c +n ,a +1;1\right)
 = 
\frac{ \Gamma \left(1-b +n -m \right) \Gamma \left(c +n \right)}{\Gamma \left(n \right) \Gamma \left(1+n-m -b +c  \right)} \\
&\times {}_{3}^{}{\moversetsp{}{\mundersetsp{}{F_{2}^{}}}}\left(1-m ,c ,1+a-b  ;1+c +n-m -b  ,a +1;1\right)\,
\label{Cx1b3}
\end{align}
an identity equivalent to that derived by Karp and Kalmykov \eqref{KarpKalmy}.
Setting $a=0$, then $c:=a$ in \eqref{Cx1b3} produces a special case of the (constrained) challenge \eqref{HAll2Ha}, and a limiting case of \eqref{KarpKalmy}, that is
\begin{align} \nonumber
{}_{3}^{}{\moversetsp{}{\mundersetsp{}{F_{2}^{}}}}&\left(a,b ,m  ;1,a +n;1\right)
 =
\frac{ \Gamma \left(1+n-b  -m \right) \Gamma \left(a +n \right)}{\Gamma \left(n \right) \Gamma \left(1+a +n-m -b  \right)}
{}_{3}^{}{\moversetsp{}{\mundersetsp{}{F_{2}^{}}}}\left(1-m ,a ,1-b ;1,1+a +n-m -b;1\right)\\
&=
\frac{\Gamma \left(m \right)  \Gamma \left(1+n-b -m  \right) \Gamma \left(a +n \right)}{\Gamma \left(a \right) \Gamma \left(1-b \right) \Gamma \left(n \right)}\moverset{m -1}{\munderset{k =0}{\sum}}\frac{\left(-1\right)^{k}\,\Gamma \left(a +k \right) \Gamma \left(1-b +k \right)}{\Gamma \left(m -k \right) \Gamma \left(1+k \right)^{2}\,\Gamma \left(1-m -b +a +n +k \right)}\,,
\label{Excb4}
\end{align}
the right-hand side equalities being terminating $n-$balanced sums. The relationship between the two identities \eqref{Excb4} and \eqref{HAll2Ha} with $e=1$, is pursued in the following  subsection.

\subsubsection{To evaluate an Euler-type sum} \label{sec:EulerSum}

Although it is a somewhat arduous calculation, it can be shown that evaluation of the limit $e\rightarrow 1$ in \eqref{HAll2Ha} indeed reduces that identity to \eqref{Excb4}. With the application of this limit, the last term (includes 4F3(1)) in \eqref{HAll2Ha} vanishes due to the presence of the term $1/\Gamma(e-m)$ and the other terms require a limiting process. When this is performed and the right-hand sides of \eqref{Excb4} and the so-limited \eqref{HAll2Ha} are equated, we eventually identify the following Euler-type sum
\begin{align} \nonumber
\moverset{m}{\munderset{k =0}{\sum}}&\frac{\left[\psi \left(1-k +m \right)-\psi \left(1+k -b -a \right)\right] \Gamma \left(1+k -b -a \right)}{\Gamma \left(1-m -b +k \right) \Gamma \left(1-a -m +k \right) \Gamma \left(1-k +m \right)^{2}\,\Gamma \left(k +1\right)}\\ \nonumber
 &= 
\frac{ \Gamma \left(1-a +m \right) \Gamma \left(m +1-b \right) \Gamma \left(1-b -a \right)}{\Gamma \left(1-b \right)^{2}\,\Gamma \left(1-a \right)^{2}\,\Gamma \left(1+m \right)^{2}}\\
&\times\left[\psi \left(a \right)+\psi \left(b \right)+\psi \left(1+m \right)-\psi \left(1-b -a \right)-\psi \left(a -m \right)-\psi \left(b -m \right)\right]\,,
\label{Fxcb9}
\end{align}
after a laborious journey through the jungle of computer simplification. See Appendix \ref{sec:Fxcb9Proof} for a derivation and a related sum; for an analysis of Euler-type sums convertible into pFq(1) see \cite{zbMATH08059632}. 
\subsection{Limiting cases of \eqref{Hnew2}} \label{sec:LimCase}

As suggested previously (Section \ref{sec:Hnew2}), it is interesting to evaluate \eqref{Hnew2} in the limiting case $c=b+m$ corresponding to the special case $c=n$ of \eqref{Hnew3}. From the form of the right-hand side, it is clear that straightforward substitution is invalid, so a limit must be evaluated. When this is done, one arrives at a complicated identity for ${}_{3}^{}{\moversetsp{}{\mundersetsp{}{F_{2}^{}}}}\left(a ,b ,n ;a +n ,b +m ;1\right)$ -- for the full expression see \eqref{Hn4b}. Of further interest is the special case $a=j$ along with $m\geq n$. This again requires a further limiting evaluation, eventually yielding an identity of the form
\begin{align} \nonumber
{}_{3}^{}{\moversetsp{}{\mundersetsp{}{F_{2}^{}}}}&\left(b ,j ,n ;j +n ,b +m ;1\right)
 = \frac{\pi\left(-1\right)^{j}\Gamma \left(b +m \right) \Gamma \left(j +n \right) }{\sin \left(\pi \,b \right)\Gamma \left(n -b +j \right)}\\
&\times\left(H \left(m , n , j , b\right)
+\left(-1\right)^{n} \Gamma(1-b) \left[T_{1}(m,n,j,b)+T_{2}(m,n,j,b)\right]\right)
\label{H4D}
\end{align}
where the definition of $H(m,n,j,b)$ is given in \eqref{H4i} and
\begin{equation}
T_{1}(m,n,j,b)\equiv\,\frac{-1}{\Gamma \left(b +m -n \right) \Gamma \left(n \right) \Gamma \left(j \right)}\moverset{n -1}{\munderset{k =0}{\sum}}\frac{\left(-1\right)^{k}\,\Gamma \left(n +m -1-k \right) \Gamma \left(b-1 +m -k \right)}{\Gamma \left(b -j +m -k \right) \Gamma \left(k +1\right) \Gamma \left(m -k \right) \Gamma \left(n -k \right)},
\label{T1mnjb}
\end{equation}
\begin{equation}
T_{2}(m,n,j,b)\equiv\frac{1}{\Gamma \left(m \right)}\moverset{j -1}{\munderset{k =0}{\sum}}\frac{\Gamma \left(m +k \right)}{\Gamma \left(j -k \right) \Gamma \left(k +1\right)^{2}\,\Gamma \left(b +m +1-j -n +k \right) \Gamma \left(n -k \right)}\,.
\end{equation}
Since the form of \eqref{H4D} naturally suggests that we consider the case $b=p$ where $p$ is a positive integer, we see that the term that multiplies $\Gamma(1-b)$ must vanish when $b=p$, because the left-hand side of \eqref{H4D} is convergent with a parametric excess equal to $m>0$ and there is no possibility that the term $H(m,n,j,b)$ can contribute to order $(b-p)^{-1}$ (see \eqref{H4i}) to cancel that divergence.
\subsubsection{Coefficient of $(b-p)^{-1}$ term}
For $m\geq n$, this requirement yields the identity
\begin{align} \nonumber
\moverset{J}{\munderset{k =0}{\sum}}&\frac{\Gamma \left(m +k \right)}{\Gamma \left(j -k \right) \Gamma \left(k +1\right)^{2}\,\Gamma \left(p +k \right) \Gamma \left(n -k \right)}
 =\frac{\Gamma \left(m \right) }{\Gamma \left(p -1+j \right) \Gamma \left(n \right) \Gamma \left(j \right)}\\ 
&\times\moverset{P}{\munderset{k =0}{\sum}}\frac{\left(-1\right)^{k}\,\Gamma \left(n +m -1-k \right) \Gamma \left(-2+p +j +n -k \right)}{\Gamma \left(p -1+n -k \right) \Gamma \left(k +1\right) \Gamma \left(m -k \right) \Gamma \left(n -k \right)}
\label{CxK}
\end{align}
where $J=\min \left(j -1, n -1\right)$ and $P=\min \left(m , n , p -1+n \right)-1$. This is equivalent to a hypergeometric transformation between terminating sums:
\begin{align} \nonumber
{}_{3}^{}{\moversetsp{}{\mundersetsp{}{F_{2}^{}}}}&\left(-n ,-j,m  ;1,p ;1\right)
 =
\frac{\Gamma \left(p \right) \Gamma \left(n +m \right) \Gamma \left(p +j +n \right) }{\Gamma \left(1+n \right) \Gamma \left(p +n \right) \Gamma \left(m \right) \Gamma \left(p +j \right)}\\
 &\times{}_{3}^{}{\moversetsp{}{\mundersetsp{}{F_{2}^{}}}}\left(1-m ,-n ,1-(p +n) ;1-(p+j +n) ,1-(m +n) ;1\right)
\label{CxA}
\end{align}
where the variables have been redefined such that $p:=p-m-1+j+n$, followed by $n:=n+1,~j:=j+1$.
\subsubsection{Limiting and special cases of \eqref{H4D}}
One could now proceed to evaluate \eqref{H4D} for limiting, but still general, values of $b=p$, leading to a morass of specific cases that depend on numerous permutations of inequalities among $m,n,j,p$, many of which are already listed by Prudnikov et. al. \cite[Section 7.7.4]{prudnikov}, and consist of variations of the challenge problem. However, by relaxing the generality of the possibilities, simpler evaluations are available that still generalize those found elsewhere. For example, for the case $j=1$ and non-integer values of $b$, \eqref{H4D} reduces to
\begin{align} \nonumber
{}_{3}^{}{\moversetsp{}{\mundersetsp{}{F_{2}^{}}}}&\left(1,b ,n ;b +m ,1+n ;1\right)
 = 
\frac{\,n\, \Gamma \left(b +m \right) \Gamma(1-b) }{\Gamma \left(n -b +1\right) \Gamma \left(b +m -n \right)}\\ \nonumber
 &
\times \left(\moverset{n -1}{\munderset{k =0}{\sum}}\frac{\left(-1\right)^{n +k}\,\Gamma \left(n +m -1-k \right) \left(\psi \left(b +m -k-1 \right)-\psi \left(n -k \right)\right)}{\Gamma \left(n -k \right) \Gamma \left(k +1\right) \Gamma \left(m -k \right)}\right.
\\&\left.\hspace{30pt}+\moverset{m-n -1}{\munderset{k =0}{\sum}}\frac{\Gamma \left(k +1\right) \Gamma \left(m -1-k \right)}{\Gamma \left(n +k +1\right) \Gamma \left(m -n -k \right)}\right)\,.
\label{Cj1a}
\end{align}
and the similar case $j=2$ and $m\geq n$, becomes
%
%
\begin{align} \nonumber
{}_{3}^{}{\moversetsp{}{\mundersetsp{}{F_{2}^{}}}}&\left(2,b ,n ;b +m ,2+n ;1\right)
 =\, -\,
\frac{ \pi\,\Gamma \left(b +m \right) \Gamma \left(2+n \right)}{\Gamma \left(n -b +2\right) \Gamma \left(b +m -n \right) \sin \left(\pi \,b \right) \Gamma \left(n \right) \Gamma \left(b \right)}\\ \nonumber
&\times\left(\  \moverset{m -1}{\munderset{k =n}{\sum}}\frac{\Gamma \left(-n +1+k \right) \left(b +m -k-2 \right) \Gamma \left(n +m -1-k \right)}{\Gamma \left(k +1\right) \Gamma \left(m -k \right)} +\left(-1\right)^{n} \left(2-b -m +n \right)\right. \\
&\left.+\left(-1\right)^{n} \moverset{n -1}{\munderset{k =0}{\sum}}\frac{\left(-1\right)^{k} \left(b +m -k-2 \right) \Gamma \left(n +m -1-k \right) \left(\psi \left(b +m -k-2 \right)-\psi \left(n -k \right)\right)}{\Gamma \left(n -k \right) \Gamma \left(k +1\right) \Gamma \left(m -k \right)} \right)\,.
\label{Cj2a}
\end{align}
Carried further, by evaluating \eqref{Cj1a} in the limiting case $b\rightarrow 1$, for $m\geq n$, we obtain
\begin{align} \nonumber
{}_{3}^{}{\moversetsp{}{\mundersetsp{}{F_{2}^{}}}}&\left(1,1,n ;1+m ,1+n ;1\right)\\ \nonumber
 &= 
n\left(\frac{ \Gamma \left(1+m \right) \left(-1\right)^{n}}{\Gamma \left(1+m -n \right) \Gamma \left(n \right)}\left(-\moverset{n -1}{\munderset{k =0}{\sum}}\frac{\left(-1\right)^{k}\,\psi \left(1, m -k \right) \Gamma \left(n +m -1-k \right)}{\Gamma \left(m -k \right) \Gamma \left(k +1\right) \Gamma \left(n -k \right)}\right.\right.\\&\nonumber \left.\left.
 \hspace{30pt}+\gamma^{2}+\left(\psi \left(m \right)+\psi \left(n \right)\right) \gamma +\psi \left(1+m -n \right) \left(\psi \left(m \right)+\psi \left(n \right)-\psi \left(1+m -n \right)\right)\right)
\right.\\  &\left. 
-m\left(\psi \left(1+m -n \right)+\gamma \right) \moverset{n-2}{\munderset{k =0}{\sum}}\frac{\Gamma \left(n +m -1-k \right) \left(-1\right)^{k}}{\Gamma \left(k +1\right) \Gamma \left(m -k \right) \Gamma \left(n -k \right) \left(n -1-k \right)}  \right)\,.
\label{Hbj11b}
\end{align}
\begin{rem} {\bf \phantom{x} } 
\begin{itemize}
\item{}
Under the (left-to-right ordered) mapping $m:=j,~a:=n,~e:=b+m$, \eqref{H4i} (alternatively \eqref{H4D}) is a limiting case of \eqref{Hnew};
\item{}
\eqref{Cj1a} and \eqref{Hbj11b} are limiting cases of entry 25 of \cite{Milgram447}, and \eqref{Cj2a} generalizes that same entry (see \cite[Theorem 1]{MilgRaoComment});
\item{}
For the case $n>m$ in \eqref{Hbj11b}, simply interchange $m$ and $n$ on the right-hand side;
\item{}
In the limiting case $b \rightarrow 1$, the sum of the terms of \eqref{Cj1a} that are enclosed in parentheses must vanish because of the presence of the factor $\Gamma(1-b)$ . Written in hypergeometric form with $m>n$, this identifies 
\begin{align} \nonumber
{}_{3}^{}{\moversetsp{}{\mundersetsp{}{F_{2}^{}}}}&\left(1,1,1-m +n ;1+n ,2-m ;1\right)
 = 
\frac{ \left(-1\right)^{n}\Gamma \left(1+n \right) \Gamma \left(m -n \right)}{\Gamma \left(m -1\right) }\\
&\times\moverset{n -1}{\munderset{k =0}{\sum}}\frac{\left(-1\right)^{k} \left(\psi \left(n -k \right)-\psi \left(m -k \right)\right) \Gamma \left(n +m -1-k \right)}{\Gamma \left(k +1\right) \Gamma \left(m -k \right) \Gamma \left(n -k \right)},
\label{Hbj11Id}
\end{align}
an identity that can be independently verified otherwise.
\item{}
Experimentally, it was observed that \eqref{CxA} is numerically valid for the more general case $p\Rightarrow b $, which happens to coincide with a Thomae progeny of database entry 251. Comparing the right-hand sides of both yields a large variety of transformations of terminating series. For example
\begin{align} \nonumber
{}_{3}^{}{\moversetsp{}{\mundersetsp{}{F_{2}^{}}}}&\left(-j ,-m ,a ;-(j +m) ,n +1+a ;1\right)
 \\& \nonumber
=\frac{\Gamma \left(m +n +1\right) \Gamma \left(n +a +1+m +j \right)\Gamma \left(n +1+a \right) \Gamma \left(1+j \right)}{\Gamma \left(1+n \right) \Gamma \left(a +1+m +n \right) \Gamma \left(1+m +j \right) \Gamma \left(n +1+a +j \right)}\\&
\times{}_{3}^{}{\moversetsp{}{\mundersetsp{}{F_{2}^{}}}}\left(-m ,-n ,-a -(m +n) ;-a -(m +j+n) ,-(m+n) ;1\right), 
\label{CxAb}
\end{align}
a special case of the challenge problem, both sides of which are not Thomae related.
\end{itemize}
\end{rem}

\subsection{Special cases of \eqref{HAll2Ha}} \label{sec:Specials1}

\subsubsection{To determine a special case $_{4}F_{3}(1)$} \label{sec:MEq1}

After an ordered redefinition of variables ($m:=n$ then $a:=m$) the left-hand side of \eqref{HAll2Ha} becomes $_{3}{\moversetsp{}{\mundersetsp{}{F_{2}}}}\left(m ,n, b;e ,m +n ;1\right)$, a known sum (see \eqref{Hnew}), in which case we are able to solve for the 4F3(1) appearing in the corresponding right-hand side by evaluating the limit $a\rightarrow n$ in that identity. After lengthy computation, for $m\geq n$, we eventually arrive at a summable form of the left-hand side of \eqref{HAll2Ha}, that is
%
\begin{align} \nonumber
&{}_{3}^{}{\moversetsp{}{\mundersetsp{}{F_{2}^{}}}}\left(m ,n,b ;e ,m +n ;1\right)
 = \\ \nonumber
&\frac{\Gamma \left(e \right)^{2}\,\Gamma \left(m +n \right) \left(\pi \,\cot \left(\pi\left(b-e \right)  \right)-\psi \left(e -m \right)\right) \left(-1\right)^{m}  \Gamma \left(1-e \right) \Gamma(1-b)}{ \Gamma \left(n+m-b \right) \Gamma \left(e -n \right) \Gamma \left(n \right) \Gamma \left(m \right)^{2}\,\Gamma \left(1-n -m +b \right)  \Gamma \left(1+n -e \right)}\\ \nonumber
&\times\moverset{n -1}{\munderset{k =0}{\sum}}\frac{\left(-1\right)^{k}\,\Gamma \left(1-n -m +b +k \right) \Gamma \left(n +m -1-k \right) \Gamma \left(m \right)}{\Gamma \left(m -k \right) \Gamma \left(1+k \right) \Gamma \left(n -k \right) \Gamma \left(1-m +e -n +k \right)}\\ \nonumber
&+\frac{\left(-1\right)^{m +n}\,\pi \,\Gamma \left(m +n \right) \Gamma \left(e \right) }{\Gamma \left( n +m -b\right) \Gamma \left(e -n \right)} \\ \nonumber
&\times
\left(\frac{\csc \left(\pi \,b \right) }{\Gamma \left(1-e +b \right) \Gamma \left(n \right)^{2}}
\moverset{m -1}{\munderset{k =0}{\sum}}\frac{\left(\psi \left(1-e +b +k \right)-\psi \left(1+k \right)\right) \left(-1\right)^{k}\,\Gamma \left(1-e +b +k \right) \Gamma \left(n +k \right)}{\Gamma \left(b +1+k -m \right) \Gamma \left(m -k \right) \Gamma \left(1+k \right)^{2}} \right. \\ 
& \left.  -\,\frac{\csc \left(\pi \,e \right) }{\Gamma \left(m \right) \Gamma \left(n \right) \Gamma \left(e -m \right)}\moverset{n -1}{\munderset{k =0}{\sum}}\frac{\Gamma \left(-b +1+k \right) \left(-1\right)^{k}\,\psi \left(m -k \right) \Gamma \left(n +m -1-k \right)}{\Gamma \left(m -k \right) \Gamma \left(1+k \right) \Gamma \left(n -k\right) \Gamma \left(2-e +k \right)}\right)\,.
\label{Eq479d}
\end{align}
\begin{rem}

\begin{itemize}
\item{For the case $n> m$, simply interchange $m$ and $n$ on the right-hand side of \eqref{Eq479d};} 

\item{} As happens often (e.g. \eqref{Hid}), while evaluating a limit $a\rightarrow n$ in \eqref{Hnew}, a term of order $(a-n)^{-1}$ arose, thereby indicating that the coefficient of that term must identically vanish since the sum is convergent under that condition. This yields the transformation identity where $e-n>0$
\begin{align} \nonumber
{}_{3}^{}{\moversetsp{}{\mundersetsp{}{F_{2}^{}}}}&\left(1-m ,n ,1+b-e ;1,1-m +b ;1\right)
 \\ 
&=\frac{\Gamma \left(m +n -1\right) \Gamma \left(e -1\right)\Gamma \left(1-m +b \right)}{\Gamma \left(m \right) \Gamma \left(n \right) \Gamma \left(e -m \right) \Gamma \left(b \right)} {}_{3}^{}{\moversetsp{}{\mundersetsp{}{F_{2}^{}}}}\left(1-n,1-m ,1-b  ;2-e ,2-m -n ;1\right) 
\label{HId479}
\end{align}
providing considerable simplification in the derivation of \eqref{Eq479d}.
\end{itemize}
\end{rem}
By substituting \eqref{Eq479d} into \eqref{HAll2Ha} along with $a=n$ it is possible to solve for the $4F3(1)$ in that identity  and obtain a complicated finite series representation for the hypergeometric function ${}_{4}^{}{\moversetsp{}{\mundersetsp{}{F_{3}^{}}}}\left(1,1,b ,e +m ;b +m ,n +1,e +1;1\right)$ after redefining $b:=b+m-1,e:=e+m$ to render the result into canonical form. For general values of $n$ the identity is unwieldy; for the case $n=1,m\geq 1$ we obtain
\begin{align} \nonumber
{}_{4}^{}&{\moversetsp{}{\mundersetsp{}{F_{3}^{}}}}\left(1,1,b ,e +m ;2,b +m ,e +1;1\right)
 = \frac{e \left(b +m -1\right) \left(\psi \left(e \right)-\psi \left(e +m -1\right)+\psi \left(b +m -1\right)\right)}{\left(b-1 \right) \left(e +m -1\right)}\\
&-\frac{\Gamma \left(e +1\right) \Gamma \left(b +m \right) \Gamma \left(m \right) }{\left(b-1 \right) \Gamma \left(e +m \right) \Gamma \left(b-e  \right)}\moverset{m -1}{\munderset{k =0}{\sum}}\frac{\left(-1\right)^{k}\,\Gamma \left( b +k-e \right) \psi \left(1+k \right)}{\Gamma \left(b +k \right) \Gamma \left(1+k \right) \Gamma \left(m -k \right)} \,.
\label{Q3c}
\end{align}
\begin{rem}
\begin{itemize}
\item{}
The identity \eqref{Q3c} has been obtained by first differentiating the elementary identity
\begin{equation}
\moverset{m-1}{\munderset{k =0}{\sum}}\frac{\left(-1\right)^{k}\,
\Gamma \left( b +k -e\right)}{\Gamma \left(b +k \right) \Gamma \left(1+k \right) \Gamma \left(m -k \right)}
 = 
\frac{\Gamma \left(b-e \right) \Gamma \left(e +m -1\right)}{\Gamma \left(m \right) \Gamma \left(e \right) \Gamma \left(b +m -1\right)}
\label{S2F1}
\end{equation}
with respect to the variable $e$, giving
\begin{align} \nonumber
\moverset{m-1}{\munderset{k =0}{\sum}}&\frac{\left(-1\right)^{k}\,\psi \left(b+k-e \right) \Gamma \left( b +k-e \right)}{\Gamma \left(b +k \right) \Gamma \left(1+k \right) \Gamma \left(m -k \right)}\\
& = 
\frac{\left(\psi \left(e \right)+\psi \left(b-e  \right)-\psi \left(e +m -1\right)\right) \Gamma \left(e +m -1\right) \Gamma \left(b-e \right)}{\Gamma \left(m \right) \Gamma \left(e \right) \Gamma \left(b +m -1\right)}
\label{SX}
\end{align}
and then with respect to the variable $b$, giving
\begin{equation}
\moverset{m -1}{\munderset{k =0}{\sum}}\frac{\left(-1\right)^{k}\,\Gamma \left(1 +k -e\right) \psi \left(1+k \right)}{\Gamma \left(1+k \right)^{2}\,\Gamma \left(m -k \right)}
 = 
\frac{\pi\,\left(\psi \left(e \right)-\psi \left(e +m -1\right)+\psi \left(m \right)\right) \Gamma \left(e +m -1\right) }{\sin \left(\pi \,e \right) \Gamma \left(e \right)^{2}\,\Gamma \left(m \right)^{2}}
\label{B3m1}
\end{equation}
in the limit $b\rightarrow 1$; 

\item{} If $e=n,n>0$ in \eqref{B3m1}, we find
\begin{equation}
\moverset{N}{\munderset{k =0}{\sum}}\frac{\psi \left(1+k \right)}{\Gamma \left(n -k \right)\,\Gamma \left(m -k \right) \Gamma \left(1+k \right)^{2}}
 = 
\frac{\left(\psi \left(n \right)-\psi \left(n +m -1\right)+\psi \left(m \right)\right) \Gamma \left(n +m -1\right)}{\Gamma \left(n \right)^{2}\,\Gamma \left(m \right)^{2}}\,,
\label{B3en}
\end{equation}
where $N\equiv \mathrm{min}(m-1,n-1)$. Setting $e=0$ in \eqref{B3m1} with $m:=m+1$ yields  the Euler-type sum
\begin{equation}
\moverset{m}{\munderset{k =0}{\sum}}\left(-1\right)^{k}\,{\binom{m}{k}}\,\psi \left(k +1\right)
 = -\frac{1}{m}\,.
\label{B3e0}
\end{equation}
\end{itemize}
\end{rem}

\subsubsection{Limiting Cases of \eqref{Q3c}}

If m=1, \eqref{Q3c} reduces to
\begin{equation}
{}_{3}^{}{\moversetsp{}{\mundersetsp{}{F_{2}^{}}}}\left(1,1,b ;2,b +1;1\right)
 = \frac{b \left(\gamma +\psi \left(b \right)\right)}{b-1}\,,
\label{Q3ca}
\end{equation}
a special case of \cite[Eq. 3.13(42)]{Luke}, which itself is a special case of the well-known \cite[Eq. 3.13(41)]{Luke} general form ($1+c-a-b>0$):
\begin{equation}
{}_{3}^{}{\moversetsp{}{\mundersetsp{}{F_{2}^{}}}}\left(a,b,1 ;c,2;1\right)=\frac{(c-1)}{(a-1)(b-1)}\left(\frac{\Gamma(c-1)\Gamma(c-a-b+1)}{\Gamma(c-a)\Gamma(c-b)}-1\right);
\label{LukeGen}
\end{equation}
See also \eqref{Lukep1111} below.

In \eqref{Q3c}, also consider the special case $b=1$; this involves a limit similar to previous calculations and we eventually arrive at
\begin{align} \nonumber
{}_{4}^{}{\moversetsp{}{\mundersetsp{}{F_{3}^{}}}}&\left(1,1,1,e +m ;2,e +1,1+m ;1\right)
 \\ \nonumber
&=\frac{m\,e\, \Gamma \left(m \right)^{2}\,\Gamma \left(e \right)^{2} \sin \left(\pi \,e \right)}{\pi \,\Gamma \left(e +m \right)}\moverset{m -1}{\munderset{k =0}{\sum}}\frac{\left(-1\right)^{k}\,\psi \left(1+k \right)^{2}\,\Gamma \left(1-e +k \right)}{\Gamma \left(1+k \right)^{2}\,\Gamma \left(m -k \right)} \\
&-\frac{m\,e \left(\left[\psi \left(e \right)-\psi \left(e +m -1\right)+\psi \left(m \right)\right]^{2}-\psi \left(1, e\right)-\psi \left(1, m\right)+\psi \left(1, e +m -1\right)\right)}{e +m -1},
\label{Bn1b1}
\end{align}
yielding simpler special cases such as
\begin{equation}
{}_{4}^{}{\moversetsp{}{\mundersetsp{}{F_{3}^{}}}}\left(1,1,1,e +2;2,3,e +1;1\right)
 = \frac{\pi^{2}\,e -6\,e +6}{3\,e +3}\,,
\label{Bn1b1m3}
\end{equation}
if $m=2$. We also note that \eqref{Bn1b1} both specializes and extends identities studied by Karp and Prilepkina \cite{karpandPril2022}.
\subsubsection{Limiting Cases of \eqref{HAll2Ha}} \label{sec:Lim3p15}
Alternatively to the above, if we first set $m=1$ in \eqref{HAll2Ha}, followed by the choice $n=1$ and then let $a=n$, after considerable calculation involving limiting cases we obtain
\begin{align} \nonumber
{}_{3}^{}&{\moversetsp{}{\mundersetsp{}{F_{2}^{}}}}\left(1,b ,n ;e ,n +1;1\right)
 = \frac{n\,\Gamma \left(b -n \right) \Gamma \left(e \right)}{\Gamma \left(e -n \right) \Gamma \left(b \right)} \left[\psi \left(e -n \right)-\psi \left(e -b \right)+\psi \left(b \right)\right] \\ \nonumber
&-\frac{ n\,\left(e -1\right)\Gamma \left(b -n +1\right) \Gamma \left(e \right)}{\Gamma \left(1+b \right) \Gamma \left(1+e -n \right)} {}_{4}^{}{\moversetsp{}{\mundersetsp{}{F_{3}^{}}}}\left(1,1,e ,1+b -n ;2,1+b ,1+e -n ;1\right)\\
&-\frac{ n \,  \Gamma \left(1+n -e \right) \Gamma \left(1-b \right)}{\Gamma \left(1-b +n \right) \Gamma \left(1-e \right)}\psi \left(n \right)\,.
\label{Qm1n1}
\end{align}
Since the left-hand side of \eqref{Qm1n1} can be identified with the help of \eqref{Hnew} and therefore transformed into the summed form
\begin{align} \nonumber
{}_{3}^{}&{\moversetsp{}{\mundersetsp{}{F_{2}^{}}}}\left(1,b ,n ;e ,n +1;1\right)
 = 
\frac{\left(-1\right)^{n}\,n\,\pi  \left(e -1\right)}{\sin \left(\pi \,e \right) \Gamma \left(e -n \right) \Gamma \left(1-b +n \right)} \moverset{n -1}{\munderset{k =1}{\sum}}\frac{\Gamma \left(1-b +k \right)}{k\,\Gamma \left(2-e +k \right)}\\
&-\frac{\left(-1\right)^{n}\,n\,\Gamma \left(1-b \right) \Gamma \left(e \right) \left(-\psi \left(e -1\right)+\psi \left(e -b \right)\right)}{\Gamma \left(e -n \right) \Gamma \left(1-b +n \right)}\,,
\label{Xt0}
\end{align}
it is possible then to solve \eqref{Qm1n1} to find
\begin{align} \nonumber
{}_{4}^{}{\moversetsp{}{\mundersetsp{}{F_{3}^{}}}}&\left(1,1,b ,e +n -1;2,e ,b +n ;1\right)
 = 
\frac{\left(1-e \right) \sin \left(\pi \,b \right) \Gamma \left(b +n \right) }{\sin \left(\pi \,e \right) \left(b-1 \right) \Gamma \left(e +n -1\right)}\moverset{n -1}{\munderset{k =1}{\sum}}\frac{\Gamma \left(2+k-b -n \right)}{k\,\Gamma \left(3-e -n +k \right)} \\
&+\frac{\left(1-e \right) \left(b +n -1\right) }{\left(b-1 \right) \left(e +n -2\right)}\left(\psi \left(n \right)+\psi \left(e +n -2\right)-\psi \left(b +n -1\right)-\psi \left(e -1\right)\right)\,,
\label{Q2}
\end{align}
yielding some interesting special cases, particularly relative to general ones studied by Karp and Prilepkina \cite{karpandPril2022}. That is, if $n=2$, then
\begin{align} \nonumber
{}_{4}^{}{\moversetsp{}{\mundersetsp{}{F_{3}^{}}}}&\left(1,1,b ,e +1;2,e ,b +2;1\right)
  \\&
=\frac{\left(1+b \right) \left(e -1\right) \gamma}{e \left(b-1\right)}+\frac{\left(1+b \right) \left(e -1\right) \psi \left(b \right)}{e \left(b-1 \right)}+\frac{\left(1+b \right) \left(b^{2}-b\,e +e -1\right)}{b\,e \left(b-1 \right)}
\label{Q2_2}
\end{align}
and if $n=3$, then
\begin{align} \nonumber
{}_{4}^{}{\moversetsp{}{\mundersetsp{}{F_{3}^{}}}}&\left(1,1,b ,e +2;2,e ,b +3;1\right)\\
 &= 
\frac{\left(e -1\right) \left(b +2\right) \left(\gamma +\psi \left(b \right)\right)}{\left(b-1 \right) \left(e +1\right)}+\frac{\left(b^{2}+3\,b\,e +2\,e \right) \left(1-e +b \right) \left(b +2\right)}{2 \left(e +1\right) e\,b \left(1+b \right)}\,.
\label{Q2-3}
\end{align}
Further, if $e=m$ and $m,n\geq 2$, then
\begin{align} \nonumber
{}_{4}^{}{\moversetsp{}{\mundersetsp{}{F_{3}^{}}}}&\left(1,1,b ,m +n -1;2,m ,b +n ;1\right)\\ \nonumber
&
=\frac{\left(-1\right)^{n}\,\left(1-m\right) \sin \left(\pi \,b \right)  \Gamma \left(b +n \right)}{\pi  \left(b-1 \right) \Gamma \left(m +n -1\right)}\moverset{n -1}{\munderset{k =1}{\sum}}\frac{\left(-1\right)^{k}\,\Gamma \left(2+k-b -n  \right) \Gamma \left(m +n -k -2\right)}{k}\\
&+\frac{\left(1-m\right)\left(b +n -1\right) \left[\psi \left(n \right)+\psi \left(n -2+m \right)-\psi \left(m -1\right)-\psi \left(b +n -1\right)\right] }{\left(b-1 \right) \left(n +m-2 \right)}
\label{Qne1}
\end{align}
and, if $m=1,n\geq2$ then the 4F3(1) reduces to a 3F2(1) yielding
\begin{equation}
{}_{3}^{}{\moversetsp{}{\mundersetsp{}{F_{2}^{}}}}\left(1,b ,n ;2,b +n ;1\right)
 = 
\frac{b +n -1}{\left(1-b \right) \left(n -1\right)}+\frac{\Gamma \left(b +n \right)}{\left(b-1 \right) \left(n -1\right) \Gamma \left(n \right) \Gamma \left(b \right)}\,,
\label{Q2-M1}
\end{equation}
coinciding with special cases of a large number of database entries, including the known (see \eqref{Q3ca}) special case $b\rightarrow 1$
\begin{equation}
{}_{3}^{}{\moversetsp{}{\mundersetsp{}{F_{2}^{}}}}\left(1,1,n ;2,1+n ;1\right)
 = 
\moverset{\infty}{\munderset{k =0}{\sum}}\frac{n}{\left(1+k \right) \left(n +k \right)}
 = \frac{n \left(\gamma +\psi \left(n \right)\right)}{n -1}\,.
\label{Q2b1}
\end{equation}

\subsection{Based on ${}_{4}^{}{\moversetsp{}{\mundersetsp{}{F_{3}^{}}}}\left(1,1,e ,d ;2,2,2;1\right)$} 

We begin with the known identity \cite[Eq. (7.1)]{MilgVar}
\begin{align} \nonumber
{}_{4}^{}&{\moversetsp{}{\mundersetsp{}{F_{3}^{}}}}\left(1,1,e ,d ;2,2,2;1\right)
 =\,-\frac{2\,\gamma +\psi \left(e -1\right)+\psi \left(2-d \right)}{\left(e -1\right) \left(d -1\right)}\\  
&-\frac{\sin \left(\pi \,e \right) \Gamma \left(e \right) \Gamma \left(1-d \right) {}}{\Gamma \left(1+e -d \right) \left(e -1\right)^{3}\,\pi}{_{3}^{}\moversetsp{}{\mundersetsp{}{F_{2}^{}}}}\left(e -1,e -1,e -1;e ,1+e -d ;1\right)\,. 
\label{Eq7p1}
\end{align}
Set $e=a+1$ and $d=3-a-b$ to yield
\begin{align} \nonumber
{}_{4}^{}{\moversetsp{}{\mundersetsp{}{F_{3}^{}}}}&\left(1,1,a+1 ,3-a -n ;2,2,2;1\right)
 = -\frac{2\,\gamma +\psi \left(a \right)+\psi \left(a -1+n \right)}{a \left(2-a -n \right)}\\ 
&+\frac{\sin \left(\pi  a \right) \Gamma \left(1+a \right) \Gamma \left(a +n-2 \right)}{\pi\,\Gamma \left(2\,a +n-1 \right) a^{3}} {}_{3}^{}{\moversetsp{}{\mundersetsp{}{F_{2}^{}}}}\left(a ,a ,a ;1+a ,2\,a-1 +n ;1\right)
\label{eq7p1a}
\end{align}
and note that the 3F2(1) on the right-hand side of \eqref{eq7p1a} is contiguous to Watsons's theorem (see Appendix \ref{sec:Contig}), i.e.
\begin{align} \nonumber
{}_{3}^{}{\moversetsp{}{\mundersetsp{}{F_{2}^{}}}}&\left(a ,a ,a ;1+a ,2\,a -1+n ;1\right)=W_{1,n}(a,a,a) 
 \\&= 
\,\frac{a\,\Gamma \left(2\,a +n -1\right) \Gamma \left(n \right) }{2\,\Gamma \left(a \right)^{2}}\moverset{n -1}{\munderset{k=0}{\sum}}\frac{\left(-1\right)^{k} \left(2\,a -1+2\,k \right) \left(\psi \left(1, \frac{a+k}{2}\right)-\psi \left(1, \frac{1+a+k}{2}\right)\right) \Gamma \left(2\,a -1+k \right)}{\Gamma \left(n -k \right) \Gamma \left(2\,a +n +k -1\right) \Gamma \left(1+k \right)}\,.
\label{E3}
\end{align}
Substituting the above into \eqref{HAll2Ha} eventually reveals a new database identity \eqref{data494}, a very special case of the Miller/Paris challenge. 
\begin{rem}
Since the left-hand side of \eqref{Eq7p1} is symmetric under the interchange $e\leftrightarrow d$ and the right-hand side isn't, it is interesting to try $d=a+1$ and $e=3-a-b$ in \eqref{Eq7p1} and repeat the above calculation, eventually arriving at the same result, although this requires a considerable amount of computational fortitude.
\end{rem}
Setting $a=1$ in \eqref{E3}, produces the identity
\begin{equation}
{}_{3}^{}{\moversetsp{}{\mundersetsp{}{F_{2}^{}}}}\left(1,1,1;2,1+n ;1\right)
 = 
n \Gamma \left(n \right)^{2}\moverset{n -1}{\munderset{k =0}{\sum}}\frac{\left(-1\right)^{k+1} \left(\frac{1}{2}+k \right) \left(\psi \left(1, 1+\frac{k}{2}\right)-\psi \left(1, \frac{1}{2}+\frac{k}{2}\right)\right)}{\Gamma \left(n -k \right) \Gamma \left(1+n +k \right)},
\label{Z3a}
\end{equation} 
which can be compared with the known identity \cite[Eq. 3.13(43)]{Luke} 
\begin{equation}
{}_{3}^{}{\moversetsp{}{\mundersetsp{}{F_{2}^{}}}}\left(1,1,1;2,c ;1\right)=(c-1)\psi(1,c-1)\hspace{30pt} c>1\,,
\label{Lukep1111}
\end{equation}
to yield the (new?) identity
\begin{equation}
\moverset{n -1}{\munderset{j =0}{\sum}}\frac{\left(-1\right)^{j} \left(j +\frac{1}{2}\right) \left(\psi \left(1, 1+\frac{j}{2}\right)-\psi \left(1, \frac{1}{2}+\frac{j}{2}\right)\right)}{\Gamma \left(n -j \right) \Gamma \left(1+n +j \right)}
 = \,-\frac{\psi \left(1, n\right)}{\Gamma \left(n \right)^{2}}\,.
\label{L43}
\end{equation}
For comparison, setting $n=1$ in \eqref{data494} yields the known identity
\begin{equation}
{}_{3}^{}{\moversetsp{}{\mundersetsp{}{F_{2}^{}}}}\left(a ,1,1;a +1,a +1;1\right)
 =\, 
-\frac{a^{2} \left(\psi \left(1, \frac{1}{2}+\frac{a}{2}\right)-\psi \left(1, \frac{a}{2}\right)\right)}{2}
\label{EnEq1}
\end{equation}
that is one positive unit (in both directions) contiguous to both Whipple and Dixon's theorems.

\section{Other Exotics} \label{sec:Applics}
In the following sections, we review recent additions to the literature concerning 3F2(1). Among other things, this includes a study of so-called exotic sums of the 3F2(1) genre. The purpose is to show that many such identities can be simply derived and generalized by the use of the methods demonstrated in the previous sections.

\subsection{Chen and Chu \cite{ChuChen1} } \label{sec:ChenChu}
In a series of papers motivated by the work of Asakura, Otsubo and Terasoma \cite{AsakuraEtAl} who considered a subset of the Miller/Paris challenge (viz. $_{3}^{}{\moversetsp{}{\mundersetsp{}{F_{2}^{}}}}\! \left(a,b,x;a+b,x+1;1\right)$), Chen and Chu studied the case $a=y,~b=1-y$ and noted that this led to a large number of closed forms for $x=1/2$, with $y$ taking on particular fractional values. These they labelled ``exotic" in concert with earlier work by Chu \cite{Chu700} where $\sim 700$ individual instances were listed, some of which overlap with identities found by  Krupnikov and K{\"o}lbig  \cite[Table 2]{KrupKol}. Here we shall demonstrate that the Chen and Chu results can be reproduced and generalized by the simple expedient of identifying them as special cases contiguous to the Watson identity \eqref{Watson} -- which happens to coincide with the starting point of Asakura et.al. 
\subsubsection{Chen and Chu (basic)}

In the first paper under consideration, Chu and Chen \cite{ChuChen1} note that the derivation of Asakura et. al. requires a ``lengthy proof of five pages". They then provide a simplified proof involving integral representations of ${}_{3}^{}{\moversetsp{}{\mundersetsp{}{F_{2}^{}}}}\left(y ,1-y,\frac{1}{2} ;1,\frac{3}{2};1\right)$, occupying about two pages, that reproduces the four exemplar identities of Asakura et. al., who consider $y=\{1/2,1/3,1/4,1/6\}$ to obtain:
\begin{align}
&{}_{3}^{}{\moversetsp{}{\mundersetsp{}{F_{2}^{}}}}\left(\frac{1}{2},\frac{1}{2},\frac{1}{2};1,\frac{3}{2};1\right)
 = \frac{4\,\mathit{G}}{\pi};
\label{AOT1}\\
&{}_{3}^{}{\moversetsp{}{\mundersetsp{}{F_{2}^{}}}}\left(\frac{1}{3},\frac{2}{3},\frac{1}{2};1,\frac{3}{2};1\right)
 = \frac{3\,\sqrt{3}\,\ln \left(2\right)}{\pi};
\label{AOT2}\\
&{}_{3}^{}{\moversetsp{}{\mundersetsp{}{F_{2}^{}}}}\left(\frac{1}{4},\frac{3}{4},\frac{1}{2};1,\frac{3}{2};1\right)
 = \frac{4\,\ln \left(1+\sqrt{2}\right)}{\pi};
\label{AOT3}\\
&{}_{3}^{}{\moversetsp{}{\mundersetsp{}{F_{2}^{}}}}\left(\frac{1}{6},\frac{5}{6},\frac{1}{2};1,\frac{3}{2};1\right)
 = \frac{3\,\sqrt{3}\,\ln \left(2+\sqrt{3}\right)}{2\,\pi}.
\label{AOT4}
\end{align}
Here, we observe that the fundamental element under consideration can be identified as a special case contiguous to Watson's identity \eqref{Watson}, evaluated \cite{ChuW} years earlier by Chu himself in  the form of a terminating sum:
\begin{equation}
{}_{3}^{}{\moversetsp{}{\mundersetsp{}{F_{2}^{}}}}\! \left(y ,1-y ,1/2 ;1  ,\frac{3}{2};1\right)=\munderset{c \rightarrow \frac{1}{2}}{\mathrm{lim}} W_{1 ,0}\! \left(y , 1-y , c\right). 
\label{AOT}
\end{equation}
A complexity now arises because the straightforward substitution $c=1/2$ in \eqref{AOT} generates an exceptional case during its evaluation, so a simple limit must be employed. When this is done, we obtain, in only a few lines, the generalization 
\begin{equation}
{}_{3}^{}{\moversetsp{}{\mundersetsp{}{F_{2}^{}}}}\left(\frac{1}{2},y ,1-y ;1,\frac{3}{2};1\right)
 = 
\frac{1}{2\,y -1}-\left(\psi \left(\frac{1}{2}+\frac{y}{2}\right)-\psi \left(\frac{y}{2}\right)\right)\frac{\sin \left(\pi\,y \right)}{\pi  \left(2\,y -1\right)}\,. 
\label{AOTP}
\end{equation}
To reproduce \eqref{AOT1}, we evaluate the limit $y\rightarrow 1/2$ in \eqref{AOTP}, producing
\begin{equation}
{}_{3}^{}{\moversetsp{}{\mundersetsp{}{F_{2}^{}}}}\left(\frac{1}{2},\frac{1}{2},\frac{1}{2};1,\frac{3}{2};1\right)
 = 
\frac{2\,G}{\pi}+\frac{\pi^{2}-\psi \left(1, \frac{3}{4}\right)}{4\,\pi}\,,
\label{AOT1a}
\end{equation}
reducing to \eqref{AOT1} if we take into account the known identity \cite{Kolbig1996}
\begin{equation}
\psi \left(1, 3/4\right) = \pi^{2}-8\,G\,.
\label{PsiId}
\end{equation}
The remaining reductions are straightforward, involving the evaluation of known  digamma identities, based on Gauss' well-known general rules \cite[Eq. (2.21)]{Coffey2010}:
\begin{equation}
\psi \left(\frac{p}{q}\right) = 
-\gamma -\frac{\pi \,\cot \left(\frac{\pi \,p}{q}\right)}{2}-\ln \left(q \right)+2 \moverset{{\lfloor \frac{q}{2}\rfloor}}{\munderset{n =1}{\sum}}\cos \left(\frac{2\,\pi \,n\,p}{q}\right) \ln \left(2\,\sin \left(\frac{\pi \,n}{q}\right)\right)
\label{PpqO}
\end{equation}
and
\begin{align} \nonumber
\psi \left(\frac{p}{q}\right) =& 
\,-\gamma -\frac{\pi \,\cot \left(\frac{\pi \,p}{q}\right)}{2}-\ln \left(q \right)+2 \moverset{{\lfloor \frac{q}{2}\rfloor}-1}{\munderset{n =1}{\sum}}\cos \left(\frac{2\,\pi \,n\,p}{q}\right) \ln \left(2\,\sin \left(\frac{\pi \,n}{q}\right)\right)\\
&+\cos \left(\frac{2\,\pi \,{\lfloor \frac{q}{2}\rfloor}\,p}{q}\right) \ln \left(2\,\sin \left(\frac{\pi \,{\lfloor \frac{q}{2}\rfloor}}{q}\right)\right)
\label{PpqE}
\end{align}
according to whether $q$ is odd or even respectively. For the other special instances under discussion, we then have
 $$\psi \left(1/6\right)-\psi \left(2/3\right) = 
-2\,\ln \left(2\right)-\frac{2\,\pi \,\sqrt{3}}{3},
$$
\begin{equation}\psi \left(5/8\right)-\psi \left(1/8\right) = 
\left(\pi -2\,\ln \left(\sqrt{2}-1\right)\right) \sqrt{2}
\label{Psim58}
\end{equation}
and
$$\psi \left(1/12\right)-\psi \left(7/12\right) = 
-2\,\pi +2\,\sqrt{3}\,\ln \left(2-\sqrt{3}\right),$$
corresponding to the choices $y=1/3,y=1/4,y=1/6$ in \eqref{AOTP} respectively. 

For other choices of $y$, I submit that \eqref{AOTP} is a more fundamental expression than those obtained by both Asakura et. al. and Chen and Chu, since most expressions for the digamma function with general rational arguments will not reduce as simply as has been demonstrated here, and this property will be carried over into the tractability of their derivation methods. In turn, this observation suggests the following proposal:
\begin{definition}
The term ``exotic" in the sense introduced by Chu \cite{Chu700} as it pertains to specialized cases of 3F2(1), applies to all instances where the fractional parameters (e.g. $y= \pi p/q$) correspond to angles $y$ that are constructable using only straight edge and compass and therefore are amenable to analysis using \eqref{PpqO} and \eqref{PpqE}. Otherwise, see Murty and Saradha \cite{RamMurty}.
\end{definition}

\begin{rem}
In 2023, Chu redefined the word "exotic" (see \cite[page 2]{ChuAxioms12030291}) to essentially refer to the challenge problem, where the difference between a top and bottom parameter is a negative integer.
\end{rem} 

\begin{rem}
I have tested a limited number of the 700 examples presented in \cite{Chu700}, all of which turn out to be special cases contiguous to one of DWW and are therefore summable. It is important however, to ask what constitutes a representation of 3F2(1) in terms of fundamental constants? The exemplars listed above are about as fundamental as possible if $y$ is expressed in terms of an exceptional rational value; in other cases, one might consider that the digamma function (e.g. $\psi(y)$) is a fundamental constant for general values of $y$. Similarly, Qureshi and Shadab \cite[Eq. (12)]{QurSha} study problems of the form \eqref{Q3ca} and so rediscover a family of summations for 3F2(1) in terms of $\psi(p/q)$. Again, and only for specific values of $p$ and $q$, do they they utilize \eqref{PpqO} and \eqref{PpqE} to generate a collection of special cases  that might satisfy the definition of ``exotic" as proposed, omitting the fact that there exists an infinite number of values of $p/q$ where $\psi(p/q)$ does not resolve in terms of fundamental constants. In  such cases, an evaluation in terms of $\psi(p/q)$ would still be of considerable utility.
\end{rem}

\subsubsection{Chen and Chu Extended}

In an extension of their paper \cite[Section 3]{ChuChen1}, Chen and Chu consider contiguous instances of the function that has been considered in the previous section, and present a table of values. Specifically, with $m\leq n$, they study the function
$${}_{3}^{}{\moversetsp{}{\mundersetsp{}{F_{2}^{}}}}\! \left(y ,1-y ,x+m ;1  ,1+x+n;1\right).$$
By replacing $x:=x-m$  and $n:=n-m$, this extension can be reduced to the simpler form
$${}_{3}^{}{\moversetsp{}{\mundersetsp{}{F_{2}^{}}}}\! \left(y ,1-y ,x ;1  ,1+x+n;1\right),$$ 
a function that is again contiguous to Watson's theorem \eqref{Watson} with $n\geq -1$. Employing the same methods discussed previously, here we will set $x=1/2$ and find
\begin{equation}
{}_{3}^{}{\moversetsp{}{\mundersetsp{}{F_{2}^{}}}}\! \left(y ,1-y,1/2;1,5/2+n ;1\right)
 = 
\munderset{c \rightarrow \frac{1}{2}}{{\mathrm{lim}}}\! \textit{W}_{2n+3,0}  \left( y , 1-y , c\right)\,.
\label{CCExt}
\end{equation}
After some simplification involving evaluation of the limit, we eventually arrive at
\begin{align} \nonumber
{}_{3}^{}{\moversetsp{}{\mundersetsp{}{F_{2}^{}}}}&\left(y ,1-y ,1/2;1,5/2+n ;1\right)
 = 
\frac{4^{2\,n +3}  \sin \left(y\,\pi \right) \Gamma \left(\frac{5}{2}+n \right)^{2} }{\pi^{2}}\\ \nonumber
&\times\left(\left(y -n -2\right)\moverset{2\,n +3}{\munderset{k =0}{\sum}}\frac{\left(-1\right)^{k} \left(\psi \left(1-\frac{y}{2}-\frac{k}{2}\right)+\psi \left(\frac{y}{2}-n-1 +\frac{k}{2}\right)\right) \Gamma \left(\frac{y}{2}+\frac{k}{2}\right)^{2}\,\Gamma \left(1-2\,y -k \right)}{\Gamma \left(5-2\,y -k +2\,n \right) \Gamma \left(\frac{y}{2}-n -1+\frac{k}{2}\right)^{2}\,\Gamma \left(2\,n +4-k \right) \Gamma \left(k +1\right)}\right.\\
&\left.-\moverset{2\,n +2}{\munderset{k =0}{\sum}}\frac{\left(-1\right)^{k}\,\Gamma \left(-2\,y -k \right) \Gamma \left(\frac{1}{2}+\frac{y}{2}+\frac{k}{2}\right)^{2} \left(\psi \left(\frac{1}{2}-\frac{y}{2}-\frac{k}{2}\right)+\psi \left(\frac{y}{2}-n-\frac{1}{2} +\frac{k}{2}\right)\right)}{\Gamma \left(k +1\right) \Gamma \left(2\,n +3-k \right) \Gamma \left(4-2\,y -k +2\,n \right) \Gamma \left(\frac{y}{2}-\frac{1}{2}-n +\frac{k}{2}\right)^{2}}\right)\,.
\label{CC1B}
\end{align}
The case $n=0$ gives
\begin{align} \nonumber
{}_{3}^{}{\moversetsp{}{\mundersetsp{}{F_{2}^{}}}}&\left(y ,1-y,\frac{1}{2} ;1,\frac{5}{2};1\right)
 = 
\frac{3 }{\pi  \left(2\,y -3\right) \left(4\,y^{2}-1\right)}\\
&\times\left(\left(\left(2\,y^{2}-2\,y -1\right) \left(\psi \left(\frac{y}{2}\right)-\psi \left(\frac{1}{2}+\frac{y}{2}\right)\right)+2\,y -1\right) \sin \left(\pi y \right)+\pi  \left(2\,y^{2}-2\,y -1\right)\right)\,;
\label{CCn0}
\end{align}
the case $n=-1$ reproduces \eqref{AOTP} and previous comments apply. 

\subsubsection{Application to the challenge problem} 

Since the left-hand side of \eqref{CC1B} is a special case of the challenge problem \eqref{HAll2}, here we consider \eqref{HAll2} employing the parameters of \eqref{CC1B} and equate the right-hand sides of both, with $n=0$ as in \eqref{CCn0}. Solving yields the identity
\begin{equation}
{}_{3}^{}{\moversetsp{}{\mundersetsp{}{F_{2}^{}}}}\left(y ,y ,y -\frac{3}{2};2\,y ,\frac{1}{2}+y ;1\right)
 = 
\frac{\Gamma \left(\frac{1}{2}+y \right) 4^{y} \left(\left(\frac{1}{2}-y^{2}+y \right) \left(\psi \left(\frac{1}{2}+\frac{y}{2}\right)-\psi \left(\frac{y}{2}\right)\right)+y -\frac{1}{2}\right)}{\sqrt{\pi} \left(2\,y +1\right) \Gamma \left(y \right)},
\label{Hn0}
\end{equation}
the left-hand side of which is contiguous to Watson's identity \cite{ChuW}, with $m=0,n=3$. 

\subsubsection{Extension}
For $a,b,c,d\in \mathbb{Z}$, a recent paper by Chu \cite{ChuAxioms12030291}  evaluates a general form ${}_{3}^{}{\moversetsp{}{\mundersetsp{}{F_{2}^{}}}}\left(1+a ,c ,\frac{1}{2}+e ;\frac{3}{4}+b ,\frac{5}{4}+d ;1\right)$ by similar means and lists exotic examples for selected values of $a,b,c,d$. 
In a related work, Chen and Chu \cite{ChuChen2} consider the even more general case ${}_{3}^{}{\moversetsp{}{\mundersetsp{}{F_{2}^{}}}}\left(n,j ,k -y ;p +x ,q -x ;1\right)$ subdivided into four categories. If we set $x:=x-p$ and redefine $p+q=m$ , eliminating one parameter without loss of generality, it turns out that the case $k=0$ and $y=1/2$, equivalent to Chen and Chu's category A, is contiguous to Whipple's theorem. Specifically, with reference to Appendix \ref{sec:Contig}, if $m-n-j+1/2>0$,
\begin{equation}
{}_{3}^{}{\moversetsp{}{\mundersetsp{}{F_{2}^{}}}}\left(n ,j,-1/2 ;x,m -x  ;1\right)
 = 
\munderset{b \rightarrow -\frac{1}{2}}{\mathrm{lim}}\munderset{a \rightarrow n}{\mathrm{lim}}\Omega_{n+j-1,m} \left(a , b , m -x \right)\,,
\label{LjL}
\end{equation}
where the limit is required because \eqref{LjL} is an exceptional case. For example, if $j=1,n=2,m=3$, from \eqref{LjL} we find
\begin{align} \nonumber
{}_{3}^{}{\moversetsp{}{\mundersetsp{}{F_{2}^{}}}}&\left(1,2,-1/2;x ,3-x ;1\right)
 = 
\frac{4 \left(-1+x \right)^{2} \left(-2+x \right)^{2}\,\csc \left(\pi \,x \right) \pi}{\left(1-2\,x \right) \left(3-2\,x \right) \left(5-2\,x \right)}\\
&+\frac{4 \left(-2+x \right)^{2} \left(-1+x \right)^{2} \left(\psi \left(\frac{1}{2}+\frac{x}{2}\right)-\psi \left(\frac{x}{2}\right)\right)}{\left(-5+2\,x \right) \left(-1+2\,x \right) \left(-3+2\,x \right)}+\frac{2 \left(-2+x \right) \left(-1+x \right)}{\left(-5+2\,x \right) \left(-3+2\,x \right)}\,,
\label{LjL123}
\end{align}
and exotic instances will arise as discussed previously for selected rational values of $x$. Special cases respectively corresponding to $x=\{3/2,5/2,9/2,1/3\}$, as chosen by Chen and Chu, follow:
\begin{align}
&{}_{3}^{}{\moversetsp{}{\mundersetsp{}{F_{2}^{}}}}\left(-\frac{1}{2},1,2;\frac{3}{2},\frac{3}{2};1\right)
 = -\frac{\mathit{G}}{4}+\frac{3}{8},
\label{Excc2}\\
&{}_{3}^{}{\moversetsp{}{\mundersetsp{}{F_{2}^{}}}}\left(-\frac{1}{2},1,2;\frac{1}{2},\frac{5}{2};1\right)
 = -\frac{9\,\mathit{G}}{8}+\frac{3}{16},
\label{Excc1}\\
&{}_{3}^{}{\moversetsp{}{\mundersetsp{}{F_{2}^{}}}}\left(-\frac{1}{2},1,2;-\frac{3}{2},\frac{9}{2};1\right)
 = -{\frac{35}{9}},
\label{Excc4}\\
&{}_{3}^{}{\moversetsp{}{\mundersetsp{}{F_{2}^{}}}}\left(-\frac{1}{2},1,2;\frac{1}{3},\frac{8}{3};1\right)
 = -\frac{800}{273}\,\ln \left(2\right)+\frac{20}{91},
\label{Excc5}
\end{align}
where use has been made of \eqref{PsiId} and \eqref{Psim58} in the simplification process.

\subsection{The exotic identities of Chu and Wang} \label{sec:CWExotic}

In their 2007 work, Chu and Wang \cite{ChuWang2007} derived two and three part contiguity relations. Many of these relations are closely related to, but do not solve, the challenge problem. Their four fundamental contiguity relations are listed in Appendix \ref{sec:Contig} (see \eqref{Eq1A} - \eqref{Eq4A}). In their study, Chu and Wang \cite[Section 2 ff]{ChuWang2007} varied combinations of parameters and applied them to their four fundamental relations, yielding a collection of contiguity relations some of which they labelled ``exotic". To shed some light on the source of these exotic identities, let us simply specify one parameter as it affects one identity and consider the consequence using the challenge problem as a template. Specifically, in \eqref{Eq4A}, let $f\rightarrow a+1$, that is
\begin{align} \nonumber
&{}_{3}^{}{\moversetsp{}{\mundersetsp{}{F_{2}^{}}}}\left(a ,b ,c ;e ,a +1;1\right) = 
\munderset{f \rightarrow a +1}{\mathrm{lim}}\left\{\frac{1}{\left(a +1-f \right) \left(1+a -e \right)}\right. \\
&\left. \times\left[{\left(1-e \right) \left(1-f \right)\, _{3}F_{2}  \left(\right.\overset{ \displaystyle {a,b-1,c-1~}}{_ {\displaystyle {e-1,f-1~
 }}}|1\left.\right)+a \left(1+a +c +b -e -f \right)\, _{3}F_{2}  \left(\right.\overset{ \displaystyle {a+1,b,c~}}{_ {\displaystyle {e,f~
 }}}|1\left.\right)}\right]\right\}
\label{Ex4A}
\end{align}
and notice that the combination of terms enclosed in square brackets must vanish in that limit, a result easily verified because both reduce to 2F1(1). However, if we instead explicitly consider the series expansion of \eqref{Ex4A} before evaluating the limit, we find that a term of order $(a+1-f)^{-1}$ arises, which of course must vanish when the limit is applied. This produces a variety of ``strange" and ``exotic" identities (see below). But first, as part of  that expansion, we note that the term of order $(a+1-f)^{0}$ also yields an interesting identity: if $1+e-b-c>0$ then
\begin{align} \nonumber
&{}_{3}^{}{\moversetsp{}{\mundersetsp{}{F_{2}^{}}}}\left(a ,b ,c ;e ,a +1;1\right)
 = 
\frac{\left(c +b -e \right) a\,\Gamma \left(e \right)}{\Gamma \left(c \right) \Gamma \left(b \right) \left(1+a -e \right)}\moverset{\infty}{\munderset{k =0}{\sum}}\frac{\Gamma \left(c +k \right) \Gamma \left(b +k \right) \psi \left(a +k +1\right)}{\Gamma \left(k +1\right) \Gamma \left(e +k \right)}\\
&+\frac{a\,\Gamma \left(e \right) }{\left(1+a -e \right) \Gamma \left(c -1\right) \Gamma \left(b -1\right)}\moverset{\infty}{\munderset{k =0}{\sum}}\frac{\Gamma \left(c +k \right) \Gamma \left(b +k \right) \psi \left(a +k +1\right)}{\Gamma \left(k +2\right) \Gamma \left(e +k \right)}\\
&+\frac{a\,\Gamma \left(-c -b +e \right) \Gamma \left(e \right)}{\Gamma \left(e -c \right) \Gamma \left(e -b \right) \left(1+a -e \right)}+\frac{\left(e -1\right) a\,\psi \left(a \right)}{1+a -e},
\label{CWT2}
\end{align}
which might simplify further.

By straightforward machine evaluation, the vanishing coefficient of the term $(a+1-f)^{-1}$ then produces the exotic identity
\begin{equation}
{}_{4}^{}{\moversetsp{}{\mundersetsp{}{F_{3}^{}}}}\left(1,b ,c ,\frac{c\,b +b +c -2\,e +1}{c +b -e};2,e ,\frac{c\,b -e +1}{c +b -e};1\right)
 = \frac{e -1}{e-c\,b -1}\,,
\label{CWT1}
\end{equation}
which can be written in the slightly simplified, but no less exotic, form
\begin{equation}
{}_{4}^{}{\moversetsp{}{\mundersetsp{}{F_{3}^{}}}}\left(1,b ,c ,e ;2,\frac{\left(e-c  -1\right) b -1+\left(e -1\right) c}{e -2},e -1;1\right)
 = 
\frac{\left(c -e +1\right) b -\left(e -1\right) \left(c -1\right)}{\left(e -1\right) \left(c -1\right) \left(b -1\right)}
\label{CWT1A}
\end{equation}
by applying the change of variables $e := 
\left(\left(b +c \right) e-\left(c +1\right) \left(b +1\right)\right)/(e -2)
$.
Arising from a computational artifact, such forms emerge when two sums are combined rather than expanded. As a result, the identity \eqref{CWT1A} can be further simplified if we first expand a summation term of the form $$\moverset{\infty}{\munderset{k =0}{\sum}}f_{k}\Rightarrow\moverset{\infty}{\munderset{k =0}{\sum}}f_{k+1}-f_{-1}$$
to yield an identity equivalent to \eqref{CWT1A}
\begin{equation}
{}_{3}^{}{\moversetsp{}{\mundersetsp{}{F_{2}^{}}}}\left(a ,b ,c ;a -1,\frac{a \left(1+c +b \right)-\left(c +1\right) \left(b +1\right)}{a -1};1\right)
 = 0,
\label{Simp}
\end{equation}
after redefinition of variables. 
\begin{rem} \label{sec:Learn}
To understand how such exotic identities arise, consider a typical summation
\begin{equation}
\moverset{\infty}{\munderset{k =0}{\sum}}\frac{\left(a \right)_{k} \left(b \right)_{k} \left(\alpha \,k +\beta \right)}{\left(c \right)_{k}\,\Gamma \left(k +1\right)}
 = 
\alpha  \moverset{\infty}{\munderset{k =0}{\sum}}\frac{\left(a \right)_{k} \left(b \right)_{k}\,k}{\left(c \right)_{k}\,\Gamma \left(k +1\right)}+\beta  \moverset{\infty}{\munderset{k =0}{\sum}}\frac{\left(a \right)_{k} \left(b \right)_{k}}{\left(c \right)_{k}\,\Gamma \left(k +1\right)}
\label{Str1}
\end{equation}
that can be written in two different ways exemplified by the left and right-hand sides
\begin{equation}
\beta \,{}_{3}^{}{\moversetsp{}{\mundersetsp{}{F_{2}^{}}}}\left(a ,b ,\frac{\alpha +\beta}{\alpha};c ,\frac{\beta}{\alpha};1\right)
 = 
\frac{\alpha \,a\,b\,}{c}{}_{2}^{}{\moversetsp{}{\mundersetsp{}{F_{1}^{}}}}\left(b +1,a +1;c +1;1\right)+\beta \,{}_{2}^{}{\moversetsp{}{\mundersetsp{}{F_{1}^{}}}}\left(a ,b ;c ;1\right)\,,
\label{Str2}
\end{equation}
thereby demonstrating the genesis of ``strange" and ``exotic" identities of this genre. It is reminiscent of known identities that have been labelled ``strange"  in the literature (see \cite[refs. 3 and 4]{Milgram447}, for example \cite[entry 91]{Milgram447}):
\begin{align} 
{}_{3}^{}{\moversetsp{}{\mundersetsp{}{F_{2}^{}}}}\left(a ,b ,c ;a +2,\frac{\left(a +1\right) \left(b +c -1\right)-c\,b}{a};1\right)
 = 
\frac{\left(a +1\right) \Gamma \left(\frac{\left(a +1\right) \left(b +c -1\right)-c\,b}{a}\right) \Gamma \left(\frac{\left(c -1\right) \left(1-b \right)+a}{a}\right)}{\Gamma \left(\frac{a\,c +\left(b -1\right) \left(-c +1\right)}{a}\right) \Gamma \left(\frac{a\,b +\left(b -1\right) \left(-c +1\right)}{a}\right)}
\label{E91}
\end{align}
due to Chu and Gessel and Stanton. 
\end{rem}

\subsection{Campbell and Abrarov }
We continue with Campbell and Abrarov \cite{CampAbr} where many of the principles discussed previously can be demonstrated.
\subsubsection{Theorem 8} \label{sec:CandAT8}
 In their paper, Campbell and Abrarov obtain the identity \cite[Eqs. (16) and (17)]{CampAbr}
\begin{equation}
\pi \,\ln \! \left(2\right)-\frac{16}{15}\,{}_{3}^{}{\moversetsp{}{\mundersetsp{}{F_{2}^{}}}}\! \left(1,1,\frac{3}{2};\frac{7}{4},\frac{9}{4};1\right)
 = 
-\moverset{\infty}{\munderset{n =1}{\sum}}\! \frac{\left({\mstack{{1}/{2}}{n}}\right) 2^{-n}\,\pi \,\Gamma \! \left(n +1\right) \left(h({\frac{n}{2}})+2\,\ln \! \left(2\right)\right)}{2\,\Gamma \! \left(\frac{n}{2}+1\right)^{2}}
\label{Ca16}
\end{equation}
where $h(\frac{n}{2})$ is the harmonic series defined by
\begin{equation}
h(\frac{n}{2}) = \psi \! \left(\frac{n}{2}+1\right)+\gamma \,,
\label{Harm}
\end{equation}
and lament that ``Mathematica is not able to evaluate the $_{3}F_{2}$ series...". It will now be shown that that particular $_{3}F_{2}$ is a special case of a known identity and therefore evaluable, in which case a simplified form of the series can be identified. However, we first note that the right-hand side of \eqref{Ca16} can we expanded as
\begin{align} \nonumber
\moverset{\infty}{\munderset{n =1}{\sum}}\! &\frac{\left({\mstack{{1}/{2}}{n}}\right) 2^{-1-n}\,\pi \,\Gamma \! \left(n +1\right) \left(h(\frac{n}{2})+2\,\ln \! \left(2\right)\right)}{\Gamma \! \left(\frac{n}{2}+1\right)^{2}}
 \\&= 
\frac{\pi^{\frac{3}{2}}}{4} \moverset{\infty}{\munderset{n =1}{\sum}}\! \frac{2^{-n}\,h({\frac{n}{2}})}{\Gamma \left(\frac{3}{2}-n \right) \Gamma \left(\frac{n}{2}+1\right)^{2}}
+\frac{\pi^{\frac{3}{2}}\,\ln \! \left(2\right)}{2}\moverset{\infty}{\munderset{n =1}{\sum}}\! \frac{2^{-n}}{\Gamma \left(\frac{3}{2}-n \right) \Gamma \left(\frac{n}{2}+1\right)^{2}}
\label{Ca16A}
\end{align}
and, courtesy of Maple,
\begin{equation}
\frac{\pi^{\frac{3}{2}}\,\ln \! \left(2\right)}{2} \moverset{\infty}{\munderset{n =1}{\sum}}\! \frac{2^{-n}}{\Gamma \left(\frac{3}{2}-n \right) \Gamma \left(\frac{n}{2}+1\right)^{2}}
 = \left(4-\pi \right) \ln \! \left(2\right)\,.
\label{Msum}
\end{equation}
Putting it all together along with \eqref{Harm} eventually yields the identity

\begin{equation}
{}_{3}^{}{\moversetsp{}{\mundersetsp{}{F_{2}^{}}}}\! \left(1,1,\frac{3}{2};\frac{7}{4},\frac{9}{4};1\right)
 = 
\frac{15\,\pi^{\frac{3}{2}}}{64} \moverset{\infty}{\munderset{n =1}{\sum}}\! \frac{2^{-n}\,\psi \left(\frac{n}{2}+1\right)}{\Gamma \left(\frac{3}{2}-n \right) \Gamma \left(\frac{n}{2}+1\right)^{2}}-\frac{15 \left(\pi -4\right) \gamma}{32}+\frac{15\,\ln \! \left(2\right)}{4}\,.
\label{Ca16c}
\end{equation}
We now note that the hypergeometric function on the left-hand side of \eqref{Ca16c} is a special case contiguous to Whipple's theorem, specifically

\begin{equation}
{}_{3}^{}{\moversetsp{}{\mundersetsp{}{F_{2}^{}}}}\! \left(1,1,\frac{3}{2};\frac{7}{4},\frac{9}{4};1\right)
 =\munderset{b \rightarrow \frac{3}{2}}{\mathrm{lim}}\munderset{a \rightarrow 1}{\mathrm{lim}}\Omega_{1,0}(a,b ,7/4)
\label{Campb1}
\end{equation}
as defined in \eqref{Whipple}. Since a closed form, finite summation for the right-hand side of \eqref{Campb1} is available in \cite{ChuW} after evaluating the limits, we obtain

\begin{equation}
{}_{3}^{}{\moversetsp{}{\mundersetsp{}{F_{2}^{}}}}\! \left(1,1,\frac{3}{2};\frac{7}{4},\frac{9}{4};1\right)
 = 
-\frac{15\,\sqrt{2}\,\mathrm{arctanh}\! \left(\frac{2}{3}\sqrt{2}\right)}{8}+\frac{15}{2}\,.
\label{CaR2}
\end{equation}
in agreement with the result quoted by Campbell and Abrarov \cite[Theorem 8]{CampAbr} who utilize Mathematica to evaluate an equivalent integral. Some reordering and simplification finally yields the Euler-type identity
\begin{equation}
\moverset{\infty}{\munderset{n =1}{\sum}}\! \frac{2^{-n}\,\psi \! \left(\frac{n}{2}+1\right)}{\Gamma \! \left(\frac{3}{2}-n \right) \Gamma \! \left(\frac{n}{2}+1\right)^{2}}
 = 
\frac{-8\,\sqrt{2}\,\ln \! \left(3+2\,\sqrt{2}\right)+32+2\,\pi \,\gamma -8\,\gamma -16\,\ln \! \left(2\right)}{\pi^{\frac{3}{2}}}\,,
\label{Ca16d}
\end{equation}
which can be written in the generic Euler series form

\begin{align} \nonumber
\sqrt{\pi} \moverset{\infty}{\munderset{n =0}{\sum}}\! \frac{2^{-2\,n}\,\psi \! \left(n +\frac{1}{2}\right) \Gamma \! \left(2\,n -\frac{3}{2}\right)}{\Gamma \! \left(n +\frac{1}{2}\right)^{2}}
 &= 
2 \gamma\left(\sqrt{2}-\frac{8}{3}\right)  +4 \left(3\,\sqrt{2}-\frac{8}{3}\right) \ln \! \left(2\right)\\
&-4\,\sqrt{2} \left(\ln \! \left(3+2\,\sqrt{2}\right)+2\right)+16
\label{Ca16e}
\end{align}
by dissecting the left-hand side of \eqref{Ca16d} into even and odd indices. See Appendix \ref{sec:Ca16eProof}.

\begin{rem}
Why Maple (``simplify") and Mathematica (``FullSimplify") both fail the evaluation of the hypergeometric function on the left-hand side of \eqref{CaR2} according to the most modern WZ methods (see \cite{AeqB}) is likely due to the fact that careful limits are required and therefore digamma functions are embedded in the final evaluations (viz. {\it arctanh} is a limiting digamma value). The presence of digamma functions means that corresponding expressions will not satisfy the criterion to be categorized as ``hypergeometric identities" and therefore related WZ algorithms are invalid.

\end{rem}
\subsubsection{Theorem 15} \label{sec:CandAT15}
With their Theorem 15, Campbell and Abrarov explore ``the evaluation of 3F2(1) series that would otherwise seem to have no possible closed-form expression". Their identity

\begin{equation}
\frac{{}_{3}^{}{\moversetsp{}{\mundersetsp{}{F_{2}^{}}}}\! \left(\frac{3}{2},\frac{1}{2}-a ,\frac{3}{2}-a ;1,\frac{5}{2}-a ;1\right)}{{}_{3}^{}{\moversetsp{}{\mundersetsp{}{F_{2}^{}}}}\! \left(-\frac{1}{2},a ,1-a ;\frac{1}{2},\frac{3}{2};1\right)}
 = 
\frac{4\,\Gamma \! \left(\frac{5}{2}-a \right) \Gamma \! \left(a \right)}{\pi^{\frac{3}{2}}}
\label{Car23}
\end{equation}
provides one such interesting result, more so if it is rewritten in the form of a two-part transformation 
\begin{equation}
{}_{3}^{}{\moversetsp{}{\mundersetsp{}{F_{2}^{}}}}\! \left(\frac{3}{2},\frac{1}{2}-a ,\frac{3}{2}-a ;1,\frac{5}{2}-a ;1\right)
 = 
\frac{4\,\Gamma \! \left(\frac{5}{2}-a \right) \Gamma \! \left(a \right)}{\pi^{\frac{3}{2}}} {}_{3}^{}{\moversetsp{}{\mundersetsp{}{F_{2}^{}}}}\! \left(-\frac{1}{2},a ,1-a ;\frac{1}{2},\frac{3}{2};1\right)\,.
\label{Car23a}
\end{equation}
By inspection, we discover that \eqref{Car23a} is a special case of Thomae identity \eqref{Thom5}. Furthermore, both of the 3F2(1) that appear in \eqref{Car23a} are special cases of known results listed by Chu \cite{ChuW}. Specifically, the left-hand side is a special case contiguous to Watson's theorem \eqref{Watson} (see \cite[Entry 5]{Milgram447}) with $m=1,~n=2,~c=3/2,b\rightarrow 3/2-a$ and $a:=1/2-a$ so that

\begin{align} \nonumber
&{}_{3}^{}{\moversetsp{}{\mundersetsp{}{F_{2}^{}}}}\! \left(\frac{3}{2},\frac{1}{2}-a ,\frac{3}{2}-a;1,\frac{5}{2}-a ;1\right)
 = \munderset{b \rightarrow \frac{3}{2}-a}{\mathrm{lim}}W_{2 ,-2}\! \left(a , b, 3/2\right)\\
&=\frac{\left(2\,a -1\right)  \Gamma \! \left(\frac{5}{2}-a \right) \Gamma \! \left(a \right)}{\sqrt{\pi}}\left(\frac{ \cos \left(\pi \,a \right)}{\pi}\left(\psi \left(\frac{3}{4}-\frac{a}{2}\right)-\psi \left(\frac{5}{4}-\frac{a}{2}\right)-\frac{4\,a}{\left(2\,a -1\right)^{2}}\right)-1\right)\,.
\label{CaR3b}
\end{align}
The denominator function can then be found from either \eqref{Car23a} or from the observation that it is a special case of Whipples Theorem \eqref{Whipple} (see \cite[entry 13]{Milgram447}) to eventually arrive at
\begin{align}\nonumber
&{}_{3}^{}{\moversetsp{}{\mundersetsp{}{F_{2}^{}}}}\! \left(-\frac{1}{2},a ,-a +1;\frac{1}{2},\frac{3}{2};1\right)
 =\munderset{b \rightarrow -\frac{1}{2}}{\mathrm{lim}} \Omega_{0,2}\! \left(a:=1-a , b, 1/2\right)\\
&=
\frac{\left(a -\frac{1}{2}\right)}{2} \left(\left(\psi \! \left(\frac{1}{4}+\frac{a}{2}\right)-\psi \! \left(\frac{3}{4}+\frac{a}{2}\right)+\frac{a -1}{\left(a -\frac{1}{2}\right)^{2}}\right) \cos \! \left(\pi \,a \right)+\pi \right)\,.
\label{CaR3a}
\end{align}
Substituting $a=1/4$ in \eqref{Car23a} will reproduce an identity equivalent to Campbell and Abrarov Corollary (4) and sheds light on their following comments; substitution of $a=3/4$ in \eqref{CaR3b} will reproduce their Corollary (5).

\subsubsection{Theorem 16} \label{sec:CandAT16}

Similar to the previous section, Campbell and Abrarov's Theorem (16) is equivalent to 
\begin{equation}
{}_{3}^{}{\moversetsp{}{\mundersetsp{}{F_{2}^{}}}}\! \left(\frac{5}{2},\frac{1}{2}-\alpha ,\frac{5}{2}-\alpha ;2,\frac{7}{2}-\alpha ;1\right)
=
\frac{4\, \left(\frac{5}{2}-\alpha \right) \Gamma \! \left(\alpha \right) \Gamma \! \left(\frac{3}{2}-\alpha \right)}{3 \left(\frac{1}{2}+\alpha \right) \pi^{\frac{3}{2}}}{}_{3}^{}{\moversetsp{}{\mundersetsp{}{F_{2}^{}}}}\! \left(-\frac{3}{2},\alpha ,1-\alpha ;-\frac{1}{2},\frac{3}{2};1\right)
\label{Thm16}
\end{equation}
and, surprisingly, \eqref{Thm16} is not reproduced by any of the independent Thomae relations presented in Appendix \ref{sec:Thomae} when applied to the left-hand side. However, if the Thomae relation \eqref{Thom1} is applied to the right-hand side of \eqref{Thm16}, an exceptional case that requires the invocation of a limiting calculation arises. Specifically,

\begin{equation}
{}_{3}^{}{\moversetsp{}{\mundersetsp{}{F_{2}^{}}}}\! \left(-\frac{3}{2},\alpha ,1-\alpha ;-\frac{1}{2},\frac{3}{2};1\right)
 = 
\frac{ \Gamma \! \left(\frac{3}{2}\right)^{2}\,\Gamma \! \left(-\frac{1}{2}\right)}{\Gamma \! \left(\alpha \right)\,\Gamma \! \left(0\right)\,\Gamma \! \left(\frac{5}{2}-\alpha \right)}{}_{3}^{}{\moversetsp{}{\mundersetsp{}{F_{2}^{}}}}\! \left(\frac{3}{2},-\frac{1}{2}-\alpha ,\frac{3}{2}-\alpha ;0,\frac{5}{2}-\alpha ;1\right)\,,
\label{Th1}
\end{equation}
where the presence of ``$0$" in the fourth parameter of the 3F2(1) and $\Gamma(0)$ in the coefficient, indicates that a limit is required. In this particular example, applying \eqref{Thom1} to the left-hand side of \eqref{Th1} with the fourth parameter $-1/2:=e$, then evaluating the limit $e\rightarrow -1/2$ to the right-hand side, will reproduce the equivalent of \eqref{Thm16} and so, in this sense, \eqref{Thm16} does obey a Thomae relation. The details are left as an exercise for the reader.

Considering its individual components, we find that the left-hand side of \eqref{Thm16} is contiguous to Watson's theorem (\cite[entry (6)]{Milgram447} with $m=1,n=3$, giving, after simplification and invocation of limits
\begin{align} \nonumber
&{}_{3}^{}{\moversetsp{}{\mundersetsp{}{F_{2}^{}}}}\! \left(\frac{5}{2},\frac{1}{2}-a ,\frac{5}{2}-a ;2,\frac{7}{2}-a ;1\right)=\munderset{b \rightarrow \frac{5}{2}-a}{\mathrm{lim}} W_{3,-3}\! \left(a:=\frac{1}{2}-a , b , c=\frac{5}{2}\right)\\
&= 
-\frac{2 \left(a -\frac{1}{2}\right)  \Gamma \! \left(a \right)}{3\,\Gamma \! \left(-\frac{5}{2}+a \right) \sqrt{\pi}}\left(\frac{\pi}{\cos \left(\pi \,a \right)}+\psi \! \left(\frac{1}{4}+\frac{a}{2}\right)-\psi \! \left(\frac{3}{4}+\frac{a}{2}\right)+\frac{-2\,a^{2}+\frac{1}{4}\,a +a^{3}+1}{\left(a -\frac{1}{2}\right)^{2} \left(a -\frac{3}{2}\right) \left(a +\frac{1}{2}\right)}\right)
\label{Cf3d}
\end{align}
reducing to 
\begin{equation}
{}_{3}^{}{\moversetsp{}{\mundersetsp{}{F_{2}^{}}}}\! \left(\frac{1}{4},\frac{9}{4},\frac{5}{2};2,\frac{13}{4};1\right)
= 
\frac{15\,\sqrt{\pi}\,\ln \! \left(17-12\,\sqrt{2}\right)}{128\,\Gamma \! \left(\frac{3}{4}\right)^{2}}+\frac{61\,\sqrt{\pi}\,\sqrt{2}}{32\,\Gamma \! \left(\frac{3}{4}\right)^{2}}
\label{Tnq}
\end{equation}
when $a=1/4$, equivalent to Campbell and Abrarov Corollary (6). Similarly, the right-hand side function is available from either \eqref{Thm16} and \eqref{Tnq} or from Whipple's theorem \cite[entry (13)]{Milgram447} with $m=0,n=3$, yielding
\begin{align}\nonumber
&{}_{3}^{}{\moversetsp{}{\mundersetsp{}{F_{2}^{}}}}\! \left(1-a ,-\frac{3}{2},a ;-\frac{1}{2},\frac{3}{2};1\right)= \munderset{b \rightarrow -\frac{3}{2}}{\mathrm{lim}}\Omega_{0,3}\! \left(a:=1-a,b, -\frac{1}{2}\right)\\ \nonumber&
 = 
\frac{1}{2}\left(a -\frac{3}{2}\right) \left(-\frac{1}{4}+a^{2}\right) \left(\left(\psi \! \left(\frac{3}{4}+\frac{a}{2}\right)-\psi \! \left(\frac{1}{4}+\frac{a}{2}\right)\right) \cos \! \left(\pi \,a \right)-\pi \right)\\
&+\frac{\left(-4\,a^{3}+8\,a^{2}-a -4\right) }{8\,a -4}\cos \! \left(\pi \,a \right)\,,
\label{CAd2}
\end{align}
reducing to 
\begin{equation}
{}_{3}^{}{\moversetsp{}{\mundersetsp{}{F_{2}^{}}}}\! \left(-\frac{3}{2},\frac{1}{4},\frac{3}{4};-\frac{1}{2},\frac{3}{2};1\right)
 = 
-\frac{15\,\ln \! \left(1+\sqrt{2}\right)}{64}+\frac{61\,\sqrt{2}}{64}
\label{C16d}
\end{equation}
when $a=1/4$.

\subsubsection{Theorem 17} \label{sec:CandAT17}

In the same manner as the previous sections, we find that Campbell and Abrarov Theorem (17), written as
\begin{equation}
{}_{3}^{}{\moversetsp{}{\mundersetsp{}{F_{2}^{}}}}\! \left(1-a,a ,b ;\frac{3}{2},b +1;1\right)
 = 
\frac{2\,b\,\cos \! \left(\pi \,a \right) {}_{3}^{}{\moversetsp{}{\mundersetsp{}{F_{2}^{}}}}\! \left(1,\frac{3}{2},\frac{3}{2}-b ;\frac{5}{2}-a ,a +\frac{3}{2};1\right)}{\left(2\,a -3\right) \left(2\,a -1\right) \left(2\,a +1\right)}
\label{CaR17}
\end{equation}
corresponds to the Thomae relation \eqref{Thom2}. If we set $a=b=1/4$ we identify the corresponding 3F2(1) from \cite[Entry 17]{Milgram447} as being  contiguous to Dixon's theorem \eqref{Dixon}, specifically
\begin{equation}
{}_{3}^{}{\moversetsp{}{\mundersetsp{}{F_{2}^{}}}}\! \left(1,\frac{5}{4},\frac{3}{2};\frac{7}{4},\frac{9}{4};1\right)
 = \munderset{c \rightarrow 1}{\mathrm{lim}}X_{1,1}\! \left(\frac{5}{4}, \frac{3}{2}, c\right)=
-\frac{15}{2}-\frac{15}{8}\,\left(\psi \! \left(\frac{1}{8}\right)-\psi \! \left(\frac{5}{8}\right)\right),
\label{CaR17c}
\end{equation}
equivalent to Campbell and Abrarov Corollary (7) with the use of \eqref{Psim58}.

\begin{rem}
Further, although neither of the 3F2(1) appearing in \eqref{CaR17} can be evaluated from the database \cite{Milgram447}, in the case that $b=n$, a non-negative integer, we find, from \cite[entry 26]{Milgram447}, the identity

\begin{align} \nonumber
&{}_{3}^{}{\moversetsp{}{\mundersetsp{}{F_{2}^{}}}}\! \left(a ,n ,1-a ;\frac{3}{2},n +1;1\right)
 = 
 \frac{\Gamma \! \left(-a -n +1\right) \Gamma \! \left(n +1\right) \sqrt{\pi}}{2\,\Gamma \! \left(n +1-a \right) \Gamma \! \left(\frac{3}{2}-n \right)}\\
&+
\frac{\Gamma \! \left(a -n \right) \Gamma \! \left(n +1\right) \pi}{4\,\Gamma \! \left(-a +1\right) \Gamma \! \left(a +\frac{1}{2}\right)} \moverset{n -1}{\munderset{k =0}{\sum}}\! \frac{\Gamma \left(-k -a \right) \left(-1\right)^{k}}{\Gamma \left(-a +\frac{3}{2}-k \right) \Gamma \left(n -k \right) \Gamma \left(a +1+k -n \right)},
\label{CaR17a}
\end{align}
thereby reducing an infinite sum to a terminating one. Alternatively, if $a=m$, a positive integer, we identify, from \cite[entry (33)]{Milgram447}

\begin{equation}
{}_{3}^{}{\moversetsp{}{\mundersetsp{}{F_{2}^{}}}}\! \left(m ,b ,1-m ;\frac{3}{2},b +1;1\right)
 = 
\frac{\left(-1\right)^{m}\,b  \Gamma \! \left(b -\frac{1}{2}\right) \Gamma \! \left(m \right)}{4\,\Gamma \! \left(b +m \right) \Gamma \! \left(\frac{1}{2}+m \right)}\moverset{m -1}{\munderset{k =0}{\sum}}\! \frac{\Gamma \left(k -\frac{1}{2}\right) \Gamma \left(b +k \right)}{\Gamma \left(k +1\right) \Gamma \left(\frac{1}{2}+b -m +k \right)}
\label{CaR17b}
\end{equation}
and corresponding identities from \eqref{CaR17}. 
\end{rem}

\subsection{Mirzoev and Safonova} \label{sec:MzSf}

In a lengthy analysis stretching many pages, Mirzoev and Safonova \cite{MirzSaf} evaluate two sums of interest
$${}_{3}^{}{\moversetsp{}{\mundersetsp{}{F_{2}^{}}}}\left(\frac{1}{2}-\frac{a}{2},\frac{a}{2}+\frac{1}{2},\frac{1}{2};1,\frac{3}{2};1\right)~
\mathrm{and}~~
{}_{3}^{}{\moversetsp{}{\mundersetsp{}{F_{2}^{}}}}\left(1-a ,1+a ,\frac{1}{2};2,\frac{3}{2};1\right)$$
employing integral representations and special sums. Since generalizations of both of these sums are contiguous to Watson's theorem \cite{ChuW}, we can extend Mirzoev and Safonova's solution by writing 
\begin{align} \nonumber
&{}_{3}^{}{\moversetsp{}{\mundersetsp{}{F_{2}^{}}}}\left(\frac{1}{2}-\frac{a}{2},\frac{a}{2}+\frac{1}{2},c ;2\,c ,\frac{3}{2};1\right)
 = 
W_{1,0} \left( \frac{1}{2}-\frac{a}{2}, \frac{a}{2}+\frac{1}{2}, c\right)\\
&=\frac{2\,\Gamma \left(c +\frac{1}{2}\right)^{2} }{a \left(2\,c -1\right)}\left(-\frac{\cos \left(\frac{\pi  \left(1+a \right)}{4}\right)}{\Gamma \left(\frac{1}{4}-\frac{a}{4}+c \right) \Gamma \left(-\frac{1}{4}+\frac{a}{4}+c \right)}+\frac{\sin \left(\frac{\pi  \left(1+a \right)}{4}\right)}{\Gamma \left(\frac{1}{4}+\frac{a}{4}+c \right) \Gamma \left(-\frac{1}{4}-\frac{a}{4}+c \right)}\right)
\label{Ms1a}
\end{align} 

and

\begin{align} \nonumber
&{}_{3}^{}{\moversetsp{}{\mundersetsp{}{F_{2}^{}}}}\left(1-a ,1+a ,c ;2\,c +1,\frac{3}{2};1\right)
 = W_{0,1} \left(1 -a , 1+a , c\right)\\
 &=\frac{\Gamma \left(c +\frac{1}{2}\right)^{2} }{2\,c -1}\left(-\frac{\cos \left(\frac{\pi \,a}{2}\right)}{\Gamma \left(\frac{1}{2}+c +\frac{a}{2}\right) \Gamma \left(\frac{1}{2}+c -\frac{a}{2}\right)}+\frac{2\,\sin \left(\frac{\pi \,a}{2}\right)}{a\,\Gamma \left(\frac{a}{2}+c \right) \Gamma \left(-\frac{a}{2}+c \right)}\right)\,.
\label{Ms2a}
\end{align}
It is a straightforward task to now evaluate the limit $c\rightarrow 1/2$ in both \eqref{Ms1a} and \eqref{Ms2a} to arrive at 
\begin{equation}
{}_{3}^{}{\moversetsp{}{\mundersetsp{}{F_{2}^{}}}}\left(\frac{1}{2},-\frac{a}{2}+\frac{1}{2},\frac{a}{2}+\frac{1}{2};1,\frac{3}{2};1\right)
 = 
\frac{1}{a}+\frac{\sin \left(\frac{\pi  \left(1+a \right)}{2}\right) }{\pi \,a}\left(\psi \left(\frac{1}{4}+\frac{a}{4}\right)-\psi \left(\frac{3}{4}+\frac{a}{4}\right)\right)
\label{Ms1b}
\end{equation}
and
\begin{equation}
{}_{3}^{}{\moversetsp{}{\mundersetsp{}{F_{2}^{}}}}\left(\frac{1}{2},1+a ,-a +1;\frac{3}{2},2;1\right)
 = \frac{1}{a}+
\frac{ \sin \left(\pi \,a \right)}{a^{2}\,\pi}\left(\psi \left(\frac{a}{2}\right) a -\psi \left(\frac{a}{2}+\frac{1}{2}\right) a +1\right)\,,
\label{Ms2b}
\end{equation}
after which setting $a=p/q$ for selected values of $p$ and $q$ will yield exotic identities similar to those discussed previously that are equivalent to those obtained Mirzoev and Safonova. It is interesting to speculate whether an equivalent analysis exists for two sums 4F3(1) also considered by Mirzoev and Safonova.

\subsection{Zaidi and Almuthaybiri }

A recent paper by Zaidi and Almuthaybiri \cite{ZaidiAbdel} provides, among other things, an educational example that demonstrates how a complex derivation can be simplified. In about 15 pages, the authors evaluate three instances of 3F2(1):
\begin{align}
&K(a)\equiv {}_{3}^{}{\moversetsp{}{\mundersetsp{}{F_{2}^{}}}}\left(1-a,1,1+a ,;1+a,a+2 ;1\right)\,,\hspace{20pt}a>1/2;
\label{K}
\\
&G(a)\equiv {}_{3}^{}{\moversetsp{}{\mundersetsp{}{F_{2}^{}}}}\left(1-a,1,2+a ,;1+a,3+a ;1\right)\,,\hspace{20pt}a>1/2;
\label{Gz}\\
&H(a)={}_{3}^{}{\moversetsp{}{\mundersetsp{}{F_{2}^{}}}}\left(1/2-n,1,2\,n +1 ;3/2+n,2\,n +2 ;1\right)\,,\hspace{20pt} n=0,1\dots.
\label{Za5}
\end{align}
Some comments follow:
\begin{itemize}
\item{} The specified limits in \eqref{K} and \eqref{Gz} should be $a>0$;
\item{}
In \eqref{K}, since the parameter $c=e$ (see \eqref{3F2Def}), the 3F2(1) trivially reduces to 2F1(1) and immediately reproduces Zaidi and Almuthaybiri's result
\begin{equation}
{}_{3}^{}{\moversetsp{}{\mundersetsp{}{F_{2}^{}}}}\left(1-a,1,1+a ;1+a ,2+a ;1\right)
 = 
{}_{2}^{}{\moversetsp{}{\mundersetsp{}{F_{1}^{}}}}\left(1,1-a ;2+a ;1\right)
 = \frac{1+a}{2\,a} \, ,a > 0,
\label{ZaK}
\end{equation}
from a well-known theorem due to Gauss (1812);
\item
The function $G(a)$ in \eqref{Gz} satisfies  the Karlsson identity \eqref{Karl8} with $n=1$,  immediately reproducing, in one line, the equivalent of Zaidi and Almuthaybiri's result 
\begin{equation}
{}_{3}^{}{\moversetsp{}{\mundersetsp{}{F_{2}^{}}}}\left(1,2+a ,1-a ;a +3,1+a ;1\right)
 = 
\frac{\left(2+a \right) \left(2\,a^{2}+a +1\right)}{2\,a \left(2\,a +1\right) \left(1+a \right)}\,;
\label{ZaG}
\end{equation}

\item{}
The  function $H(a)$ in \eqref{Za5} happens to be contiguous to Dixon's theorem \eqref{Dixon}, but a simpler evaluation can be found by first considering its Thomae progeny and noticing that
\begin{align} \nonumber
{}_{3}^{}{\moversetsp{}{\mundersetsp{}{F_{2}^{}}}}&\left(\frac{1}{2}-n,1,2\,n +1 ;\frac{3}{2}+n,2\,n +2 ;1\right)\\&
 = 
\frac{ \Gamma \left(2\,n +1\right) \Gamma \left(2\,n +2\right)}{\Gamma \left(\frac{3}{2}+3\,n \right) \Gamma \left(\frac{3}{2}+n \right)}
\, {}_{3}^{}{\moversetsp{}{\mundersetsp{}{F_{2}^{}}}}\left(\frac{1}{2}+n ,\frac{1}{2}-n ,\frac{1}{2}-n ;\frac{3}{2}+n ,\frac{3}{2}+n ;1\right)
\label{Za8}
\end{align}
after applying the Thomae operator \eqref{Thom8}. The function on the right-hand side of \eqref{Za8} is a special case of the identity \eqref{L2p2}, yielding
\begin{align} \nonumber
&{}_{3}^{}{\moversetsp{}{\mundersetsp{}{F_{2}^{}}}}\left(\frac{1}{2}+n ,\frac{1}{2}-n ,\frac{1}{2}-n ;\frac{3}{2}+n ,\frac{3}{2}+n ;1\right)\\ \nonumber
& = 
\munderset{h \rightarrow \frac{1}{2}}{\mathrm{lim}}\left(\frac{\left(2\,n +1\right) \Gamma \left(\frac{1}{2}-2\,n -h \right)  \Gamma \left(\frac{3}{2}+3\,n \right) \Gamma \left(h +n +1\right) \Gamma \left(\frac{3}{2}+n \right)}{\Gamma \left(\frac{1}{2}-n \right) \Gamma \left(2\,n +2\right)}\right. \\ \nonumber & \times\left.\moverset{2\,n}{\munderset{k =0}{\sum}}\frac{\Gamma \left(-k -n -\frac{1}{2}\right) \left(-1\right)^{k}}{\Gamma \left(h +n -k \right) \Gamma \left(2\,n +1-k \right) \Gamma \left(\frac{3}{2}-2\,n -h +k \right)}
\right. \\ &\left.
+\frac{\Gamma \left(-\frac{1}{2}-3\,n \right) \Gamma \left(h -\frac{1}{2}\right) \Gamma \left(\frac{3}{2}+3\,n \right) \Gamma \left(h +n +1\right) \Gamma \left(\frac{3}{2}+n \right)}{\Gamma \left(h +2\,n +\frac{1}{2}\right)^{2}\,\Gamma \left(\frac{1}{2}-n \right)^{2}}\right)
\label{Za1}
\end{align}
where again a limit is required because this is an exceptional case. It is a straightforward, computer-aided exercise, to evaluate the limit in \eqref{Za1} and eventually obtain a simpler form of the general result sought by Zaidi and Almuthaybiri:
\begin{align} \nonumber
{}_{3}^{}{\moversetsp{}{\mundersetsp{}{F_{2}^{}}}}&\left(\frac{1}{2}-n,1,2\,n +1 ;\frac{3}{2}+n,2\,n +2 ;1\right)
 \\\nonumber
 &\, 
=-\frac{2 \left(-1\right)^{n}\,\Gamma \left(\frac{3}{2}+n \right)^{2} }{\pi}\moverset{2\,n -1}{\munderset{k =0}{\sum}}\frac{\Gamma \left(k +\frac{1}{2}-3\,n\right) \Gamma \left(k +1\right)}{\Gamma \left(\frac{3}{2}-n +k \right) \Gamma \left(2+k \right)}\\
&-\frac{\Gamma \left(\frac{3}{2}+n \right)^{3} \left(-1\right)^{n} }{\pi  \left(\frac{1}{2}+n \right) \Gamma \left(\frac{3}{2}+3\,n \right)}\left(\psi \left(1/2+n \right)-\psi \left(1+n \right)-2\,\ln \left(2\right)\right)\,.
\label{za5}
\end{align}
The identity \eqref{za5} numerically reproduces the results quoted by Zaidi and Almuthaybiri for the cases $n=1,~n=2$, but disagrees for the case $n=3$, where we find the correct identity
\begin{equation}
{}_{3}^{}{\moversetsp{}{\mundersetsp{}{F_{2}^{}}}}\left(1,7,-\frac{5}{2};8,\frac{9}{2};1\right)
 = \frac{518231}{831402}-\frac{70}{138567}\,\ln \left(2\right),
\label{Za5n3}
\end{equation}
a result verified numerically;
\item{}
It should be pointed out that the authors' reference 25 cites a paper that has long been known (see \cite{MilgExt}) to be incorrect and should never be cited.
\end{itemize}

\subsection{Campbell, D’Aurizio and Sondow }

In a paper \cite{CaDaSo} devoted to the evaluation of ${}_{3}^{}{\moversetsp{}{\mundersetsp{}{F_{2}^{}}}}\left(-\frac{1}{2},\frac{1}{4},\frac{3}{4};\frac{1}{2},1;1\right)$, relevant to the properties of a parbelos (q.v.), Campbell, D'Aurizio and Sondow regretfully noted that the computer algebra programs of that day (Maple and Mathematica) were unable to evaluate that particular sum\footnote{That observation remains true up to the present day.}. 
After considering the likelihood, they found that ``Watson's  theorem cannot be directly applied" to the sum under consideration, but unfortunately they omitted the possibility that it was {\bf contiguous} to Watson's theorem, after which they proceeded to formulate five independent proofs of that evaluation. Here, I add to their collection, by showing that, for general values of $c$, with $a=1/2-s,~b=1/2+s$ and contiguity parameters $m=0,~n=2$, a more general sum can be found by utilizing Chu's identity for contiguous Watson cases (see Appendix \ref{sec:Contig}, Section \ref{sec:MzSf} and \cite{ChuW}):
\begin{align} \nonumber
{}_{3}^{}{\moversetsp{}{\mundersetsp{}{F_{2}^{}}}}&\left(c ,\frac{1}{2}-s ,\frac{1}{2}+s ;\frac{1}{2},2\,c +2;1\right)=W_{0,2}( 1/2 - s, 1/2 + s, c)\\
& = 
\frac{4\,\Gamma \left(\frac{5}{2}+c \right)^{2} }{\left(c +1\right) \left(2\,c +1\right) \left(2\,c +3\right)^{2}}\left(\frac{\left(1-2\,s +2\,c \right) \cos \left(\frac{\pi  \left(2\,s +1\right)}{4}\right)}{\Gamma \left(\frac{5}{4}-\frac{s}{2}+c \right) \Gamma \left(\frac{3}{4}+\frac{s}{2}+c \right)}+\frac{\sin \left(\frac{\pi  \left(2\,s +1\right)}{4}\right) \left(1+2\,s +2\,c \right)}{\Gamma \left(\frac{5}{4}+\frac{s}{2}+c \right) \Gamma \left(\frac{3}{4}-\frac{s}{2}+c \right)}\right)\,.
\label{CDSa}
\end{align}
After evaluating the limit $c\rightarrow -1/2$, \eqref{CDSa} identifies
\begin{equation}
{}_{3}^{}{\moversetsp{}{\mundersetsp{}{F_{2}^{}}}}\left(-\frac{1}{2},\frac{1}{2}-s ,\frac{1}{2}+s ;\frac{1}{2},1;1\right)
 = 
2\,s +\frac{2\,\cos \left(\pi \,s \right)}{\pi}\left(s\,\psi \left(\frac{1}{4}+\frac{s}{2}\right)-s\psi \left(\frac{3}{4}+\frac{s}{2}\right) +1\right) \,.
\label{CDSb}
\end{equation}
Setting $s=1/4$ and making use of \eqref{PpqE}, i.e.
\begin{equation}
\psi \left(\frac{7}{8}\right)-\psi \left(\frac{3}{8}\right) = 
\left(\pi -2\,\ln \left(1+\sqrt{2}\right)\right) \sqrt{2}\,,
\label{psi3m7}
\end{equation}
\eqref{CDSb} yields
\begin{equation}
{}_{3}^{}{\moversetsp{}{\mundersetsp{}{F_{2}^{}}}}\left(-\frac{1}{2},\frac{1}{4},\frac{3}{4};\frac{1}{2},1;1\right)
 = \frac{\ln \left(1+\sqrt{2}\right)+\sqrt{2}}{\pi}\,,
\label{CDS}
\end{equation}
an ``exotic" identity in agreement with the several proofs presented by Campbell, D’Aurizio and Sondow.

\subsection{Xiaoxia Wang} \label{sec:XWang}

In 2013, Xiaoxia Wang \cite{WangX}, ``inspired by the work of Chu and Wang who established Whipple-like fomulas for the pattern of the terminating series"
\begin{equation}
{}_{3}^{}{\moversetsp{}{\mundersetsp{}{F_{2}^{}}}}\left(-m ,m +\epsilon +\delta ,x ;1+\epsilon ,\delta +2\,x ;1\right),\label{XWa} 
\end{equation}
{\bf where $\epsilon$ and $\delta$ are integers}, $|\delta|\leq 5$ and $0\leq\epsilon\leq 5$,  extended that work by employing the Gould–Hsu inversion
technique (q.v.). This required three pages of analysis to derive closed general expressions for a larger range of both $\epsilon$ and $\delta$ characterized by $\epsilon=2t$ and $\epsilon=2t+1$, {\bf where $t$ is also an integer} and $\epsilon + \delta=1$. This was later extended to $|\epsilon+\delta|\leq2$ with examples. For the case where $m$ is even, Wang finds
\begin{equation}
{}_{3}^{}{\moversetsp{}{\mundersetsp{}{F_{2}^{}}}}\left(x ,-m ,1+m ;1+2\,j ,1-2\,j +2\,x ;1\right)
 = 
\frac{\left(j -x \right) \left(\frac{1}{2}-j \right)_{\frac{m}{2}} \left(1-x +j \right)_{\frac{m}{2}}}{\left(2\,j -x \right) \left(1+j \right)_{\frac{m}{2}} \left(x +\frac{1}{2}-j \right)_{\frac{m}{2}}}
\label{Wang3P}
\end{equation}
where the variable $t$ has been replaced by $j$ for clarity.

At this point, we point out that the template sum \eqref{XWa} is itself a special case of a more general, but still-terminating sum, characterized by an extended set of parameters. It can be recognized as being contiguous to Whipple's theorem (see \eqref{Whipple}). In short, we have
\begin{equation}
{}_{3}^{}{\moversetsp{}{\mundersetsp{}{F_{2}^{}}}}\left(a ,x ,k-a  ;1+2\,t +n , 2\,x +p -2\,t;1\right)
 = 
\Omega_{k-1,n+p} \left(a , x ,  2\,x-2\,t +p \right)
\label{WxWpp}
\end{equation}
where, in terms of Wang's notation, we have generalized $-m:=a$, $\epsilon + \delta=k$, $\delta:=p-2t$ and $\epsilon:=n+2t$, introduced a new independent variable $k$ and no longer require that $m$, $\epsilon$, $\delta$ and $t$ be integers. The full evaluation of \eqref{WxWpp}, being rather lengthy, is relegated to Appendix \ref{sec:XWangFull}.
We now consider the case $n=0,p=1,k=1$ corresponding to a generalization of Wang's Theorem 3, using straightforward computer-aided substitution into \eqref{LT1}, to obtain
\begin{align} \nonumber
&{}_{3}^{}{\moversetsp{}{\mundersetsp{}{F_{2}^{}}}}\left(a ,x ,1-a ;1+2\,t ,1-2\,t +2\,x ;1\right)
 = 
\frac{\pi \,4^{-x}\,\Gamma \left(1-2\,t +2\,x \right) \Gamma \left(1+2\,t \right)}{2t-x } \\ \nonumber
&\times \left(-\frac{1}{\Gamma \left( x +\frac{1}{2}-t-\frac{a}{2}\right) \Gamma \left( \frac{a}{2}+x-t \right) \Gamma \left(\frac{a}{2}+t +\frac{1}{2}\right) \Gamma \left(1+t -\frac{a}{2}\right)}\right. \\ 
&\left.\hspace{30pt}+\frac{1}{\Gamma \left(\frac{a}{2}+t \right) \Gamma \left(\frac{1}{2}+t -\frac{a}{2}\right) \Gamma \left(1+x -\frac{a}{2}-t \right) \Gamma \left(-t +\frac{a}{2}+\frac{1}{2}+x \right)}\right)\,,
\label{T3}
\end{align}
where we emphasize now that neither $a$ nor $t$ need be integers, and note that \eqref{T3} is symmetric under the interchange $a\leftrightarrow 1-a$. Introducing the integer $q$, where $q=m/2$ in terms of the variables used by Wang, in the case that $a=-2q$, $t=j$ and $j-q\leq 0$, we find the shorter form
\begin{align} \nonumber
{}_{3}^{}{\moversetsp{}{\mundersetsp{}{F_{2}^{}}}}&\left(x ,-2\,q ,2\,q +1;1+2\,j ,1-2\,j +2\,x ;1\right)
 \\
 &= 
-\frac{\pi \,4^{-x}\,\Gamma \left(1-2\,j +2\,x \right) \Gamma \left(1+2\,j \right)}{\left(2\,j -x\right) \Gamma \left( x -j+\frac{1}{2}+q \right) \Gamma \left(x-j -q  \right) \Gamma \left(j +\frac{1}{2}-q \right) \Gamma \left(1+q +j \right)}\,.
\label{T3b}
\end{align}
 In the case that $a=-2q-1$ and again using $t=j$, we find, if $j-q\leq 0$,
\begin{align} \nonumber
{}_{3}^{}{\moversetsp{}{\mundersetsp{}{F_{2}^{}}}}&\left(x ,-2\,q -1,2\,q +2;1+2\,j ,1-2\,j +2\,x ;1\right)
 \\
&=\frac{\pi \,4^{-x}\,\Gamma \left(1-2\,j +2\,x \right) \Gamma \left(1+2\,j \right)}{\left(2\,j -x \right) \Gamma \left(j-q -\frac{1}{2} \right) \Gamma \left(1+q +j \right) \Gamma \left(\frac{3}{2}+x +q -j \right) \Gamma \left(x-j -q \right)}\,.
\label{T3a}
\end{align}
The former of these is equivalent to \eqref{Wang3P}, where we note that Wang's solution \eqref{Wang3P} also fails unless $m/2-|t|\geq0$ (in Wang's notation where $t$ is an integer).

We also consider \cite[Eq.(4)]{WangX} 
\begin{align} \nonumber
{}_{3}^{}{\moversetsp{}{\mundersetsp{}{F_{2}^{}}}}&\left(x ,-m ,1+m ;2+2\,j ,2\,x-2\,j  ;1\right)\\
 &= 
\frac{\left(1+2\,j \right) \Gamma \left(\frac{1}{2}-j +\frac{m}{2}\right) \Gamma \left(1-x +j +\frac{m}{2}\right) \Gamma \left(1+j \right) \Gamma \left(x +\frac{1}{2}-j \right)}{2 \left(1+2\,j -x \right) \Gamma \left(\frac{1}{2}-j \right) \Gamma \left(1-x +j \right) \Gamma \left(1+\frac{m}{2}+j \right) \Gamma \left(\frac{1}{2}+x +\frac{m}{2}-j \right)},
\label{Wang4}
\end{align}
where the original variable $t$ has been rewritten $t:=j$ for clarity. \eqref{Wang4} is valid for even values of $m$ and, also requires $j<m/2$. The case corresponding to odd values of $m$ is given in \cite{WangX}. For general values of $a,t$ and $x$, we set $n=1,p=0,k=1$ in \eqref{LT1} and find
\begin{align} \nonumber
&{}_{3}^{}{\moversetsp{}{\mundersetsp{}{F_{2}^{}}}}\left(a ,x ,1-a;2+2\,t ,2\,x -2\,t ;1\right)
 =\frac{\pi \,4^{-x}  \Gamma \left(2\,x -2\,t \right) \Gamma \left(2+2\,t \right)}{1+2\,t -x}\\ \nonumber
&\times
 \left(\frac{1}{\Gamma \left(\frac{a}{2}-t +x \right) \Gamma \left(1+t -\frac{a}{2}\right) \Gamma \left(t +\frac{a}{2}+\frac{1}{2}\right) \Gamma \left(\frac{1}{2}-t +x -\frac{a}{2}\right)}\right.\\ &\left.\hspace{40pt}-\frac{4}{\Gamma \left(\frac{a}{2}-t -\frac{1}{2}+x \right) \Gamma \left(x -\frac{a}{2}-t \right) \Gamma \left(\frac{3}{2}+t -\frac{a}{2}\right) \Gamma \left(t +\frac{a}{2}+1\right)}\right)
\label{T4}
\end{align}
reducing to a generalization of the two parts of \cite[Eq.(4)]{WangX}, if $a=m$ even or odd, and $t=j$.

\section{Summary}

It has been shown that mainstay identities, developed over many years and commonly employed in 3F2(1) analysis, can be derived in a much more simple way than previously known. When the associated techniques are then applied to the usual identities themselves rather than their derivation, it becomes possible to develop generalizations and extensions, again in a simple way. Although some headway was made on the Miller/Paris challenge problem itself, by relaxing the ground-rules, several useful identities were developed that appear to be new, and a simpler derivation was found for two recent additions (i.e. \eqref{Hnew3} and \eqref{Hnew}) to the panoply of tools commonly available for analysts' use. 

With respect to exotic evaluations in the sense of Krattenthaler and Rivoal, many examples were added to the list of two-part hypergeometric 3F2(1) transformations that do not satisfy the Thomae transformations, and it was noted that series reversal should be included in the set of those transformations when terminating series are involved. In consideration of exotic evaluations in the sense of Chu, Chen and others, it was shown that the exotic property resides in the fact that certain parameters assume rational values that can be mapped onto angles where the associated trigonometric functions reduce to simple irrational, but fundamental constants. Otherwise the examples they put forward can still be summed in terms of digamma functions by recognizing that many such cases are simply contiguous to at least one of DWW and are therefore summable. Whether non-reducible instances of $\psi(p/q)$ are acceptably ``exotic" for general values of $p/q$ lies in the eyes of the beholder.

In a final section, a selection of recent papers were revisited to demonstrate that derivations they contained could be simplified and generalized using the methods presented here.

\section{Acknowledgements}
I am grateful to Prof. Roy Hughes, (University of New South Wales) for bringing the paper \cite{CaDaSo} to my attention and obtaining for me a copy of the paper \cite{MirzSaf} .

This paper represents the sole effort of the author. It was performed without financial input from any granting agency.

\appendix
\section{Appendix: Dixon/Whipple/Watson ($m=n=0$) and their contiguity counterparts} \label{sec:Contig}

Here are the three relations contiguous to Watson/Whipple/Dixon respectively, as evaluated by Chu \cite{ChuW} and defined by non-zero values of integers $m$ and $n$:\newline

Dixon
\begin{equation}
X_{m ,n}\! \left(a , b , c\right) = 
{}_{3}^{}{\moversetsp{}{\mundersetsp{}{F_{2}^{}}}}\! \left(a ,b ,c ;1+a -b +m ,1+a -c +n ;1\right)
\label{Dixon}
\end{equation}

Whipple:
\begin{equation}
\Omega_{m ,n}\! \left(a , b , c\right) = 
{}_{3}^{}{\moversetsp{}{\mundersetsp{}{F_{2}^{}}}}\! \left(a ,b ,1-a +m ;c ,1+2\,b -c +n ;1\right)
\label{Whipple}
\end{equation}

Watson:
\begin{equation}
W_{m ,n}\! \left(a , b , c\right) = 
{}_{3}^{}{\moversetsp{}{\mundersetsp{}{F_{2}^{}}}}\! \left(a ,b ,c ;2\,c +n ,\frac{1}{2}+\frac{a}{2}+\frac{b}{2}+\frac{m}{2};1\right)
\label{Watson}
\end{equation}

\section{Appendix: The non-trivial Thomae relations} \label{sec:Thomae}
\subsection{Non-terminating cases} \label{sec:nonTerm}

The nine non-trivial Thomae transformations for non-terminating 3F2(1) can be easily obtained by noting a Theorem of Hardy, as presented by Krattenthaler and Rivoal \cite[Eq. (2.3)]{KrattRiv}:

\begin{theorem} \label{sec:T21}
Let
\begin{equation}
s = a+b+c-e-f\,.
\end{equation}
Then the function
\begin{align} 
&\frac{1}{\Gamma \left(s \right) \Gamma \left(2\,e \right) \Gamma \left(2\,f \right)}\;\;_{3}F_{2}  \left(\overset{ \displaystyle {2a-s,2b-s,2c-s\,}}{_ {\displaystyle {2e,2f}}}|1\right)
:=\frac{1}{\Gamma \left(-s\right) \Gamma \left(e \right) \Gamma \left(f \right)}\,_{3}F_{2}  \left(\overset{ \displaystyle {a,b,c\,}}{_ {\displaystyle {e,f}}}|1\right)
\label{HT}
\end{align}
is a symmetric function of permutations among $\{a,b,c,e,f\}$, where the symbol ``:=" signifies the left ordered symbolic multi-variable replacements: $a:=s,c:=-c/2, c:=c+a-2e-2f,b:=-b/2,b:=b+a-2e-2f,e:=e/2,f:=f/2$
 .
\end{theorem}
\begin{cor} 
\label{sec:CorTom}
The set of permutations among $\{a,b,c,e,f\}$ applied to \eqref{HT} is closed.
\end{cor}
Written in full, the Thomae transformations are:

\begin{align}\nonumber
&{}_{3}^{}{\moversetsp{}{\mundersetsp{}{F_{2}^{}}}}\! \left(a ,b ,c ;e ,f ;1\right)
\\ \nonumber
&=\frac{ \Gamma \! \left(-s \right) \Gamma \! \left(f \right) \Gamma \! \left(e \right)}{\Gamma \! \left(c \right) \Gamma \! \left(e+f -b  -c \right) \Gamma \! \left(e +f -a -c \right)}\\
&\hspace{30pt}\times{}_{3}^{}{\moversetsp{}{\mundersetsp{}{F_{2}^{}}}}\! \left(-s ,f -c ,e -c ;e +f-b  -c ,e +f -a -c ;1\right)\,;\label{Thom1}
\\
\nonumber
&=\frac{\Gamma \! \left(-s\right) \Gamma \! \left(f \right) \Gamma \! \left(e \right) }{\Gamma \! \left(b \right) \Gamma \! \left(e+f -b -c \right) \Gamma \! \left(e +f-b  -a \right)}\\
&\hspace{30pt}\times{}_{3}^{}{\moversetsp{}{\mundersetsp{}{F_{2}^{}}}}\! \left(-s,e -b ,f-b ;e+f -b  -c ,e+f -b  -a ;1\right)\,\label{Thom2};
\\
&=\frac{\Gamma \! \left(-s\right) \Gamma \! \left(e \right) }{\Gamma \! \left(e -a \right) \Gamma \! \left(e+f -b  -c \right)}\times{}_{3}^{}{\moversetsp{}{\mundersetsp{}{F_{2}^{}}}}\! \left(a ,f-b  ,f -c ;f ,e+f -b  -c ;1\right)\,\label{Thom3};
\\
&=\frac{\Gamma \! \left(-s\right) \Gamma \! \left(f \right)}{\Gamma \! \left(f -a \right) \Gamma \! \left(e+f -b  -c \right)}\times {}_{3}^{}{\moversetsp{}{\mundersetsp{}{F_{2}^{}}}}\! \left(a ,e -c ,e -b ;e ,e+f -b  -c ;1\right)\,\label{Thom4};
\\ \nonumber
&=\frac{\Gamma \! \left(-s \right) \Gamma \! \left(f \right) \Gamma \! \left(e \right) }{\Gamma \! \left(a \right) \Gamma \! \left(e +f -a -c \right) \Gamma \! \left(e+f -b  -a \right)}\\ 
&\hspace{30pt}\times{}_{3}^{}{\moversetsp{}{\mundersetsp{}{F_{2}^{}}}}\! \left(-s ,e -a ,f -a ;e +f-b  -a ,e +f -a -c ;1\right)\,;\label{Thom5}
 \\
&=\frac{\Gamma \! \left(-s \right) \Gamma \! \left(e \right) }{\Gamma \! \left(e -b \right) \Gamma \! \left(e +f -a -c \right)}\times{}_{3}^{}{\moversetsp{}{\mundersetsp{}{F_{2}^{}}}}\! \left(b ,f -c ,f -a ;f ,e +f -a -c ;1\right)\,;\label{Thom6}
\\
&=\frac{\Gamma \! \left(-s \right) \Gamma \! \left(f \right) }{\Gamma \! \left(f-b  \right) \Gamma \! \left(e +f -a -c \right)}\times{}_{3}^{}{\moversetsp{}{\mundersetsp{}{F_{2}^{}}}}\! \left(b ,e -c ,e -a ;e ,e +f -a -c ;1\right)\,;\label{Thom7}
\\
&=\frac{\Gamma \! \left(-s \right) \Gamma \! \left(e \right) }{\Gamma \! \left(e -c \right) \Gamma \! \left(e+f -a  -b \right)}\times{}_{3}^{}{\moversetsp{}{\mundersetsp{}{F_{2}^{}}}}\! \left(c ,f -a ,f-b  ;f ,e+f -b  -a ;1\right)\,;\label{Thom8}
\\
&=\frac{\Gamma \! \left(-s \right) \Gamma \! \left(f \right)}{\Gamma \! \left(f -c \right) \Gamma \! \left(e+f -b  -a \right)}\times {}_{3}^{}{\moversetsp{}{\mundersetsp{}{F_{2}^{}}}}\! \left(c ,e -a ,e -b ;e ,e +f-b  -a ;1\right)\,.
\label{Thom9}
\end{align}
\begin{rem} \label{sec:RemXfm}
These transformations arise from all possible non-trivial permutations of exchanges among variables in the left-hand side of \eqref{HT} between $\{e,f\}\leftrightarrow \{a,b,c\}$ individually (6 possibilities) and in pairs (3 possibilities). They are summarized in Table \ref{sec:TabThom} below and, in terms of a group theory permutation matrix, by Beyer et. al. \cite[Section III]{BeyLouStein1987}.

Note that the right-hand side of \eqref{Thom1} diverges unless one of $e-c=-n$ or $f-c=-n$ so that it terminates; this is the original  basis of Minton's identity \eqref{Minton} -- see Appendix (\ref{sec:MintsExTheorem})

\end{rem}
\begin{table}
\centering
\begin{tabular}{|c|c|}
\hline
\rule[-1.5ex]{0pt}{4.5ex} Eq'n no. & interchange \\ 
\hline 
\rule[-1ex]{0pt}{2.5ex} \eqref{Thom1} & $\{e,f\}\leftrightarrow\{a,c\}$ \\ 
\hline 
\rule[-1ex]{0pt}{2.5ex} \eqref{Thom2} & $\{e,f\}\leftrightarrow\{a,b\}$ \\ 
\hline 
\rule[-1ex]{0pt}{2.5ex} \eqref{Thom3} & $\{e\}\leftrightarrow\{a\}$ \\ 
\hline 
\rule[-1ex]{0pt}{2.5ex} \eqref{Thom4} & $\{f\}\leftrightarrow\{a\}$ \\ 
\hline 
\rule[-1ex]{0pt}{2.5ex} \eqref{Thom5} & $\{e,f\}\leftrightarrow\{b,c\}$ \\ 
\hline 
\rule[-1ex]{0pt}{2.5ex} \eqref{Thom6} & $\{e\}\leftrightarrow\{c\}$ \\ 
\hline 
\rule[-1ex]{0pt}{2.5ex} \eqref{Thom7} & $\{f\}\leftrightarrow\{c\}$ \\ 
\hline 
\rule[-1ex]{0pt}{2.5ex} \eqref{Thom8} & $\{e\}\leftrightarrow\{b\}$ \\ 
\hline 
\rule[-1ex]{0pt}{2.5ex}\eqref{Thom9} & $\{f\}\leftrightarrow\{b\}$ \\ 
\hline 
\end{tabular} 
\caption{Correspondence between each of the Thomae transformations and the interchange of variables on the right-hand side of \eqref{HT}.}
\label{sec:TabThom}
\end{table}
\subsection{Terminating cases} \label{sec:termCase}
%
For completeness, the 7 basic transformations of Rao, Van der Jeugt, Raynal, Jagannathan and Rajeswari (RJRJR) \cite{RaoJeu1992} are reproduced below:
\begin{align} \nonumber
{}_{3}^{}&{\moversetsp{}{\mundersetsp{}{F_{2}^{}}}}\left(-n,a ,b  ;d ,e ;1\right)\\
& = 
\frac{ \left(d -a\right)_ n }{{\left(d\right)_ n}}\,{}_{3}^{}{\moversetsp{}{\mundersetsp{}{F_{2}^{}}}}\left(-n ,a, e-b  ;e ,1+a -d -n ;1\right)\,,
\label{Ra1p}\\
& = 
\left(-1\right)^{n}\frac{ \left(1-\sigma \right)_{n}\,}{\left(d \right)_{n}}\,{}_{3}^{}{\moversetsp{}{\mundersetsp{}{F_{2}^{}}}}\left(-n ,e-b  ,e -a ;e ,\sigma -n ;1\right)\,,
\label{Ra2p}\\
 &= 
\frac{\left(d -a \right)_{n} \left(e -a \right)_{n}\,}{\left(d \right)_{n} \left(e \right)_{n}}{}_{3}^{}{\moversetsp{}{\mundersetsp{}{F_{2}^{}}}}\left(-n ,a ,1-\sigma ;1+a -d -n ,1-e +a -n ;1\right)\,,
\label{Ra3}\\
&= 
\frac{\left(d -a \right)_{n} \left(b \right)_{n}\,}{\left(d \right)_{n} \left(e \right)_{n}}\,{}_{3}^{}{\moversetsp{}{\mundersetsp{}{F_{2}^{}}}}\left(-n ,e-b  ,1-d -n ;1-b -n ,1+a -d -n ;1\right)\,,
\label{Ra4}\\
& = \left(-1\right)^{n} 
\frac{\left(1-\sigma \right)_{n} \left(b \right)_{n}}{\left(d \right)_{n} \left(e \right)_{n}}\,{}_{3}^{}{\moversetsp{}{\mundersetsp{}{F_{2}^{}}}}\left(-n ,e-b ,d -b ;1-b -n ,\sigma -n ;1\right)\,,
\label{Ra6}\\
& = \left(-1\right)^{n}
\frac{ \left(d -a \right)_{n} \left(d -b \right)_{n}}{\left(d \right)_{n} \left(e \right)_{n}}\,{}_{3}^{}{\moversetsp{}{\mundersetsp{}{F_{2}^{}}}}\left(-n ,1-\sigma ,1-d -n ;1+b -d -n ,1+a -d -n ;1\right)\,,
\label{Ra8}\\
& = \left(-1\right)^{n}
\frac{ \left(a \right)_{n} \left(b \right)_{n}}{\left(d \right)_{n} \left(e \right)_{n}}\,{}_{3}^{}{\moversetsp{}{\mundersetsp{}{F_{2}^{}}}}\left(-n ,1-d -n ,-e -n +1;1-b -n ,1-a -n ;1\right)\,,
\label{RaB}
\end{align}
where $\sigma\equiv -a - b + e + d + n$. The remaining transformation are obtained by elementary permutations within the elements of $\{a,b\}$ and $\{d,e\}$. \eqref{RaB} embodies the reversal of series.

\section{Appendix: A selection of contiguity relations} \label{sec:select}

This Appendix lists a selection of contiguity relations for 3F2(1). The first four are reproduced from Rainville's book \cite[page 82]{Rainville}:

\begin{equation}
{}_{3}^{}{\moversetsp{}{\mundersetsp{}{F_{2}^{}}}}\left(a ,b ,c ;e ,f ;1\right)
 = 
\frac{a}{a-b}\,{}_{3}^{}{\moversetsp{}{\mundersetsp{}{F_{2}^{}}}}\left(b ,c ,a +1;e ,f ;1\right)-\frac{b}{a -b}\,{}_{3}^{}{\moversetsp{}{\mundersetsp{}{F_{2}^{}}}}\left(a ,c ,b +1;e ,f ;1\right)
\label{RainvEq14a}
\end{equation}

\begin{equation}
{}_{3}^{}{\moversetsp{}{\mundersetsp{}{F_{2}^{}}}}\left(a ,b ,c ;e ,f ;1\right)
 = 
\frac{\left(1-e \right)}{a -e +1} {}_{3}^{}{\moversetsp{}{\mundersetsp{}{F_{2}^{}}}}\left(a ,b ,c ;f ,e -1;1\right)+\frac{a\,}{a -e +1}{}_{3}^{}{\moversetsp{}{\mundersetsp{}{F_{2}^{}}}}\left(b ,c ,a +1;e ,f ;1\right)
\label{RainvEq15a}
\end{equation}

\begin{equation}
{}_{3}^{}{\moversetsp{}{\mundersetsp{}{F_{2}^{}}}}\left(a ,b ,c ;e ,f ;1\right)
 = 
-\frac{U_{2}\,}s\,{}_{3}^{}{\moversetsp{}{\mundersetsp{}{F_{2}^{}}}}\left(a ,b ,c ;e ,f +1;1\right)-\frac{U_{1}}{s}\,{}_{3}^{}{\moversetsp{}{\mundersetsp{}{F_{2}^{}}}}\left(a ,b ,c ;f ,e +1;1\right)
\label{RainvEq19a}
\end{equation}

\begin{equation}
{}_{3}^{}{\moversetsp{}{\mundersetsp{}{F_{2}^{}}}}\left(a ,b ,c ;e ,f ;1\right)
 = 
-W_{1}\;{}_{3}^{}{\moversetsp{}{\mundersetsp{}{F_{2}^{}}}}\left(b ,c ,a +1;f ,e +1;1\right)-W_{2}\;{}_{3}^{}{\moversetsp{}{\mundersetsp{}{F_{2}^{}}}}\left(b ,c ,a +1;e ,f +1;1\right)
\label{RainvEq21a}
\end{equation}
where 
\begin{align}
&U_{1} = 
\frac{\left(a -e \right) \left(b -e \right) \left(c -e \right)}{e \left(f -e \right)}\,,\\
&U_{2} = 
\frac{\left(a -f \right) \left(b -f \right) \left(c -f \right)}{f \left(e -f \right)}\,,\\
&W_{1}\equiv\frac{\left(b -e \right) \left(c -e \right)}{e \left(f -e \right)}\,,\\
&W_{2}\equiv\frac{\left(b -f \right) \left(c -f \right)}{f \left(e -f \right)}
\end{align}
and $s\equiv a+b+c-e-f$.

From a long-lost source, valid when $z=1$:
\begin{equation}
{}_{3}^{}{\moversetsp{}{\mundersetsp{}{F_{2}^{}}}}\left(a ,b ,c ;e ,f ;z \right)
 = 
{}_{3}^{}{\moversetsp{}{\mundersetsp{}{F_{2}^{}}}}\left(a -1,b ,c ;e ,f ;z \right)+\frac{b\,c\,z}{e\,f}\,{}_{3}^{}{\moversetsp{}{\mundersetsp{}{F_{2}^{}}}}\left(a ,b +1,c +1;e +1,f +1;z \right)
\label{Unkn1}
\end{equation}
and
\begin{equation}
{}_{3}^{}{\moversetsp{}{\mundersetsp{}{F_{2}^{}}}}\left(a ,b ,c ;e ,f ;z \right)
 = 
{}_{3}^{}{\moversetsp{}{\mundersetsp{}{F_{2}^{}}}}\left(a ,b ,c ;e +1,f ;z \right)+\frac{a\,b\,c\,z}{e \left(e +1\right) f}\,{}_{3}^{}{\moversetsp{}{\mundersetsp{}{F_{2}^{}}}}\left(a +1,b +1,c +1;e +2,f +1;z \right)\,.
\label{Unkn2}
\end{equation}

From Karp and Prilepkina\cite[Theorem 3.4]{Karp&Prep2018}, 
\begin{align} \nonumber
&{}_{3}^{}{\moversetsp{}{\mundersetsp{}{F_{2}^{}}}}\left(a ,b ,c +n ;b +p ,c ;1\right)
 = 
\frac{\left(b +p -1\right) \left(p -a -1\right)}{\left(p -1\right) \left(b +p -a -1\right)} {}_{3}^{}{\moversetsp{}{\mundersetsp{}{F_{2}^{}}}}\left(a ,b ,c +n ;b +p -1,c ;1\right)\\
&+\frac{a\,b \left({}_{3}^{}{\moversetsp{}{\mundersetsp{}{F_{2}^{}}}}\left(a +1,b +1,c +n ;b +p ,c ;1\right)-\frac{\left(c +n \right) }{c}{}_{3}^{}{\moversetsp{}{\mundersetsp{}{F_{2}^{}}}}\left(a +1,b +1,c +n +1;b +p ,c +1;1\right)\right)}{\left(p -1\right) \left(b +p -a -1\right)}
\label{Th3p4}
\end{align}
From Krattenthaler and Rivoal \cite[Eqs. (1.1) and (1.2)]{KrattRiv}
\begin{align} \nonumber
{}_{3}^{}{\moversetsp{}{\mundersetsp{}{F_{2}^{}}}}\left(1+a ,b +1,c ;\right.&\left.a +2\,b +1,2\,a +b +1;1\right)\\
 &= 
\frac{2 \left(a +b \right) }{2\,a +2\,b -c}{}_{3}^{}{\moversetsp{}{\mundersetsp{}{F_{2}^{}}}}\left(a ,b ,c ;a +2\,b +1,2\,a +b +1;1\right)\,
\label{HidbKR}
\end{align}
and
\begin{align} \nonumber
{}_{3}^{}{\moversetsp{}{\mundersetsp{}{F_{2}^{}}}}\left(a ,b ,c ;a +1,c +\frac{a \left(a -c +1\right)}{b -1}+1;1\right)
 &= 
\frac{\left(a -b +2\right) \left(a^{2}-a\,c +b\,c +a -c \right)}{\left(a +1\right) \left(a^{2}-a\,b -a\,c +b\,c +2\,a -c \right)}\\
&\times {}_{3}^{}{\moversetsp{}{\mundersetsp{}{F_{2}^{}}}}\left(c ,a +1,b -1;a +2,c +\frac{a \left(a -c +1\right)}{b -1};1\right)\,.
\label{Eq1p2}
\end{align}
%
From Chu and Wang, \cite[Eqs. (1a) - (4a)]{ChuWang2007}:
\begin{align*} \nonumber
{}_{3}^{}{\moversetsp{}{\mundersetsp{}{F_{2}^{}}}}&\left(a ,b,c ;e ,f ;1\right)
 = 
\frac{\left(\left(a +1-f \right) e +c\,b \right)}{\left(a +1-f \right) e} {}_{3}^{}{\moversetsp{}{\mundersetsp{}{F_{2}^{}}}}\left(a +1,b,c ;e +1,f ;1\right) \nonumber
\\&\hspace{10pt}+\frac{\mathfrak{s}\, \left(a +1\right) c\,b}{e \left(a +1-f \right) \left(e +1\right) f}\,{}_{3}^{}{\moversetsp{}{\mundersetsp{}{F_{2}^{}}}}\left(a +2,b+1,c +1;2+e ,f +1;1\right); \numberthis
\label{Eq1A}\\\\
\nonumber
& = 
\frac{\left(1+c +b -f \right) \left(1-f \right) }{\left(1+c -f \right) \left(1+b -f \right)}{}_{3}^{}{\moversetsp{}{\mundersetsp{}{F_{2}^{}}}}\left(a -1,b,c ;e ,f -1;1\right) \\
&\hspace{10pt}+\frac{\mathfrak{s}\, c\,b}{\left(1+c -f \right) \left(f -b -1\right) e}\,{}_{3}^{}{\moversetsp{}{\mundersetsp{}{F_{2}^{}}}}\left(a ,b +1,c +1;e +1,f ;1\right); \numberthis
\label{Eq2A}\\\\
 \nonumber 
& = 
\frac{\left(e\,f -a \left(1+c +b \right)\right) }{e\,f}{}_{3}^{}{\moversetsp{}{\mundersetsp{}{F_{2}^{}}}}\left(a ,c +1,b +1;e +1,f +1;1\right)\\
&\hspace{10pt}-\frac{a\,\mathfrak{s} \left(c +1\right) \left(b +1\right) }{\left(e +1\right) \left(f +1\right) e\,f}\,{}_{3}^{}{\moversetsp{}{\mundersetsp{}{F_{2}^{}}}}\left(a +1,b+2 ,c +2;2+e ,2+f ;1\right)\,; \numberthis
\label{Eq3A}\\ 
\nonumber
&= 
\frac{\left(1-e \right) \left(1-f \right) }{\left(1+a -e \right) \left(a +1-f \right)}{}_{3}^{}{\moversetsp{}{\mundersetsp{}{F_{2}^{}}}}\left(a ,c -1,b -1;e -1,f -1;1\right)\\
&\hspace{10pt}+\frac{a\,\mathfrak{s}\,}{\left(1+a -e \right) \left(a +1-f \right)}{}_{3}^{}{\moversetsp{}{\mundersetsp{}{F_{2}^{}}}}\left(a +1,c ,b ;e ,f ;1\right)\,, \numberthis
\label{Eq4A}
\end{align*}
where $\mathfrak{s}=1 + a + c + b - e - f$. For an extensive collection of consequent contiguity identities refer to \cite{ChuWang2007} and Section \ref{sec:CWExotic}.

\section{Appendix: Proofs} \label{sec:Proofs}
\subsection{Proof of Minton (Eq. \eqref{Minton}) with Extensions}\label{sec:MintsExTheorem}

\begin{theorem} 
For all positive integers $m,n$, it is true  that 
\begin{align} \nonumber
{}_{3}^{}{\moversetsp{}{\mundersetsp{}{F_{2}^{}}}}&\left(-m ,a ,c +n ;c ,n +a -m ;1\right)
 = 
\left(-1\right)^{m}\,\frac{\Gamma \left(c \right) \Gamma \left(n +a -m \right) \Gamma \left(1+m \right)}{\Gamma \left(c +n \right) \Gamma \left(a \right)}\\ 
&+\frac{\Gamma \left(c \right) \Gamma \left(n +a -m \right) \Gamma \left(1+n \right) }{\Gamma \left(c +n \right) \Gamma \left(a -m -c \right)}\moverset{n}{\munderset{k =1+m}{\sum}}\frac{\Gamma \left(-c +a -m +k \right) \left(-1\right)^{k}}{\Gamma \left(1+n -k \right) k\,\Gamma \left(a +k \right) \Gamma \left(k-m  \right)}\,.
\label{M3C}
\end{align}

\proof
Throughout the following, we require that $a,c$ are not integers, and begin with \eqref{Thom1} after substitutions:
\begin{align} \nonumber
{}_{3}^{}{\moversetsp{}{\mundersetsp{}{F_{2}^{}}}}\left(-m ,a ,c ;e ,f ;1\right)
& =
\,\frac{\Gamma \left(e \right) \Gamma \left(f \right) \Gamma \left(e +f -a +m -c \right)}{\Gamma \left(c \right) \Gamma \left(e +f +m -c \right) \Gamma \left(e +f -a -c \right)}
\\
&\times {}_{3}^{}{\moversetsp{}{\mundersetsp{}{F_{2}^{}}}}\left(e -c ,f -c ,e +f -a +m -c ;e +f +m -c ,e +f -a -c ;1\right)
\label{M2b}
\end{align} 
Let $e=c-n$, $c:=c+n$ and $f=s + n + a - m$ to obtain
\begin{align}
{}_{3}^{}{\moversetsp{}{\mundersetsp{}{F_{2}^{}}}}\left(a ,-m ,c +n ;c ,s +n +a -m ;1\right)
& = 
\left(-1\right)^{m}\,\frac{ \sin \left(\pi \,s \right)\Gamma \left(c \right) \Gamma \left(s +n +a -m \right) \Gamma \left(1+n \right) }{\Gamma \left(c +n \right) \pi \,\Gamma \left(a -m +s -c \right)}\\
&\times\moverset{n}{\munderset{k =0}{\sum}}\frac{\Gamma \left(s -c +a -m +k \right) \Gamma \left(s +k \right) \Gamma \left(1-s +m -k \right)}{\Gamma \left(1+n -k \right) \Gamma \left(k +1\right) \Gamma \left(s +a +k \right)}\,.
\label{M2C}
\end{align}
Note that the upper limit of the sum has been reduced to $n$ because of the factor $\Gamma \left(1+n -k \right)$ in the denominator. Now we split the sum into two parts, yielding

\begin{align} \nonumber
{}_{3}^{}{\moversetsp{}{\mundersetsp{}{F_{2}^{}}}}&\left(a ,-m ,c +n ;c ,s +n +a -m ;1\right)
 = 
\frac{\Gamma \left(c \right) \Gamma \left(s +n +a -m \right) \left(-1\right)^{m}\,\Gamma \left(1-s +m \right)}{\Gamma \left(c +n \right) \Gamma \left(1-s \right) \Gamma \left(s +a \right)}\\
&+\frac{\Gamma \left(c \right) \Gamma \left(s +n +a -m \right) \Gamma \left(1+n \right) }{\Gamma \left(c +n \right) \Gamma \left(a -m +s -c \right)}\moverset{n}{\munderset{k =1}{\sum}}\frac{\Gamma \left(s -c +a -m +k \right) \Gamma \left(s +k \right) \left(-1\right)^{k}}{\Gamma \left(1+n -k \right) \Gamma \left(k +1\right) \Gamma \left(s +a +k \right) \Gamma \left(s +k-m  \right)}
\label{M3A}
\end{align}
and consider two cases. If $s=0$ and $m\geq n$, the sum in \eqref{M3A} vanishes identically because of the factor $1/\Gamma(k-m)$. If $s=0$ and $m<n$, then the sum in \eqref{M3A} only vanishes when $k<m<n$. Putting the two possibilities together gives
\eqref{M3C}.
Because a sum vanishes if a lower limit exceeds an upper limit, \eqref{M3C} covers both cases. 
\end{theorem}
\begin{cor} \label{sec:MintonExCor}

We now consider the case where $s=p$, a positive integer. By the same reasoning as above, we obtain
\begin{align} \nonumber
{}_{3}^{}{\moversetsp{}{\mundersetsp{}{F_{2}^{}}}}&\left(-m ,a ,c +n ;c ,p +n +a -m ;1\right)
 = 
\frac{\Gamma \left(c \right) \Gamma \left(p +n +a -m \right) \Gamma \left(1+n \right) }{\Gamma \left(c +n \right) \Gamma \left(a -m +p -c \right)}\\
&\times\moverset{n}{\munderset{k =  1-p +m \geq 0}{\sum}}\frac{\Gamma \left(p -c +a -m +k \right) \Gamma \left(p +k \right) \left(-1\right)^{k}}{\Gamma \left(1+n -k \right) \Gamma \left(k +1\right) \Gamma \left(p +a +k \right) \Gamma \left(p -m +k \right)}\,.
\label{M3E}
\end{align}
Notice that, due to a quirk of the notation, \eqref{M3E} does not reduce to \eqref{M3C} if $p=0$.
\end{cor}

\subsection{Proof of \eqref{Fxcb9}} \label{sec:Fxcb9Proof}
By equating the right-hand sides of \eqref{Excb4} and \eqref{HAll2Ha} in the limit $e\rightarrow 1$, we obtain
\begin{align}\nonumber
{}_{3}^{}{\moversetsp{}{\mundersetsp{}{F_{2}^{}}}}&\left(-m ,a ,b ;1,a +b-m  ;1\right)
 =\,-\frac{\left(-1\right)^{m}\,\Gamma \left(1+m \right) \Gamma \left( a +b-m \right) }{\Gamma \left(a \right) \Gamma \left(b \right)}\\ \nonumber
&\times\moverset{m}{\munderset{k =0}{\sum}}\frac{\left[2\,\psi \left(1-k +m \right)-\psi \left(1-a -m +k \right)-\psi \left(1-m -b +k \right)\right) \Gamma \left(1-b -a +k \right)}{\Gamma \left(1-m -b +k \right] \Gamma \left(1-a -m +k \right) \Gamma \left(1-k +m \right)^{2}\,\Gamma \left(k +1\right)}\\
&+\left(-1\right)^{m}\,\frac{\left[2\,\psi \left(1+m \right)-\psi \left(b -m \right)-\psi \left(a -m \right)\right] \Gamma \left(1-b -a \right) \Gamma \left( a +b -m\right)}{\Gamma \left(b -m \right) \Gamma \left(a -m \right) \Gamma \left(1-b \right) \Gamma \left(1+m \right) \Gamma \left(1-a \right)}\,,
\label{FxcBb1}
\end{align}
after an appropriate redefinition of variables to transform the result into canonical form. Since the left-hand side of \eqref{FxcBb1} is a known 1-balanced Saalsch{\"u}tzian sum (see \eqref{Saal} with $c=1$), we can solve \eqref{FxcBb1} and thereby evaluate the sum
\begin{align} \nonumber
\moverset{m}{\munderset{k =0}{\sum}}&\frac{\left(2\,\psi \left(1-k +m \right)-\psi \left(1-a -m +k \right)-\psi \left(1-m -b +k \right)\right) \Gamma \left(1-b -a +k \right)}{\Gamma \left(1-m -b +k \right) \Gamma \left(1-a -m +k \right) \Gamma \left(1-k +m \right)^{2}\,\Gamma \left(k +1\right)}\\ \nonumber
 &=\, 
-\frac{\sin \left(\pi \,b \right)\,\sin \left(\pi \,a \right) \Gamma \left(a \right)^{2}\Gamma \left(b \right)^{2}\,\Gamma \left(1-b -a \right) \left[\psi \left(a -m \right)-2\,\psi \left(1+m \right)+\psi \left(b -m \right)\right] }{\pi^{2}\,\Gamma \left(b -m \right) \Gamma \left(a -m \right) \Gamma \left(1+m \right)^{2}}\\
&-\frac{\sin \left(\pi  \left(a +b \right)\right) \Gamma \left(a \right)^{2}\, \Gamma \left(b \right)^{2}\Gamma \left(1-b -a \right)}{\pi \,\Gamma \left(b -m \right) \Gamma \left(a -m \right) \Gamma \left(1+m \right)^{2}}\,.
\label{Fxcb6}
\end{align}
Further, the basic entity embedded in the left-hand side of \eqref{Fxcb6} is also a specialization of a known, terminating,  $1-$balanced sum (see \eqref{Saal}), that is
\begin{align} \nonumber
&\frac{\Gamma \left(1-b -a \right) }{\Gamma \left(1-b-m  \right) \Gamma \left(1-a -m \right) \Gamma \left(m +1\right)^{2}}\;{}_{3}^{}{\moversetsp{}{\mundersetsp{}{F_{2}^{}}}}\left(-m ,-m ,1-a -b ;1-b-m  ,1-a -m ;1\right)\\ \nonumber
&=\moverset{m}{\munderset{k =0}{\sum}}\frac{\Gamma \left(1-b -a +k \right)}{\Gamma \left(1-m -b +k \right) \Gamma \left(1-a -m +k \right) \Gamma \left(1-k +m \right)^{2}\,\Gamma \left(k +1\right)}\\
& = 
\frac{\Gamma \left(1-b -a \right) \Gamma \left(m +1-b \right) \Gamma \left(-a +m +1\right)}{\Gamma \left(1-b \right)^{2}\,\Gamma \left(1-a \right)^{2}\,\Gamma \left(1+m \right)^{2}}\,,
\label{FF}
\end{align}
and therefore, by differentiating \eqref{FF} first with respect to $a$, then with respect to $b$, and adding the results to \eqref{Fxcb6} with appropriate factors, we arrive at \eqref{Fxcb9}.
\begin{rem}
Since a 2-balanced Saalsch{\"u}tzian sum is also known (see \eqref{XidbTwo}), it is possible  to repeat the above exercise using $n=2$ in \eqref{Excb4} and eventually arrive at
\begin{align} \nonumber
&\moverset{m}{\munderset{k =0}{\sum}}\frac{\left[\psi \left( k-a -b \right)-\psi \left(1-k +m \right)\right]\, \Gamma \left(k-a -b  \right)}{\Gamma \left(1-m -b +k \right) \Gamma \left(1-a -m +k \right) \Gamma \left(1-k +m \right)^{2}\,\Gamma \left(k +1\right)}\\ \nonumber
 &= 
\frac{ \Gamma \left(1-b -a \right) \Gamma \left(m-b \right) \Gamma \left(m-a  \right)}{\Gamma \left(1-b \right)^{2}\,\Gamma \left(1-a \right)^{2}\,\Gamma \left(1+m \right)^{2}}
\left(1-\frac{a\,b}{\left(a +b \right)^{2}}+\left(m -\frac{a\,b}{a +b}\right)\right. \\
&\left.\times\left[\psi \left(a +1-m \right)+\psi \left(b +1-m \right)+\psi \left(1-b -a \right)-\psi \left(a \right)-\psi \left(b \right)-\psi \left(1+m \right)\right]\right)\,.
\label{Excb9}
\end{align}
\end{rem}

\subsection{Proof of \eqref{Ca16e}} \label{sec:Ca16eProof}
Dissect the left-hand side of modified \eqref{Ca16d} into its even and odd components to yield
\begin{align} \nonumber
\moverset{\infty}{\munderset{n =1}{\sum}}\! \frac{2^{-n}}{\Gamma \! \left(\frac{3}{2}-n \right) \Gamma \! \left(\frac{n}{2}+1\right) \Gamma \! \left(\frac{n}{2}+x \right)}
& = 
\frac{8}{3\,\pi \,\Gamma \! \left(x -\frac{1}{2}\right)}\moverset{\infty}{\munderset{n =0}{\sum}}\! \frac{3\,\Gamma \left(n -\frac{3}{4}\right) \sqrt{2}\,\Gamma \left(n -\frac{1}{4}\right) \Gamma \left(x -\frac{1}{2}\right)}{32\,\sqrt{\pi}\,\Gamma \left(n +\frac{1}{2}\right) \Gamma \left(x -\frac{1}{2}+n \right)}\\
&+\frac{2\,\Gamma \! \left(x \right)}{\sqrt{\pi}\,\Gamma \! \left(-\frac{1}{4}+x \right) \Gamma \! \left(\frac{1}{4}+x \right)}-\frac{8}{3\,\pi \,\Gamma \! \left(x -\frac{1}{2}\right)}-\frac{2}{\sqrt{\pi}\,\Gamma \! \left(x \right)}\,,
\label{PdoA}
\end{align}
where the elementary identity
\begin{equation}
{}_{2}^{}{\moversetsp{}{\mundersetsp{}{F_{1}^{}}}}\! \left(-\frac{1}{4},\frac{1}{4};x ;1\right)
 = 
\frac{2^{2\,x -\frac{3}{2}}\,\Gamma \! \left(x \right)^{2}}{\sqrt{\pi}\,\Gamma \! \left(-\frac{1}{2}+2\,x \right)}
\label{Heq}
\end{equation}
has been employed. Differentiate \eqref{PdoA} with respect to $x$ and then set $x=1$ to obtain
\begin{align}
\pi^{\frac{3}{2}} \moverset{\infty}{\munderset{n =1}{\sum}}\! \frac{2^{-n}\,\psi \! \left(\frac{n}{2}+1\right)}{\Gamma \! \left(\frac{3}{2}-n \right) \Gamma \! \left(\frac{n}{2}+1\right)^{2}}
& = 
\frac{\sqrt{2}}{4}\moverset{\infty}{\munderset{n =0}{\sum}}\! \frac{\Gamma \left(n -\frac{3}{4}\right) \Gamma \left(n -\frac{1}{4}\right) \psi \left(n +\frac{1}{2}\right)}{\Gamma \left(n +\frac{1}{2}\right)^{2}}+\left(\frac{16}{3}-24\,\sqrt{2}\right) \ln \! \left(2\right)\\
&+\left(2\,\pi -4\,\sqrt{2}+\frac{8}{3}\right) \gamma +16\,\sqrt{2}\,.
\label{PdoAd1}
\end{align}
Some minor rearrangement and simplification will produce the result \eqref{Ca16e}.

\section{Appendix: Lengthy identities} \label{sec:Hcbm}
The following subsections list a number of lengthy identities for reference.
\subsection{With reference to Section \ref{sec:LimCase}, the special case of \eqref{Hnew2} when $c=b+m$ and $m\geq n$ is}
\begin{align} \nonumber
&{}_{3}^{}{\moversetsp{}{\mundersetsp{}{F_{2}^{}}}}\left(a ,b ,n ;a +n ,b +m ;1\right)
 =\\ \nonumber
 & 
\frac{\Gamma \left(1-a +b -n \right)  \sin \left(\pi \,a \right) \Gamma \left(b +m \right) \Gamma \left(a +n \right) }{\pi \,\Gamma \left(m \right) \Gamma \left(b \right)}\left[\psi \left(b +m -n \right)-\psi \left(a \right)+\pi \cot \left(\pi  \left(a -b \right)\right)\right]
\\ \nonumber
&\hspace{20pt}\times \moverset{n -1}{\munderset{k =0}{\sum}}\frac{\left(-1\right)^{k}\Gamma \left(m +k \right) \Gamma \left(1-a +k \right) }{\Gamma \left(k +1\right)^{2}\,\Gamma \left(b +m +1-a -n +k \right) \Gamma \left(n-k  \right)}
\\ \nonumber 
&+\frac{\Gamma \left(1-b\right) \Gamma \left(b +m \right) \Gamma \left(a +n \right) }{\Gamma \left(a -b +n \right) \Gamma \left(b +m -n \right) \Gamma \left(a \right) \Gamma \left(n \right)}\moverset{m -n -1}{\munderset{k =0}{\sum}}\frac{\Gamma \left(k +1\right) \Gamma \left(a -b -m +1+k +n \right) \Gamma \left(m -1-k \right)}{\Gamma \left(k +n +1\right) \Gamma \left(m -k -n \right) \Gamma \left(2-b -m +k +n \right)}
\\ \nonumber
&+\frac{\Gamma \left(b +m \right) \Gamma \left(a +n \right) \sin \left(\pi \,a \right) \Gamma \left(1-a +b -n \right)}{\pi \,\Gamma \left(m \right) \Gamma \left(b \right)} \\ \nonumber
&\hspace{20pt}\times\moverset{n -1}{\munderset{k =0}{\sum}}\frac{\left(-1\right)^{k} \Gamma \left(1-a +k \right) \Gamma \left(m +k \right)}{\Gamma \left(k +1\right)^{2}\,\Gamma \left(b +m +1-a -n +k \right) \Gamma \left(n-k \right)} \left[\psi \left(k +1\right)-\psi \left(b +m +1-a -n +k \right)\right]
\\ \nonumber
&+\frac{\left(-1\right)^{n} \Gamma \left(b +m \right) \Gamma \left(a +n \right) \Gamma \left(1-b \right)}{\Gamma \left(a -b +n \right) \Gamma \left(b +m -n \right) \Gamma \left(a \right) \Gamma \left(n \right)}\\
&\hspace{10pt}\times\moverset{M}{\munderset{k =0}{\sum}}\frac{\left(-1\right)^{k }  \Gamma \left(n +m -1-k \right) \Gamma \left(a -b -m +1+k \right)}{\Gamma \left(n-k  \right) \Gamma \left(k +1\right) \Gamma \left(m -k \right) \Gamma \left(2-b-m +k \right)}\left[\psi \left(a -b -m +1+k \right)-\psi \left(n-k  \right)\right]\,
\label{Hn4b}
\end{align}
where $M=\mathit{min}(n-1,m-1)$.
\subsection{An element of the identity \eqref{H4D} in full.}
\begin{align} \nonumber
&H\left(m,n,j,b\right)  
=\frac{  \left(-1\right)^{n} }{\Gamma \left(b \right) \Gamma \left(m \right)}\moverset{j -1}{\munderset{k =0}{\sum}}\frac{\left(\psi \left(b +m +1-j -n +k \right)-\psi \left(k +1\right)\right) \Gamma \left(m +k \right)}{\Gamma \left(b +m +1-j -n +k \right) \Gamma \left(n -k \right) \Gamma \left(k +1\right)^{2}\,\Gamma \left(j -k \right)}\\  \nonumber
&-\frac{1}{\Gamma \left(b \right) \Gamma \left(b +m -n \right) \Gamma \left(n \right) \Gamma \left(j \right)}  \moverset{m-n  -1}{\munderset{k =0}{\sum}}\frac{\Gamma \left(k +1\right) \Gamma \left(m -1-k \right) \Gamma \left(b +m -n -k-1 \right)}{\Gamma \left(b -j +m -n -k \right) \Gamma \left(n +k +1\right) \Gamma \left(m -n -k \right)}
\\ \nonumber
& -\frac{\left(-1\right)^{n}  \left(\psi \left(b +m -n \right)-\psi \left(j \right)\right)  }{ \Gamma \left(1 +b -j-n \right) \Gamma \left(b  +m-j \right) \Gamma \left(n \right) \Gamma \left(j \right)}
\\ \nonumber
&\hspace{30pt}\times \moverset{n -1}{\munderset{k =0}{\sum}}\frac{\left(-1\right)^{k} \,\Gamma \left(b +1-j -n +k \right) \Gamma \left(n +m -1-k \right)}{\Gamma \left(k +1\right) \Gamma \left(m -k \right) \Gamma \left(n -k \right) \Gamma \left(b +1+k -n \right)}\\ \nonumber
& +\frac{\left(-1\right)^{n  }  }{\Gamma \left(b \right)  \Gamma \left(b +m -n \right) \Gamma \left(n \right) \Gamma \left(j \right)}
\\ & 
\times \moverset{n -1}{\munderset{k =0}{\sum}}\frac{\left(-1\right)^{k } \left(\psi \left(n -k \right)-\psi \left(b+m -j  -k \right)\right) \Gamma \left(n +m -1-k \right) \Gamma \left(b +m -k-1 \right)}{\Gamma \left(b  +m -k -j\right) \Gamma \left(n -k \right) \Gamma \left(m -k \right) \Gamma \left(k +1\right)}
\label{H4i}
\end{align}
\subsection{With reference to Section \ref{sec:XWang}} \label{sec:XWangFull}

The following evaluate the Whipple sum sought in Section  \ref{sec:XWang} in full generality based on \cite{ChuW}:

Specifically, if $k\geq 1+n+p$ then
\begin{align} \nonumber
&{}_{3}^{}{\moversetsp{}{\mundersetsp{}{F_{2}^{}}}}\left(a ,x ,k-a ;1+2\,t +n , 2\,x +p-2\,t ;1\right)
 =\\ \nonumber
&\,-\frac{2^{2\,k-2\,a }\,\Gamma \left(1+a -k \right) \Gamma \left(n +p +1\right) \Gamma \left(k -n -p \right) \Gamma \left( 2\,x +p-2\,t \right) \Gamma \left(1+2\,t +n \right)}{4\,\Gamma \left(a \right)\Gamma \left(1-a +2\,t +n \right) \Gamma \left( 2\,x +p -a-2\,t \right) }\\
&\times\moverset{n +p}{\munderset{i =0}{\sum}}\frac{\left(-1\right)^{i}\,\left(2\,x-4\,t  -1-2\,n +2\,i \right) \Gamma \left(1+4\,t +n -2\,x -p -i \right) }{\Gamma \left(1+i \right) \Gamma \left(n +p -i +1\right) \Gamma \left(2+4\,t +2\,n -2\,x -i \right)}S_{1}(i)\,.
\label{GT1}
\end{align}
where
\begin{align} \nonumber
S_{1}\left(i\right)\equiv& 
\moverset{k -1-n -p}{\munderset{j =0}{\sum}}\frac{\left(-1\right)^{j} \Gamma \left(x-t +\frac{p}{2}-\frac{a}{2}+\frac{i}{2}+\frac{j}{2}\right) \Gamma \left(t-\frac{a}{2} +n +\frac{p}{2}-\frac{i}{2}+\frac{j}{2}+\frac{1}{2}\right)}{\Gamma \left(1+j \right) \Gamma \left(k -n -p -j \right) \Gamma \left(2\,x -k +1+n +p +j \right)  }\\
&\times\frac{\left(2\,x -2\,k +1+2\,n +2\,p +2\,j \right)\Gamma \left(2\,x -2\,k +1+2\,n +2\,p +j \right)  }{\Gamma \left(\frac{3}{2}+t +\frac{p}{2}+\frac{a}{2}-\frac{i}{2}+\frac{j}{2}-k +n \right)\Gamma \left(1+\frac{a}{2}-t +\frac{p}{2}+\frac{i}{2}+\frac{j}{2}+x -k \right)}\,.
\label{sGT1}
\end{align}

If $k<1+n+p$, then
\begin{align} \nonumber
&{}_{3}^{}{\moversetsp{}{\mundersetsp{}{F_{2}^{}}}}\left(a ,x ,k-a  ;1+2\,t +n ,2\,x +p-2\,t  ;1\right)=\,-\left(-1\right)^{n +p}2^{2\,p -2\,a +2\,n}\,\Gamma \left(n +p +1\right)
\\  \nonumber
&\times\frac{ \Gamma \left(2+n +p -k \right) \Gamma \left(2\,x-2\,t  +p \right)\Gamma \left(1-a \right)\Gamma \left(1+2\,t +n \right)\Gamma \left(x -k +1+n +p \right) }{\Gamma \left(2\,x +p-2\,t  -a \right) \Gamma \left(1-a +2\,t +n \right)\Gamma \left(x \right)  \Gamma(1-a+n+p)}\\
&\times\moverset{n +p}{\munderset{i =0}{\sum}}\frac{\left(-1\right)^{i} \left(2\,x-4\,t  -1-2\,n +2\,i \right) \Gamma \left(1+4\,t +n -2\,x -p -i \right) 
}{\Gamma \left(2+4\,t +2\,n -2\,x -i \right) \Gamma \left(1+i \right) \Gamma 
\left(n +p -i +1\right)} S_{2}(i)\,.
\label{LT1}
\end{align}
where
\begin{align} \nonumber
S_{2}(i)\equiv&\moverset{1 +n +p-k}{\munderset{j =0}{\sum}}\frac{ \Gamma \left(x-t +\frac{p}{2}-\frac{a}{2}+\frac{i}{2}+\frac{j}{2}\right) \Gamma \left(t+n-\frac{a}{2} +\frac{p}{2}-\frac{i}{2}+\frac{j}{2}+\frac{1}{2}\right) }{\Gamma \left(1+j \right) \Gamma \left(2+n +p -k -j \right) \Gamma \left(2\,x -k +1+n +p +j \right)  }\\
&\times\frac{(2\,x-1+2\,j)\Gamma \left(2\,x -1+j \right)}{\Gamma \left(\frac{a}{2}-t -n -\frac{p}{2}+\frac{i}{2}+\frac{j}{2}+x \right)\Gamma \left(\frac{1}{2}+t -\frac{p}{2}+\frac{a}{2}-\frac{i}{2}+\frac{j}{2}\right)}\,.
\label{sLT1}
\end{align}

\section{Appendix: Selected additions to the database} \label{sec:DataNew}

Listed below are a selection of new additions to the database originally described in \cite{Milgram447}.
\begin{align}
470:~~{}_{3}^{}{\moversetsp{}{\mundersetsp{}{F_{2}^{}}}}\! \left(a ,b ,b ;b +1,b -n +2;1\right)
 = 
\frac{\left(-1\right)^{n} b \pi  \csc \! \left(\pi  a \right) \Gamma \! \left(n -1\right) \Gamma \! \left(b -n +2\right)}{\Gamma \! \left(a \right) \Gamma \! \left(b -a +1\right)}
\label{data470}
\end{align}
\begin{align}
471:~~{}_{3}^{}{\moversetsp{}{\mundersetsp{}{F_{2}^{}}}}\! \left(a ,b ,2-n ;b +1,2-n +a ;1\right)
 = 
\frac{\Gamma \! \left(n -1\right) \Gamma \! \left(b +1\right) \Gamma \! \left(a -b \right) \Gamma \! \left(2-n +a \right)}{\Gamma \! \left(a \right) \Gamma \! \left(b -1+n \right) \Gamma \! \left(a -n -b +2\right)}
\label{data471}
\end{align}
\begin{align}
472:~~{}_{3}^{}{\moversetsp{}{\mundersetsp{}{F_{2}^{}}}}\! \left(a ,b ,b -1+n ;b +1,b +1;1\right)
 = 
\frac{\Gamma \! \left(n -1\right) \Gamma \! \left(b +1\right)^{2} \Gamma \! \left(-a +1\right)}{\Gamma \! \left(b -a +1\right) \Gamma \! \left(b -1+n \right)}
\label{data472}
\end{align}
\begin{align}
473:~~{}_{3}^{}{\moversetsp{}{\mundersetsp{}{F_{2}^{}}}}\! \left(a ,b ,c ;b -n ,c +1;1\right)
 = 
\frac{\Gamma \! \left(-a +1\right) \Gamma \! \left(c +1-b +n \right) \Gamma \! \left(-b +1\right) \Gamma \! \left(c +1\right)}{\Gamma \! \left(c +1-b \right) \Gamma \! \left(c +1-a \right) \Gamma \! \left(n +1-b \right)}
\label{data473}
\end{align}
\begin{align} \nonumber
474:~~{}_{3}^{}{\moversetsp{}{\mundersetsp{}{F_{2}^{}}}}\!& \left(a ,b ,-n ;c ,1+a +b -c -n ;1\right)
 \\&=
 \frac{\Gamma \! \left(1+b -c \right) \Gamma \! \left(1+a -c \right) \Gamma \! \left(1+a +b -c -n \right) \Gamma \! \left(-c -n +1\right)}{\Gamma \! \left(b -c -n +1\right) \Gamma \! \left(1+a -c -n \right) \Gamma \! \left(-c +b +1+a \right) \Gamma \! \left(-c +1\right)}
\label{data474}
\end{align}
\begin{align} \nonumber
482:~~{}_{3}^{}{\moversetsp{}{\mundersetsp{}{F_{2}^{}}}}\!& \left(a ,b ,n ;c ,b +1-m ;1\right)
 = 
\frac{\Gamma \! \left(m \right) \Gamma \! \left(c \right) \Gamma \! \left(b +1-m \right)}{\Gamma \! \left(b -m +1-n \right) \Gamma \! \left(c -a \right) \Gamma \! \left(b \right) \Gamma \! \left(n \right)}\\
&\times\moverset{m -1}{\munderset{j =0}{{\sum}}}\! \frac{\Gamma \left(b -n -m +1+j \right) \Gamma \left(-a +c -n -m +1+j \right) \Gamma \left(n +m -1-j \right)}{\Gamma \left(1+j \right) \Gamma \left(m -j \right) \Gamma \left(1+c -n -m +j \right)}
\label{data482}
\end{align}
\begin{align} \nonumber
483:~~{}_{3}^{}{\moversetsp{}{\mundersetsp{}{F_{2}^{}}}}&\left(a ,b ,c ;c +n ,1+c -m ;1\right)
\\& \nonumber = 
\frac{\Gamma \! \left(m \right) \Gamma \! \left(1+c -m \right)  \Gamma \! \left(n +1-a +c -b -m \right) \Gamma \! \left(c +n \right)}{\Gamma \! \left(1+c -b -m \right) \Gamma \! \left(1-a +c -m \right) \Gamma \! \left(c -a +n \right) \Gamma \! \left(c -b +n \right) \Gamma \! \left(n \right)}\\&
\times\moverset{m -1}{\munderset{j =0}{{\sum}}}\! \frac{\Gamma \left(1-b +c -m +j \right) \Gamma \left(1-a +c -m +j \right) \Gamma \left(n +m -1-j \right)}{\Gamma \left(1+j \right) \Gamma \left(m -j \right) \Gamma \left(1+c -m +j \right)}
\label{data483}
\end{align}
\begin{align} \nonumber
485:~~{}_{3}^{}{\moversetsp{}{\mundersetsp{}{F_{2}^{}}}}\!& \left(a ,b ,a -m +1-n ;c ,1+a -m ;1\right)\\ \nonumber &
 = 
\frac{\Gamma \! \left(m \right) \Gamma \! \left(c \right)  \Gamma \! \left(c-b  -a +n \right) \Gamma \! \left(1+a -m \right)}{\Gamma \! \left(a -m +1-n \right) \Gamma \! \left(c -a \right) \Gamma \! \left(m -1-a +c +n \right) \Gamma \! \left(a \right) \Gamma \! \left(n \right)}\\
&\times \moverset{m -1}{\munderset{j =0}{{\sum}}}\! \frac{\Gamma \left(a -m +1-n +j \right) \Gamma \left(c-a  +j \right) \Gamma \left(n +m -1-j \right)}{\Gamma \left(1+j \right) \Gamma \left(m -j \right) \Gamma \left(1-b +c -m +j \right)}
\label{data485}
\end{align}
\begin{align}
486:~~{}_{3}^{}{\moversetsp{}{\mundersetsp{}{F_{2}^{}}}}\! \left(b ,a +n ,c +m ;1+a ,c +1;1\right)
 = 
\left\{\begin{array}{cc}
0 & 0<2-b -m -n  
\\
 \infty  & \mathit{otherwise}  
\end{array}\right.
\label{data486}
\end{align}
\begin{equation}
487:~~{}_{3}^{}{\moversetsp{}{\mundersetsp{}{F_{2}^{}}}}\! \left(1,a ,b ;2,c ;1\right)
 = 
\frac{\Gamma \! \left(c \right) \Gamma \! \left(-b +c -a +1\right)-\Gamma \! \left(c -a \right) \Gamma \! \left(c -b \right) \left(c -1\right)}{\left(a -1\right) \left(b -1\right) \Gamma \! \left(c -a \right) \Gamma \! \left(c -b \right)}
\label{data487}
\end{equation}
\begin{equation}
491:~~{}_{3}^{}{\moversetsp{}{\mundersetsp{}{F_{2}^{}}}}\! \left(b ,-n ,a +n ;b +1,1+a ;1\right)
 = 
-\frac{\left(-1\right)^{n} \Gamma \! \left(b +1\right) \Gamma \! \left(b -a \right) \Gamma \! \left(1+a \right) \Gamma \! \left(n +1\right)}{\Gamma \! \left(a +n \right) \Gamma \! \left(b +n +1\right) \Gamma \! \left(1+b -a -n \right)}
\label{data491}
\end{equation}

\begin{align} \nonumber
494:~~{}_{3}^{}{\moversetsp{}{\mundersetsp{}{F_{2}^{}}}}&\left(a ,1,1;1+a ,a +n ;1\right)
 = \,-\frac{ a\,\Gamma \left(a +n \right) \Gamma \left(a +n-1 \right)}{2\,\Gamma \left(a \right)^{2}}\\&\times
\moverset{n -1}{\munderset{k =0}{\sum}}\frac{\left(-1\right)^{k} \left(2\,a+2\,k -1 \right) \left(\psi \left(1, \frac{1+a+k}{2}\right)-\psi \left(1, \frac{a+k}{2}\right)\right) \Gamma \left(2\,a +k -1\right)}{\Gamma \left(n -k \right) \Gamma \left(2\,a +n +k -1\right) \Gamma \left(1+k \right)}
\label{data494}
\end{align}

\end{document}